\documentclass{gtart}

\def\ifplaintex{\expandafter\ifx\csname documentclass\endcsname\relax}

\def\gtp{{\mathsurround=0pt\it $\cal G\mskip-2mu$eometry \&\ 
$\cal T\!\!$opology $\cal P\!$ublications}}  

\def\recd{{\small Received:\qua\receiveddate\ifx\reviseddate\relax
\else\qquad Revised:\qua\reviseddate\fi\par}} 

\def\lognumber#1{\def\thelognumber{#1}}
\def\volumenumber#1{\def\thevolumenumber{#1}}
\def\volumeyear#1{\def\thevolumeyear{#1}}
\def\papernumber#1{\def\thepapernumber{#1}}
\def\pagenumbers#1#2{\def\startpage{#1}\def\finishpage{#2}}
\def\published#1{\def\publishdate{#1}}

\def\received#1{\def\receiveddate{#1}}

\def\accepted#1{\def\accepteddate{#1}}

\def\asciiauthors#1{\def\theasciiauthors{#1}}
\def\asciiaddress#1{\def\theasciiaddress{#1}}

\def\coverauthors#1{\def\thecoverauthors{#1}}
\long\def\asciiabstract#1{\long\def\theasciiabstract{#1}}

\let\\\par\let\thelognumber\relax\let\thevolumenumber\relax
\let\thepapernumber\relax\let\thevolumeyear\relax\let\startpage\relax
\let\finishpage\relax\let\publishdate\relax\let\receiveddate\relax
\let\reviseddate\relax\let\accepteddate\relax\let\theasciititle\relax
\let\theasciiauthors\relax\let\theasciiaddress\relax
\let\theasciiabstract\relax

\let\thecoverauthors\relax\let\theasciiemail\relax

\ifplaintex
\font\logobig=cmssbx10 scaled 3836
\font\logomed=cmssbx10 scaled 2557
\else
\font\logobig=cmssbx10 scaled 4200
\font\logomed=cmssbx10 scaled 2800
\fi

\long\def\makeagttitle{   
\count0=\startpage
\agt\hfill      
\hbox to 45truept{\vbox to 0pt{\vglue -13truept{\logomed A\kern -.37em{\logobig 
T}\kern -.38em G}\vss}\hss}
\break
{\small Volume \thevolumenumber\ (\thevolumeyear)
\startpage--\finishpage\nl
Published: \publishdate}

\vglue .25truein

{\parskip=0pt\leftskip 0pt plus
1fil\def\\{\par\smallskip}{\Large\bf\thetitle}\par\medskip} \vglue
0.05truein

{\parskip=0pt\leftskip 0pt plus 1fil\def\\{\par}{\sc\theauthors}
\par\medskip}%
 
\vglue 0.03truein 

{\small\leftskip 25truept\rightskip 25truept{\bf Abstract}\stdspace\theabstract

{\bf AMS Classification}\stdspace\theprimaryclass
\ifx\thesecondaryclass\relax\else; \thesecondaryclass\fi\par
{\bf Keywords}\stdspace \thekeywords\par}\vglue 7truept

}   

\ifplaintex
\font\phead=cmsl9 scaled 950
\font\pnum=cmbx10 scaled 913
\font\pfoot=cmsl9 scaled 950
\headline{\vbox to 0pt{\vskip -4.5mm\line{\small\phead\ifnum
\count0=\startpage ISSN 1472-2739 (on-line) 1472-2747 (printed)
\hfill {\pnum\folio}\else\ifodd\count0\def\\{ }%
\ifx\theshorttitle\relax\thetitle\else\theshorttitle\fi\hfill{\pnum\folio}
\else\def\\{ and }{\pnum\folio}\hfill\ifx\theshortauthors\relax\theauthors
\else\theshortauthors\fi\fi\fi}\vss}}
\footline{\vbox to 0pt{\vglue 0mm\line{\small\pfoot\ifnum\count0=\startpage
\copyright\ \gtp\hfill\else
\agt, Volume \thevolumenumber\ (\thevolumeyear)\hfill\fi}\vss}}
\else
\font=cmsl9 scaled 1050
\font\lnum=cmbx10 
\font=cmsl9 scaled 1050
\makeatletter
\def\@oddhead{{\small\lhead\ifnum\count0=\startpage ISSN 1472-2739 
(on-line) 1472-2747 (printed)\hfill {\lnum\number\count0}\else\ifodd\count0
\def\\{ }\ifx\theshorttitle\relax \thetitle \else\theshorttitle\fi\hfill
{\lnum\number\count0}\else\def\\{ and }{\lnum\number\count0}
\hfill\ifx\theshortauthors\relax 
\theauthors\else\theshortauthors\fi\fi\fi}}\def\@evenhead{\@oddhead}
\def\@oddfoot{\small\lfoot\ifnum\count0=\startpage\copyright\ \gtp\hfill\else
\agt, Volume \thevolumenumber\ (\thevolumeyear)\hfill\fi}
\def\@evenfoot{\@oddfoot}
\makeatother
\fi
\let\maketitlepage\makeagttitle
\let\makeshorttitle\maketitlepage
\let\maketitle\maketitlepage

\newwrite\gtoutfile
\long\gdef\makeheadfile{  
{\def\\{, }\def\s{ }
\immediate\openout\gtoutfile head.xxx
\immediate\write\gtoutfile{To: math@arxiv.org}
\immediate\write\gtoutfile{Subject: put OR rep NNNNN:ppppp}
\immediate\write\gtoutfile{--text follows this line--}
\immediate\write\gtoutfile{Proxy-for: \ifx\theasciiauthors\relax
\theauthors\else\theasciiauthors\fi\s<\ifx\theasciiemail\relax\theemail\else\theasciiemail\fi>}
\immediate\write\gtoutfile{\noexpand\\}
\immediate\write\gtoutfile{Authors: \ifx\theasciiauthors\relax
\theauthors\else\theasciiauthors\fi}
{\def\\{ }\immediate\write\gtoutfile{Title: \ifx\theasciititle\relax
\thetitle\else\theasciititle\fi}}
\immediate\write\gtoutfile{Subj-class: GT or SG, GR etc}
\immediate\write\gtoutfile{MSC-class: \theprimaryclass\ifx\thesecondaryclass\relax\else, \thesecondaryclass\fi}
\immediate\write\gtoutfile{Journal-ref: Algebr. Geom. Topol. \thevolumenumber\s
(\thevolumeyear) \startpage-\finishpage}
\immediate\write\gtoutfile{Comments: Published by Algebraic and
Geometric Topology at}
\immediate\write\gtoutfile{\s\s\s  http://www.maths.warwick.ac.uk/agt/AGTVol\thevolumenumber/agt-\thevolumenumber-\thepapernumber.abs.html}
\immediate\write\gtoutfile{\noexpand\\}
\immediate\write\gtoutfile{}
\ifx\theasciiabstract\relax
\immediate\write\gtoutfile{\theabstract}\else
\immediate\write\gtoutfile{\theasciiabstract}\fi
\immediate\write\gtoutfile{}
\immediate\write\gtoutfile{\noexpand\\}
\immediate\write\gtoutfile{}
\immediate\closeout\gtoutfile}}  

\def\maketitlepage{\makeagttitle\makeheadfile}
\let\makeshorttitle\maketitlepage
\let\maketitle\maketitlepage

\def\ifplaintex{\expandafter\ifx\csname documentclass\endcsname\relax}

\def\gtp{{\mathsurround=0pt\it $\cal G\mskip-2mu$eometry \&\ 
$\cal T\!\!$opology $\cal P\!$ublications}}  

\def\recd{{\small Received:\qua\receiveddate\ifx\reviseddate\relax
\else\qquad Revised:\qua\reviseddate\fi\par}} 

\def\lognumber#1{\def\thelognumber{#1}}
\def\volumenumber#1{\def\thevolumenumber{#1}}
\def\volumeyear#1{\def\thevolumeyear{#1}}
\def\papernumber#1{\def\thepapernumber{#1}}
\def\pagenumbers#1#2{\def\startpage{#1}\def\finishpage{#2}}
\def\published#1{\def\publishdate{#1}}

\def\received#1{\def\receiveddate{#1}}

\def\accepted#1{\def\accepteddate{#1}}

\def\asciiauthors#1{\def\theasciiauthors{#1}}
\def\asciiaddress#1{\def\theasciiaddress{#1}}

\def\coverauthors#1{\def\thecoverauthors{#1}}
\long\def\asciiabstract#1{\long\def\theasciiabstract{#1}}

\let\\\par\let\thelognumber\relax\let\thevolumenumber\relax
\let\thepapernumber\relax\let\thevolumeyear\relax\let\startpage\relax
\let\finishpage\relax\let\publishdate\relax\let\receiveddate\relax
\let\reviseddate\relax\let\accepteddate\relax\let\theasciititle\relax
\let\theasciiauthors\relax\let\theasciiaddress\relax
\let\theasciiabstract\relax

\let\thecoverauthors\relax\let\theasciiemail\relax

\ifplaintex
\font\logobig=cmssbx10 scaled 3836
\font\logomed=cmssbx10 scaled 2557
\else
\font\logobig=cmssbx10 scaled 4200
\font\logomed=cmssbx10 scaled 2800
\fi

\long\def\makeagttitle{   
\count0=\startpage
\agt\hfill      
\hbox to 45truept{\vbox to 0pt{\vglue -13truept{\logomed A\kern -.37em{\logobig 
T}\kern -.38em G}\vss}\hss}
\break
{\small Volume \thevolumenumber\ (\thevolumeyear)
\startpage--\finishpage\nl
Published: \publishdate}

\vglue .25truein

{\parskip=0pt\leftskip 0pt plus
1fil\def\\{\par\smallskip}{\Large\bf\thetitle}\par\medskip} \vglue
0.05truein

{\parskip=0pt\leftskip 0pt plus 1fil\def\\{\par}{\sc\theauthors}
\par\medskip}%
 
\vglue 0.03truein 

{\small\leftskip 25truept\rightskip 25truept{\bf Abstract}\stdspace\theabstract

{\bf AMS Classification}\stdspace\theprimaryclass
\ifx\thesecondaryclass\relax\else; \thesecondaryclass\fi\par
{\bf Keywords}\stdspace \thekeywords\par}\vglue 7truept

}   

\ifplaintex
\font\phead=cmsl9 scaled 950
\font\pnum=cmbx10 scaled 913
\font\pfoot=cmsl9 scaled 950
\headline{\vbox to 0pt{\vskip -4.5mm\line{\small\phead\ifnum
\count0=\startpage ISSN 1472-2739 (on-line) 1472-2747 (printed)
\hfill {\pnum\folio}\else\ifodd\count0\def\\{ }%
\ifx\theshorttitle\relax\thetitle\else\theshorttitle\fi\hfill{\pnum\folio}
\else\def\\{ and }{\pnum\folio}\hfill\ifx\theshortauthors\relax\theauthors
\else\theshortauthors\fi\fi\fi}\vss}}
\footline{\vbox to 0pt{\vglue 0mm\line{\small\pfoot\ifnum\count0=\startpage
\copyright\ \gtp\hfill\else
\agt, Volume \thevolumenumber\ (\thevolumeyear)\hfill\fi}\vss}}
\else
\font=cmsl9 scaled 1050
\font\lnum=cmbx10 
\font=cmsl9 scaled 1050
\makeatletter
\def\@oddhead{{\small\lhead\ifnum\count0=\startpage ISSN 1472-2739 
(on-line) 1472-2747 (printed)\hfill {\lnum\number\count0}\else\ifodd\count0
\def\\{ }\ifx\theshorttitle\relax \thetitle \else\theshorttitle\fi\hfill
{\lnum\number\count0}\else\def\\{ and }{\lnum\number\count0}
\hfill\ifx\theshortauthors\relax 
\theauthors\else\theshortauthors\fi\fi\fi}}\def\@evenhead{\@oddhead}
\def\@oddfoot{\small\lfoot\ifnum\count0=\startpage\copyright\ \gtp\hfill\else
\agt, Volume \thevolumenumber\ (\thevolumeyear)\hfill\fi}
\def\@evenfoot{\@oddfoot}
\makeatother
\fi
\let\maketitlepage\makeagttitle
\let\makeshorttitle\maketitlepage
\let\maketitle\maketitlepage

\newwrite\gtoutfile
\long\gdef\makeheadfile{  
{\def\\{, }\def\s{ }
\immediate\openout\gtoutfile head.xxx
\immediate\write\gtoutfile{To: math@arxiv.org}
\immediate\write\gtoutfile{Subject: put OR rep NNNNN:ppppp}
\immediate\write\gtoutfile{--text follows this line--}
\immediate\write\gtoutfile{Proxy-for: \ifx\theasciiauthors\relax
\theauthors\else\theasciiauthors\fi\s<\ifx\theasciiemail\relax\theemail\else\theasciiemail\fi>}
\immediate\write\gtoutfile{\noexpand\\}
\immediate\write\gtoutfile{Authors: \ifx\theasciiauthors\relax
\theauthors\else\theasciiauthors\fi}
{\def\\{ }\immediate\write\gtoutfile{Title: \ifx\theasciititle\relax
\thetitle\else\theasciititle\fi}}
\immediate\write\gtoutfile{Subj-class: GT or SG, GR etc}
\immediate\write\gtoutfile{MSC-class: \theprimaryclass\ifx\thesecondaryclass\relax\else, \thesecondaryclass\fi}
\immediate\write\gtoutfile{Journal-ref: Algebr. Geom. Topol. \thevolumenumber\s
(\thevolumeyear) \startpage-\finishpage}
\immediate\write\gtoutfile{Comments: Published by Algebraic and
Geometric Topology at}
\immediate\write\gtoutfile{\s\s\s  http://www.maths.warwick.ac.uk/agt/AGTVol\thevolumenumber/agt-\thevolumenumber-\thepapernumber.abs.html}
\immediate\write\gtoutfile{\noexpand\\}
\immediate\write\gtoutfile{}
\ifx\theasciiabstract\relax
\immediate\write\gtoutfile{\theabstract}\else
\immediate\write\gtoutfile{\theasciiabstract}\fi
\immediate\write\gtoutfile{}
\immediate\write\gtoutfile{\noexpand\\}
\immediate\write\gtoutfile{}
\immediate\closeout\gtoutfile}}  

\def\maketitlepage{\makeagttitle\makeheadfile}
\let\makeshorttitle\maketitlepage
\let\maketitle\maketitlepage

\def\ifplaintex{\expandafter\ifx\csname documentclass\endcsname\relax}

\def\gtp{{\mathsurround=0pt\it $\cal G\mskip-2mu$eometry \&\ 
$\cal T\!\!$opology $\cal P\!$ublications}}  

\def\recd{{\small Received:\qua\receiveddate\ifx\reviseddate\relax
\else\qquad Revised:\qua\reviseddate\fi\par}} 

\def\lognumber#1{\def\thelognumber{#1}}
\def\volumenumber#1{\def\thevolumenumber{#1}}
\def\volumeyear#1{\def\thevolumeyear{#1}}
\def\papernumber#1{\def\thepapernumber{#1}}
\def\pagenumbers#1#2{\def\startpage{#1}\def\finishpage{#2}}
\def\published#1{\def\publishdate{#1}}

\def\received#1{\def\receiveddate{#1}}

\def\accepted#1{\def\accepteddate{#1}}

\def\asciiauthors#1{\def\theasciiauthors{#1}}
\def\asciiaddress#1{\def\theasciiaddress{#1}}

\def\coverauthors#1{\def\thecoverauthors{#1}}
\long\def\asciiabstract#1{\long\def\theasciiabstract{#1}}

\let\\\par\let\thelognumber\relax\let\thevolumenumber\relax
\let\thepapernumber\relax\let\thevolumeyear\relax\let\startpage\relax
\let\finishpage\relax\let\publishdate\relax\let\receiveddate\relax
\let\reviseddate\relax\let\accepteddate\relax\let\theasciititle\relax
\let\theasciiauthors\relax\let\theasciiaddress\relax
\let\theasciiabstract\relax

\let\thecoverauthors\relax\let\theasciiemail\relax

\ifplaintex
\font\logobig=cmssbx10 scaled 3836
\font\logomed=cmssbx10 scaled 2557
\else
\font\logobig=cmssbx10 scaled 4200
\font\logomed=cmssbx10 scaled 2800
\fi

\long\def\makeagttitle{   
\count0=\startpage
\agt\hfill      
\hbox to 45truept{\vbox to 0pt{\vglue -13truept{\logomed A\kern -.37em{\logobig 
T}\kern -.38em G}\vss}\hss}
\break
{\small Volume \thevolumenumber\ (\thevolumeyear)
\startpage--\finishpage\nl
Published: \publishdate}

\vglue .25truein

{\parskip=0pt\leftskip 0pt plus
1fil\def\\{\par\smallskip}{\Large\bf\thetitle}\par\medskip} \vglue
0.05truein

{\parskip=0pt\leftskip 0pt plus 1fil\def\\{\par}{\sc\theauthors}
\par\medskip}%
 
\vglue 0.03truein 

{\small\leftskip 25truept\rightskip 25truept{\bf Abstract}\stdspace\theabstract

{\bf AMS Classification}\stdspace\theprimaryclass
\ifx\thesecondaryclass\relax\else; \thesecondaryclass\fi\par
{\bf Keywords}\stdspace \thekeywords\par}\vglue 7truept

}   

\ifplaintex
\font\phead=cmsl9 scaled 950
\font\pnum=cmbx10 scaled 913
\font\pfoot=cmsl9 scaled 950
\headline{\vbox to 0pt{\vskip -4.5mm\line{\small\phead\ifnum
\count0=\startpage ISSN 1472-2739 (on-line) 1472-2747 (printed)
\hfill {\pnum\folio}\else\ifodd\count0\def\\{ }%
\ifx\theshorttitle\relax\thetitle\else\theshorttitle\fi\hfill{\pnum\folio}
\else\def\\{ and }{\pnum\folio}\hfill\ifx\theshortauthors\relax\theauthors
\else\theshortauthors\fi\fi\fi}\vss}}
\footline{\vbox to 0pt{\vglue 0mm\line{\small\pfoot\ifnum\count0=\startpage
\copyright\ \gtp\hfill\else
\agt, Volume \thevolumenumber\ (\thevolumeyear)\hfill\fi}\vss}}
\else
\font=cmsl9 scaled 1050
\font\lnum=cmbx10 
\font=cmsl9 scaled 1050
\makeatletter
\def\@oddhead{{\small\lhead\ifnum\count0=\startpage ISSN 1472-2739 
(on-line) 1472-2747 (printed)\hfill {\lnum\number\count0}\else\ifodd\count0
\def\\{ }\ifx\theshorttitle\relax \thetitle \else\theshorttitle\fi\hfill
{\lnum\number\count0}\else\def\\{ and }{\lnum\number\count0}
\hfill\ifx\theshortauthors\relax 
\theauthors\else\theshortauthors\fi\fi\fi}}\def\@evenhead{\@oddhead}
\def\@oddfoot{\small\lfoot\ifnum\count0=\startpage\copyright\ \gtp\hfill\else
\agt, Volume \thevolumenumber\ (\thevolumeyear)\hfill\fi}
\def\@evenfoot{\@oddfoot}
\makeatother
\fi
\let\maketitlepage\makeagttitle
\let\makeshorttitle\maketitlepage
\let\maketitle\maketitlepage

\newwrite\gtoutfile
\long\gdef\makeheadfile{  
{\def\\{, }\def\s{ }
\immediate\openout\gtoutfile head.xxx
\immediate\write\gtoutfile{To: math@arxiv.org}
\immediate\write\gtoutfile{Subject: put OR rep NNNNN:ppppp}
\immediate\write\gtoutfile{--text follows this line--}
\immediate\write\gtoutfile{Proxy-for: \ifx\theasciiauthors\relax
\theauthors\else\theasciiauthors\fi\s<\ifx\theasciiemail\relax\theemail\else\theasciiemail\fi>}
\immediate\write\gtoutfile{\noexpand\\}
\immediate\write\gtoutfile{Authors: \ifx\theasciiauthors\relax
\theauthors\else\theasciiauthors\fi}
{\def\\{ }\immediate\write\gtoutfile{Title: \ifx\theasciititle\relax
\thetitle\else\theasciititle\fi}}
\immediate\write\gtoutfile{Subj-class: GT or SG, GR etc}
\immediate\write\gtoutfile{MSC-class: \theprimaryclass\ifx\thesecondaryclass\relax\else, \thesecondaryclass\fi}
\immediate\write\gtoutfile{Journal-ref: Algebr. Geom. Topol. \thevolumenumber\s
(\thevolumeyear) \startpage-\finishpage}
\immediate\write\gtoutfile{Comments: Published by Algebraic and
Geometric Topology at}
\immediate\write\gtoutfile{\s\s\s  http://www.maths.warwick.ac.uk/agt/AGTVol\thevolumenumber/agt-\thevolumenumber-\thepapernumber.abs.html}
\immediate\write\gtoutfile{\noexpand\\}
\immediate\write\gtoutfile{}
\ifx\theasciiabstract\relax
\immediate\write\gtoutfile{\theabstract}\else
\immediate\write\gtoutfile{\theasciiabstract}\fi
\immediate\write\gtoutfile{}
\immediate\write\gtoutfile{\noexpand\\}
\immediate\write\gtoutfile{}
\immediate\closeout\gtoutfile}}  

\def\maketitlepage{\makeagttitle\makeheadfile}
\let\makeshorttitle\maketitlepage
\let\maketitle\maketitlepage

\def\ifplaintex{\expandafter\ifx\csname documentclass\endcsname\relax}

\def\gtp{{\mathsurround=0pt\it $\cal G\mskip-2mu$eometry \&\ 
$\cal T\!\!$opology $\cal P\!$ublications}}  

\def\recd{{\small Received:\qua\receiveddate\ifx\reviseddate\relax
\else\qquad Revised:\qua\reviseddate\fi\par}} 

\def\lognumber#1{\def\thelognumber{#1}}
\def\volumenumber#1{\def\thevolumenumber{#1}}
\def\volumeyear#1{\def\thevolumeyear{#1}}
\def\papernumber#1{\def\thepapernumber{#1}}
\def\pagenumbers#1#2{\def\startpage{#1}\def\finishpage{#2}}
\def\published#1{\def\publishdate{#1}}

\def\received#1{\def\receiveddate{#1}}

\def\accepted#1{\def\accepteddate{#1}}

\def\asciiauthors#1{\def\theasciiauthors{#1}}
\def\asciiaddress#1{\def\theasciiaddress{#1}}

\def\coverauthors#1{\def\thecoverauthors{#1}}
\long\def\asciiabstract#1{\long\def\theasciiabstract{#1}}

\let\\\par\let\thelognumber\relax\let\thevolumenumber\relax
\let\thepapernumber\relax\let\thevolumeyear\relax\let\startpage\relax
\let\finishpage\relax\let\publishdate\relax\let\receiveddate\relax
\let\reviseddate\relax\let\accepteddate\relax\let\theasciititle\relax
\let\theasciiauthors\relax\let\theasciiaddress\relax
\let\theasciiabstract\relax

\let\thecoverauthors\relax\let\theasciiemail\relax

\ifplaintex
\font\logobig=cmssbx10 scaled 3836
\font\logomed=cmssbx10 scaled 2557
\else
\font\logobig=cmssbx10 scaled 4200
\font\logomed=cmssbx10 scaled 2800
\fi

\long\def\makeagttitle{   
\count0=\startpage
\agt\hfill      
\hbox to 45truept{\vbox to 0pt{\vglue -13truept{\logomed A\kern -.37em{\logobig 
T}\kern -.38em G}\vss}\hss}
\break
{\small Volume \thevolumenumber\ (\thevolumeyear)
\startpage--\finishpage\nl
Published: \publishdate}

\vglue .25truein

{\parskip=0pt\leftskip 0pt plus
1fil\def\\{\par\smallskip}{\Large\bf\thetitle}\par\medskip} \vglue
0.05truein

{\parskip=0pt\leftskip 0pt plus 1fil\def\\{\par}{\sc\theauthors}
\par\medskip}%
 
\vglue 0.03truein 

{\small\leftskip 25truept\rightskip 25truept{\bf Abstract}\stdspace\theabstract

{\bf AMS Classification}\stdspace\theprimaryclass
\ifx\thesecondaryclass\relax\else; \thesecondaryclass\fi\par
{\bf Keywords}\stdspace \thekeywords\par}\vglue 7truept

}   

\ifplaintex
\font\phead=cmsl9 scaled 950
\font\pnum=cmbx10 scaled 913
\font\pfoot=cmsl9 scaled 950
\headline{\vbox to 0pt{\vskip -4.5mm\line{\small\phead\ifnum
\count0=\startpage ISSN 1472-2739 (on-line) 1472-2747 (printed)
\hfill {\pnum\folio}\else\ifodd\count0\def\\{ }%
\ifx\theshorttitle\relax\thetitle\else\theshorttitle\fi\hfill{\pnum\folio}
\else\def\\{ and }{\pnum\folio}\hfill\ifx\theshortauthors\relax\theauthors
\else\theshortauthors\fi\fi\fi}\vss}}
\footline{\vbox to 0pt{\vglue 0mm\line{\small\pfoot\ifnum\count0=\startpage
\copyright\ \gtp\hfill\else
\agt, Volume \thevolumenumber\ (\thevolumeyear)\hfill\fi}\vss}}
\else
\font=cmsl9 scaled 1050
\font\lnum=cmbx10 
\font=cmsl9 scaled 1050
\makeatletter
\def\@oddhead{{\small\lhead\ifnum\count0=\startpage ISSN 1472-2739 
(on-line) 1472-2747 (printed)\hfill {\lnum\number\count0}\else\ifodd\count0
\def\\{ }\ifx\theshorttitle\relax \thetitle \else\theshorttitle\fi\hfill
{\lnum\number\count0}\else\def\\{ and }{\lnum\number\count0}
\hfill\ifx\theshortauthors\relax 
\theauthors\else\theshortauthors\fi\fi\fi}}\def\@evenhead{\@oddhead}
\def\@oddfoot{\small\lfoot\ifnum\count0=\startpage\copyright\ \gtp\hfill\else
\agt, Volume \thevolumenumber\ (\thevolumeyear)\hfill\fi}
\def\@evenfoot{\@oddfoot}
\makeatother
\fi
\let\maketitlepage\makeagttitle
\let\makeshorttitle\maketitlepage
\let\maketitle\maketitlepage

\newwrite\gtoutfile
\long\gdef\makeheadfile{  
{\def\\{, }\def\s{ }
\immediate\openout\gtoutfile head.xxx
\immediate\write\gtoutfile{To: math@arxiv.org}
\immediate\write\gtoutfile{Subject: put OR rep NNNNN:ppppp}
\immediate\write\gtoutfile{--text follows this line--}
\immediate\write\gtoutfile{Proxy-for: \ifx\theasciiauthors\relax
\theauthors\else\theasciiauthors\fi\s<\ifx\theasciiemail\relax\theemail\else\theasciiemail\fi>}
\immediate\write\gtoutfile{\noexpand\\}
\immediate\write\gtoutfile{Authors: \ifx\theasciiauthors\relax
\theauthors\else\theasciiauthors\fi}
{\def\\{ }\immediate\write\gtoutfile{Title: \ifx\theasciititle\relax
\thetitle\else\theasciititle\fi}}
\immediate\write\gtoutfile{Subj-class: GT or SG, GR etc}
\immediate\write\gtoutfile{MSC-class: \theprimaryclass\ifx\thesecondaryclass\relax\else, \thesecondaryclass\fi}
\immediate\write\gtoutfile{Journal-ref: Algebr. Geom. Topol. \thevolumenumber\s
(\thevolumeyear) \startpage-\finishpage}
\immediate\write\gtoutfile{Comments: Published by Algebraic and
Geometric Topology at}
\immediate\write\gtoutfile{\s\s\s  http://www.maths.warwick.ac.uk/agt/AGTVol\thevolumenumber/agt-\thevolumenumber-\thepapernumber.abs.html}
\immediate\write\gtoutfile{\noexpand\\}
\immediate\write\gtoutfile{}
\ifx\theasciiabstract\relax
\immediate\write\gtoutfile{\theabstract}\else
\immediate\write\gtoutfile{\theasciiabstract}\fi
\immediate\write\gtoutfile{}
\immediate\write\gtoutfile{\noexpand\\}
\immediate\write\gtoutfile{}
\immediate\closeout\gtoutfile}}  

\def\maketitlepage{\makeagttitle\makeheadfile}
\let\makeshorttitle\maketitlepage
\let\maketitle\maketitlepage

\lognumber{35}
\volumenumber{2}
\volumeyear{2002}
\papernumber{35}
\pagenumbers{843}{895}
\received{15 September 2001}
\accepted{17 June 2002}
\published{15 October 2002}

\usepackage{amsmath,amssymb}

\numberwithin{equation}{section}
\newtheorem{theorem}[equation]{Theorem}
\newtheorem{lemma}[equation]{Lemma}
\newtheorem{proposition}[equation]{Proposition}
\newtheorem{corollary}[equation]{Corollary}
\newtheorem{completion}[equation]{Completion}
\newtheorem*{thm-A1}{Theorem A1}
\newtheorem*{thm-A2}{Theorem A2}
\newtheorem*{thm-A3}{Theorem A3}
\newtheorem*{thm-B1}{Theorem B1}
\newtheorem*{thm-B2}{Theorem B2}
\newtheorem*{thm-B3}{Theorem B3}
\newtheorem*{thm-C1}{Theorem C1}
\newtheorem*{thm-C2}{Theorem C2}
\newtheorem*{thm-C3}{Theorem C3}
\newtheorem*{com-D1}{Comment D1}
\newtheorem*{com-D2}{Comment D2}
\newtheorem*{exa-E1}{Example E1}
\newtheorem*{exa-E2}{Example E2}
\newtheorem*{exa-E3}{Example E3}
\newtheorem*{exa-E4}{Example E4}
\newtheorem*{exa-E5}{Example E5}
\newtheorem*{smith}{Smith Isomorphism Question}
\newtheorem*{laitinen}{Laitinen Conjecture}
\newtheorem*{basic-lemma}{Basic Lemma}
\newtheorem*{first}{First Rank Lemma}
\newtheorem*{second}{Second Rank Lemma}
\newtheorem*{subgroup}{Subgroup Lemma}
\newtheorem*{8-lemma}{8-condition Lemma}
\newtheorem*{8-corollary}{8-condition Corollary}
\newtheorem*{classification-thm}{Classification Theorem}
\newtheorem*{classification-cor}{Classification Corollary}
\newtheorem*{realization-thm}{Realization Theorem}
\newtheorem*{realization-cor}{Realization Corollary}
\newtheorem*{key}{Key Lemma}
\newtheorem*{ack}{Acknowledgements}

\begin{document}
\title{Smith equivalence and finite Oliver groups\\with Laitinen 
number 0 or 1}

\author{Krzysztof Pawa{\l}owski\\Ronald Solomon}
\asciiauthors{Krzysztof Pawalowski\\Ronald Solomon}
\coverauthors{Krzysztof Pawa{\noexpand\l}owski\\Ronald Solomon}
\email{kpa@main.amu.edu.pl, solomon@math.ohio-state.edu}

\address{Faculty of Mathematics and Computer Science,
Adam Mickiewicz University\\ul. Umultowska 87, 61-614 Pozna{\'n}, 
Poland\\\smallskip\\Department of Mathematics, The Ohio State 
University\\231 West 18th Avenue, Columbus, OH 43210--1174, USA}

\asciiaddress{Faculty of Mathematics and Computer Science,
Adam Mickiewicz University\\ul. Umultowska 87, 61-614 Poznan, 
Poland\\and\\Department of Mathematics, The Ohio State 
University\\231 West 18th Avenue, Columbus, OH 43210--1174, USA}

\begin{abstract}
In 1960, Paul A. Smith asked the following question. If a finite group
$G$ acts smoothly on a sphere with exactly two fixed points, is it
true that the tangent $G$-modules at the two points are always
isomorphic?  We focus on the case $G$ is an Oliver group and we
present a classification of finite Oliver groups $G$ with Laitinen
number $a_G = 0$ or $1$. Then we show that the Smith Isomorphism
Question has a negative answer and $a_G \geq 2$ for any finite Oliver
group $G$ of odd order, and for any finite Oliver group $G$ with a
cyclic quotient of order $pq$ for two distinct odd primes $p$
and~$q$. We also show that with just one unknown case, this question
has a negative answer for any finite nonsolvable gap group $G$ with
$a_G \geq 2$. Moreover, we deduce that for a finite nonabelian simple
group $G$, the answer to the Smith Isomorphism Question is affirmative
if and only if $a_G = 0$ or $1$.
\end{abstract}

\asciiabstract{
In 1960, Paul A. Smith asked the following question. If a finite group
G acts smoothly on a sphere with exactly two fixed points, is it true
that the tangent G-modules at the two points are always isomorphic?
We focus on the case G is an Oliver group and we present a
classification of finite Oliver groups G with Laitinen number a_G = 0
or 1. Then we show that the Smith Isomorphism Question has a negative
answer and a_G > 1 for any finite Oliver group G of odd order, and for
any finite Oliver group G with a cyclic quotient of order $pq$ for two
distinct odd primes p and q. We also show that with just one unknown
case, this question has a negative answer for any finite nonsolvable
gap group G with a_G > 1. Moreover, we deduce that for a finite
nonabelian simple group G, the answer to the Smith Isomorphism
Question is affirmative if and only if a_G = 0 or 1.}


\maketitle

\section*{0.1 \ \ The Smith Isomorphism Question}

Let $G$ be a finite group. By a real {\it $G$-module} we mean a finite
dimensional real vector space $V$ with a linear action of $G$.
Let $M$ be a smooth $G$-manifold with nonempty fixed point set $M^G$.
For any point $x \in M^G$, the tangent space $T_x(M)$ becomes a real
$G$-module by taking the derivatives (at the point $x$) of
the transformations $g: M \to M$, $z \mapsto gz$ for all $g \in G$.
We refer to this $G$-module $T_x(M)$ as to the
{\it tangent $G$-module at $x$}.

In 1960, Paul A.~Smith \cite[page~406]{Smith} asked the following question.

\begin{smith}
Is it true that for any smooth action of $G$ on a sphere with exactly two
fixed points, the tangent $G$-modules at the two points are isomorphic?
\end{smith}

Following \cite{Petrie:1}--\cite{Petrie:4}, two real $G$-modules $U$ and $V$
are called {\it Smith equivalent} if there exists a smooth action of $G$
on a sphere $S$ such that $S^G = \{x,y\}$ for two points $x$ and $y$
at which $T_x(S) \cong U$ and $T_y(S) \cong V$ as real $G$-modules.

In the real representation ring $RO(G)$ of $G$, we consider the subset
$Sm(G)$ consisting of the differences $U - V$ of real $G$-modules
$U$ and $V$ which are Smith equivalent. Choose a real $G$-module $W$
such that $\dim W^G = 1$. Set $S = S(W)$, the $G$-invariant unit sphere
in $W$. Then $S^G = \{x,y\}$ for the obvious points $x$ and $y$ in $S$,
and clearly as real $G$-modules, $T_x(S) \cong T_y(S) \cong W - W^G$,
the $G$-orthogonal complement of $W^G$ in $W$. As a result,
$$W - W = (W - W^G) - (W - W^G) \in Sm(G).$$
Therefore, $Sm(G)$ contains the trivial subgroup $0$ of $RO(G)$,
and the Smith Isomorphism Question can be restated as follows.
{\it Is it true that $Sm(G) = 0$}? As we shall see below, it may
happen that $Sm(G) \neq 0$, but in general, it is an open question
whether $Sm(G)$ is a subgroup of $RO(G)$.

In the following answers to the Smith Isomorphism Question,
$\mathbb{Z}_n$ is the cyclic group $\mathbb{Z}/n \mathbb{Z}$ of order $n$,
and $S_3$ is the symmetric group on three letters.

By \cite{Atiyah-Bott} and \cite{Milnor}, $Sm(\mathbb{Z}_p) = 0$ for
any prime $p$. According to \cite{Sanchez}, $Sm(\mathbb{Z}_{p^k}) = 0$
for any odd prime $p$ and any integer $k \geq 1$. By character theory,
$Sm(S_3) = 0$ and $Sm(\mathbb{Z}_n) = 0$ for $n = 2$, $4$, or $6$.
On the other hand, by \cite{Cappell-Shaneson:1}--\cite{Cappell-Shaneson:3},
$Sm(\mathbb{Z}_n) \neq 0$ for $n = 4q$ with $q \geq 2$. So, $G = \mathbb{Z}_8$
is the smallest group with $Sm(G) \neq 0$.

We refer the reader to \cite{Atiyah-Bott},
\cite{Cappell-Shaneson:1}--\cite{Cappell-Shaneson:3}, \cite{Cho},
\cite{Cho-Suh}, \cite{Dovermann-Petrie}, \cite{Dovermann-Petrie-Schultz},
\cite{Dovermann-Suh}, \cite{Dovermann-Washington}, \cite{Illman},
\cite{Laitinen-Pawalowski}, \cite{Masuda-Petrie}, \cite{Milnor},
\cite{Pawalowski:1}, \cite{Pawalowski:3}, \cite{Petrie:1}--\cite{Petrie:4},
\cite{Petrie-Randall}, \cite{Sanchez}, \cite{Schultz}, \cite{Suh}
for more related information.

If a finite group $G$ acts smoothly on a homotopy sphere $\varSigma$
with $\varSigma^G = \{x, y\}$, it follows from Smith theory that for
every $p$-subgroup $P$ of $G$ with $p \, \big| |G|$, the fixed point
set $\varSigma^P$ is either a connected manifold of dimension $\geq 1$,
or $\varSigma^P = \{x, y\}$.

Henceforth, we say that a smooth action of $G$ on a homotopy sphere
$\varSigma$ satisfies the {\it $8$-condition} if for every cyclic
$2$-subgroup $P$ of $G$ with $|P| \geq 8$, the fixed point set $\varSigma^P$
is connected (we recall that in \cite{Laitinen-Pawalowski}, such an action
of $G$ on $\varSigma$ is called $2$-proper). In particular, the action
of $G$ on $\varSigma$ satisfies the $8$-condition if $G$ has no element
of order $8$.

Now, two real $G$-modules $U$ and $V$ are called {\it Laitinen--Smith
equivalent} if there exists a smooth action of $G$ on a sphere $S$
satisfying the $8$-condition and such that $S^G = \{x, y\}$ for two
points $x$ and $y$ at which $T_x(S) \cong U$ and $T_y(S) \cong V$
as real $G$-modules.

Beside $Sm(G)$, we consider the subset $LSm(G)$ of $RO(G)$ consisting
of $0$ and the differences $U - V$ of real $G$-modules $U$ and $V$
which are Laitinen--Smith equivalent. Again, in general, if $LSm(G) \neq 0$,
it is an open question whether $LSm(G)$ is a subgroup of $RO(G)$.
Clearly, $LSm(G) \subseteq Sm(G)$.

If $G$ is a cyclic $2$-group with $|G| \geq 8$, then there are
no two real $G$-modules which are Laitinen--Smith equivalent.
Therefore $LSm(G) = 0$ while $Sm(G) \neq 0$ by
\cite{Cappell-Shaneson:1}--\cite{Cappell-Shaneson:3}.
In particular, $LSm(G) \neq Sm(G)$. However, if $G$ has no element
of order $8$, then $LSm(G) = Sm(G)$ (cf.~the 8-condition Lemma in Section~0.3).

Let $IO(G)$ be the intersection of the kernels ${\rm Ker}\,(RO(G) \to RO(P))$
of the restriction maps $RO(G) \to RO(P)$ taken for all subgroups $P$ of $G$
of prime power order. Set
$$IO(G,G) = IO(G) \cap\, {\rm Ker}\,(RO(G) \to \mathbb{Z})$$
where the map $RO(G) \to \mathbb{Z}$ is defined by
$U - V \mapsto \dim U^G - \dim V^G$.
In \cite{Laitinen-Pawalowski}, the abelian group
$IO(G,G)$ is denoted by $IO'(G)$.

According to \cite[Lemma~1.4]{Laitinen-Pawalowski}, the difference
$U - V$ of two Laitinen--Smith equivalent real $G$-modules $U$ and
$V$ belongs to $IO(G,G)$. Thus, the following lemma holds.

\begin{basic-lemma}
Let $G$ be a finite group. Then $LSm(G) \subseteq IO(G,G)$.
\end{basic-lemma}

Let $G$ be a finite group. Given two elements $g, h \in G$, $g$ is called
{\it real conjugate} to $h$ if $g$ or $g^{-1}$ is conjugate to $h$,
written $g \overset{\pm 1}{\sim} h$. Clearly,
$\overset{\pm 1}{\sim}$ is an equivalence relation in $G$.
For any $g \in G$, the resulting equivalence class $(g)^{\pm 1}$
is called the {\it real conjugacy class} of $g$. Note that
$(g)^{\pm 1} = (g) \cup (g^{-1})$,
the union of the conjugacy classes $(g)$ and $(g^{-1})$ of $g$ and $g^{-1}$,
respectively.

We denote by $a_G$ the number of real conjugacy classes $(g)^{\pm 1}$
of elements $g \in G$ not of prime power order. In 1996, Erkki Laitinen
has suggested to study the number $a_G$ while trying to answer
the Smith Isomorphism Question for specific finite groups $G$.
Henceforth, we refer to $a_G$ as to the {\it Laitinen number} of $G$.

The ranks of the free abelian groups $IO(G)$ and $IO(G,G)$ are computed
in \cite[Lemma~2.1]{Laitinen-Pawalowski} in terms of the Laitinen number $a_G$,
as follows.

\begin{first}
Let $G$ be a finite group. Then the following holds.
\begin{itemize}
\item[{\rm (1)}] ${\rm rk}\,IO(G) = a_G$. In particular, $IO(G) = 0$
      if and only if $a_G = 0$.
\item[{\rm (2)}] ${\rm rk}\,IO(G,G) = a_G - 1$ when $a_G \geq 1$, and
      ${\rm rk}\,IO(G,G) = 0$ when $a_G = 0$. In particular, $IO(G,G) = 0$
      if and only if $a_G = 0$ or $1$.
\end{itemize}
\end{first}

In 1996, Erkki Laitinen posed the following conjecture
(cf.~\cite[Appendix]{Laitinen-Pawalowski}).

\begin{laitinen}
Let $G$ be a finite Oliver group such that $a_G \geq 2$.
Then $LSm(G) \neq 0$.
\end{laitinen}

If $a_G = 0$ or $1$, $LSm(G) = 0$ by the Basic Lemma and the First Rank
Lemma. So, in the Laitinen Conjecture, the condition that $a_G \geq 2$
is necessary.

One may well conjecture that $Sm(G) \cap IO(G,G) \neq 0$ for any
finite Oliver group $G$ with $a_G \geq 2$. It is very likely that
$LSm(G) = Sm(G) \cap IO(G,G)$. Clearly, the inclusion
$LSm(G) \subseteq Sm(G) \cap IO(G,G)$ holds by the Basic Lemma.

Before we recall the notion of Oliver group, we wish to adopt the following
definition. For a given finite group $G$, a series of subgroups of $G$
of the form $P \trianglelefteq H \trianglelefteq G$ is called an
{\it isthmus series} if $|P| = p^m$ and $|G/H| = q^n$ for some
primes $p$ and $q$ (possibly $p = q$) and some integers $m, n \geq 0$,
and the quotient group $H/P$ is cyclic (possibly $H = P$).

For a finite group $G$, the following three claims are equivalent.
\begin{itemize}
\item[{\rm (1)}] $G$ has a smooth action on a sphere with exactly one fixed point.
\item[{\rm (2)}] $G$ has a smooth action on a disk without fixed points.
\item[{\rm (3)}] $G$ has no isthmus series of subgroups.
\end{itemize}

By the Slice Theorem, (1) implies (2). By the work of
Oliver \cite{Oliver:1}, (2) and (3) are equivalent, and according
to Laitinen and Morimoto \cite{Laitinen-Morimoto}, (3) implies (1).

Following Laitinen and Morimoto \cite{Laitinen-Morimoto},
a finite group $G$ is called an {\it Oliver group} if $G$ has
no isthmus series of subgroups. Recall that each finite nonsolvable
group $G$ is an Oliver group, and a finite abelian (more generally,
nilpotent) group $G$ is an Oliver group if and only if $G$ has three
or more noncyclic Sylow subgroups (cf.~\cite{Oliver:1}, \cite{Oliver:2},
and \cite{Laitinen-Morimoto-Pawalowski}).

We prove that the Laitinen Conjecture holds for large classes of
finite Oliver groups $G$ such that $a_G \geq 2$, and as a consequence,
we obtain that $Sm(G) \neq 0$. Moreover, we check that $Sm(G) = 0$
for specific classes of finite groups $G$ such that $a_G \leq 1$, and
therefore we can answer the Smith Isomorphism Question to the effect
that $Sm(G) = 0$ if and only if $a_G \leq 1$.

We wish to recall that for a finite group $G$, it may happen that
$Sm(G) \neq 0$ and $a_G \leq 1$ (the smallest group with these
properties is $G = \mathbb{Z}_8$).

\section*{0.2 \ \ Classification and Realization Theorems}

Our main algebraic theorem gives a classification of finite Oliver
groups $G$ with Laitinen number $a_G \leq 1$, and it reads as follows.

\begin{classification-thm}
Let $G$ be a finite Oliver group. Then the Laitinen number $a_G = 0$ or $1$
if and only if one of the following conclusions holds:
\begin{itemize}
\item[{\rm (1)}] $G \cong PSL(2,q)$ for some $q \in \{5, 7, 8, 9, 11, 13, 17\}$; or
\item[{\rm (2)}]  $G \cong PSL(3,3)$, $PSL(3,4)$, $Sz(8)$, $Sz(32)$, $A_7$,
	 $M_{11}$ or $M_{22}$; or
\item[{\rm (3)}]  $G \cong PGL(2,5)$, $PGL(2,7)$, $P\varSigma L(2,8)$, or $M_{10}$; or
\item[{\rm (4)}]  $G \cong PSL(3,4) \rtimes C_2
	 \cong PSL(3,4) \rtimes \langle u \rangle$; or
\item[{\rm (5)}]  $F(G) \cong C_2^2 \times C_3$ and
      $G \cong {\rm Stab}_{A_7}(\{1, 2, 3 \})$ or $C_2^2 \rtimes D_9$; or
\item[{\rm (6)}]  $F(G)$ is an abelian $p$-group for some odd prime $p$,
      $G \cong F(G) \rtimes H$ for $H < G$ with $H \cong SL(2,3)$ or
      $\hat{S_4}$, and $F(G)$ is inverted by the unique involution of $H$; or
\item[{\rm (7)}]  $F(G) \cong C_3^3$ and $G \cong F(G) \rtimes A_4$; or
\item[{\rm (8)}]  $F(G) \cong C_2^4$, $F^2(G) \cong A_4 \times A_4$, and
      $G \cong F^2(G) \rtimes C_4$; or
\item[{\rm (9)}]  $F(G) \cong C_2^8$ and $G \cong F(G) \rtimes H$ for $H < G$ with
      $H \cong PSU(3,2)$ or $C_3^2 \rtimes C_8$; or
\item[{\rm (10)}]  $F(G) \cong C_2^3$ and $G/F(G) \cong GL(3,2)$; or
\item[{\rm (11)}]  $F(G) \cong C_2^4$ and $G/F(G) \cong A_6$; or
\item[{\rm (12)}]  $F(G) \cong C_2^8$ and $G/F(G) \cong M_{10}$; or
\item[{\rm (13)}]  $F(G)$ is a non-identity elementary abelian $2$-group,
      $G/F(G) \cong SL(2,4)$, $\varSigma L(2,4)$, $SL(2,8)$, $Sz(8)$ or
      $Sz(32)$, and $C_{F(G)}(x) = 1$ for every $x \in G$ of odd order.
\end{itemize}
\end{classification-thm}

Here, we consider cyclic groups $C_q$ of order $q$, dihedral groups $D_q$
of order $2q$, elementary abelian $p$-groups $C_p^k = C_p \times \cdots \times C_p$,
alternating groups $A_n$, symmetric groups $S_n$, general linear groups $GL(n,q)$,
special linear groups $SL(n,q)$, projective general linear groups $PGL(n,q)$,
projective special linear groups $PSL(n,q)$, projective special unitary groups
$PSU(n,q)$, the Mathieu groups $M_{10}$, $M_{11}$, and $M_{22}$, and the Suzuki
groups $Sz(8)$ and $Sz(32)$. Recall that the group $PSL(3,4)$ admits an automorphism
$u$ of order $2$, referred to as a graph-field automorphism, acting as
the composition of the transpose-inverse automorphism and
the squaring map (a Galois automorphism) of the field $\mathbb{F}_4$
of four elements. The fixed points of $u$ form the group $PSU(3,2)$.

Moreover, for two finite groups $N$ and $H$, $N \rtimes H$ denotes a semi-direct product
of $N$ and $H$ (i.e., the splitting extension $G$ associated with an exact sequence
$1 \to N \to G \to H \to 1$). Also, we use the notations
$$\varSigma L(n,q) = SL(n,q) \rtimes {\rm Aut} (\mathbb{F}_q) \ \ \hbox{and} \ \
  P \varSigma L(n,q) = PSL(n,q) \rtimes {\rm Aut} (\mathbb{F}_q)$$
where ${\rm Aut} (\mathbb{F}_q)$ is the group of all automorphisms of
the field $\mathbb{F}_q$ of $q$ elements.

For $n \geq 4$ and $n \neq 6$, there exist two groups $G$ which are not
isomorphic, do not contain a subgroup isomorphic to $A_n$, and occur in
a short exact sequence $1 \to C_2 \to G \to S_n \to 1$. For $n = 4$, one
of the groups is isomorphic to $GL(2,3)$ and the other, denoted here by
$\hat{S_4}$, has exactly one element of order $2$.

Finally, for a finite group $G$, we denote by $F(G)$ the Fitting subgroup of $G$
(i.e., the largest nilpotent normal subgroup of $G$) and by $F^2(G)$
the pre-image of $F(G/F(G))$ under the quotient map $G \to G/F(G)$.

We remark that the Classification Theorem stated above extends a previous result
of Bannuscher and Tiedt \cite{Bannuscher-Tiedt} obtained for finite nonsolvable
groups $G$ such that every element of $G$ has prime power order (i.e., such that $a_G = 0$).
Our proof is largely independent of their result, but we do invoke it to establish that
$F(G)$ is elementary abelian in case (13). Moreover, their result and our cases (1)--(13)
allow us to list all finite Oliver groups $G$ with $a_G = 1$, and thus with $IO(G,G) = 0$
and $IO(G) \neq 0$ (cf.~the First Rank Lemma).

Let $G$ be a finite group. By \cite{Oliver:3}, there exists a smooth action of
$G$ on a disk with exactly two fixed points if and only if $G$ is an Oliver group.
For a finite Oliver group~$G$, two real $G$-modules $U$ and $V$ are called
{\it Oliver equivalent} if there exists a smooth action of $G$ on a disk $D$
such that $D^G = \{x,y\}$ for two points $x$ and $y$ at which $T_x(D) \cong U$
and $T_y(D) \cong V$ as real $G$-modules.

If $U - V \in RO(G)$ is the difference of two Oliver equivalent
real $G$-modules $U$ and $V$, then $U - V \in IO(G,G)$ by Smith theory and
the Slice Theorem. On the other hand, if $U - V \in IO(G,G)$, then $U$ and $V$
are isomorphic as $P$-modules for each subgroup $P$ of $G$ of prime power order,
and by subtracting the trivial summands, we may assume that
$\dim U^G = \dim V^G = 0$. Hence, by \cite[Theorem~0.4]{Oliver:3},
there exists a~smooth action of $G$ on a disk $D$ such that $D^G = \{x,y\}$
for two points $x$ and $y$ at which $T_x(D) \cong U \oplus W$ and
$T_y(D) \cong V \oplus W$ for some real $G$-module $W$ with
$\dim W^G = 0$. As in $RO(G)$,
$$U - V = (U \oplus W) - (V \oplus W),$$
the element $U - V$ is the difference of two Oliver equivalent real $G$-modules.
Consequently, $IO(G,G)$ coincides with the subset of $RO(G)$ consisting of
the differences of real $G$-modules which are Oliver equivalent. So, the
Classification Theorem and the First Rank Lemma yield the following corollary.

\begin{classification-cor}
A finite Oliver group $G$ has the property that two Oliver equivalent real
$G$-modules are always isomorphic (i.e., $IO(G,G) = 0$) if and only if $G$
is listed in cases (1)--(13) of the Classification Theorem
(i.e., the Laitinen number $a_G = 0$ or $1$).
\end{classification-cor}

For a finite group $G$, we denote by $\mathcal{P}(G)$ the family
of subgroups of $G$ consisting of the trivial subgroup of $G$
and all $p$-subgroups of $G$ for all primes $p \, \big| |G|$.

A subgroup $H$ of a finite group $G$ ($H \leq G$) is called a {\it large subgroup}
of $G$ if $O^p(G) \leq H$ for some prime $p$, where $O^p(G)$ is the smallest
normal subgroup of $G$ such that $|G/O^p(G)| = p^k$ for some integer
$k \geq 0$.

For a finite group $G$, we denote by $\mathcal{L}(G)$ the family of large
subgroups of $G$, and a real $G$-module $V$ is called {\it $\mathcal{L}$-free}
if $\dim V^H = 0$ for each $H \in \mathcal{L}(G)$, which amounts to saying
that $\dim V^{O^p(G)} = 0$ for each prime $p \, \big| |G|$.

Here, as in \cite{Morimoto-Sumi-Yanagihara}, a finite group $G$ is
called a {\it gap group} if $\mathcal{P}(G) \cap \mathcal{L}(G) = \varnothing$
and there exists a real $\mathcal{L}$-free $G$-module $V$ satisfying the
{\it gap condition} that
$$\dim V^P > 2 \dim V^H$$
for each pair $(P,H)$ of subgroups $P < H \leq G$ with $P \in \mathcal{P}(G)$.

According to \cite{Morimoto-Sumi-Yanagihara}, if $G$ is a finite group
such that $\mathcal{P}(G) \cap \mathcal{L}(G) = \varnothing$, then $G$
is a gap group under either of the following conditions:
\begin{itemize}
\item[{\rm (1)}] $O^p(G) \neq G$ and $O^q(G) \neq G$ for two distinct odd primes
      $p$ and $q$.
\item[{\rm (2)}] $O^2(G) = G$ (which is true when $G$ is of odd order or $G$
      is perfect).
\item[{\rm (3)}] $G$ has a quotient which is a gap group.
\end{itemize}
\noindent

Note that the condition (1) is equivalent to the condition that $G$ has
a cyclic quotient of order $pq$ for two distinct odd primes $p$ and $q$.
Recall that a finite group $G$ is nilpotent if and only if $G$ is the product
of its Sylow subgroups. Moreover, a finite nilpotent group $G$ is an Oliver
group if and only if $G$ has three or more noncyclic Sylow subgroups.
Therefore the condition (1) holds for any finite nilpotent Oliver group $G$.

If $G$ is a finite Oliver group, then
$\mathcal{P}(G) \cap \mathcal{L}(G) = \varnothing$ by \cite{Laitinen-Morimoto},
but it may happen that there is no real $\mathcal{L}$-free $G$-module
satisfying the gap condition. In fact, by \cite{Dovermann-Herzog} or
\cite{Morimoto-Sumi-Yanagihara}, the symmetric group $S_n$ is a gap group
if and only if $n \geq 6$. Hence, $S_5$ is an Oliver group which is not
a gap group, but $S_5$ contains $A_5$ which is both an Oliver and gap group.
We refer the reader to \cite{Morimoto-Sumi-Yanagihara},
\cite{Sumi:1} and \cite{Sumi:2} for more information about gap groups.

Let $LO(G)$ be the subgroup of $RO(G)$ consisting of the differences
$U - V$ of real $\mathcal{L}$-free $G$-modules $U$ and $V$ which are
isomorphic when restricted to any $P \in \mathcal{P}(G)$. Recall that
$IO(G)$ is the intersection of the kernels of the restriction maps
$RO(G) \to RO(P)$ taken for all $P \in \mathcal{P}(G)$, and $IO(G,G)$
is the intersection of $IO(G)$ and ${\rm Ker}\, (RO(G) \to \mathbb{Z})$
where $RO(G) \to \mathbb{Z}$ is the $G$-fixed point set dimension map.
In particular, $LO(G) \subseteq IO(G,G)$.

Now, we are ready to state our main topological theorem.

\begin{realization-thm}
Let $G$ be a finite Oliver gap group. Then any element of $LO(G)$
is the difference of two Laitinen--Smith equivalent real $G$-modules;
i.e., $LO(G) \subseteq LSm(G)$.
\end{realization-thm}

The Realization Theorem and the Basic Lemma show that
$$LO(G) \subseteq LSm(G) \subseteq IO(G,G)$$
for any finite Oliver gap group $G$. In general, $LO(G) \neq IO(G,G)$.
However, if $G$ is perfect, $O^p(G) = G$ for any prime $p$, and
hence $\mathcal{L}(G) = \{G\}$, and thus $LO(G) = IO(G,G)$.
So, the Realization Theorem and the Basic Lemma yield the following
corollary (cf.~\cite[Corollary~1.8]{Laitinen-Pawalowski} where a similar
result is obtained for the realifications of complex $G$-modules for
any finite perfect group $G$).

\begin{realization-cor}
Let $G$ be a finite perfect group. Then any element of $LO(G)$ is
the difference of two Laitinen--Smith equivalent real $G$-modules
and $LO(G) = IO(G,G)$, and thus $LO(G) = LSm(G) = IO(G,G)$.
\end{realization-cor}

\section*{0.3 \ \ Answers to the Smith Isomorphism Question}

By checking whether $Sm(G) = 0$, we answer the Smith Isomorphism Question
for large classes of finite Oliver groups $G$. In order to prove that $Sm(G) = 0$
if the Laitinen number $a_G \leq 1$, we use the Classification Theorem.
If the Laitinen number $a_G \geq 2$, we show that $LO(G) \neq 0$ and by the
Realization Theorem, we obtain that $LSm(G) \neq 0$, and thus $Sm(G) \neq 0$.

\begin{thm-A1}
Let $G$ be a finite Oliver group of odd order.
Then $a_G \geq 2$ and $LO(G) \neq 0$.
\end{thm-A1}

\begin{thm-A2}
Let $G$ be a finite group with a cyclic quotient of order $pq$ for two
distinct odd primes $p$ and $q$. Then $a_G \geq 2$ and $LO(G) \neq 0$.
\end{thm-A2}

\begin{thm-A3}
Let $G$ be a finite nonsolvable group. Then
\begin{itemize}
\item[{\rm (1)}] $LO(G) = 0$ if $a_G \leq 1$,
\item[{\rm (2)}] $LO(G) \neq 0$ if $a_G \geq 2$, except when $G \cong {\rm Aut}(A_6)$
                 or $P \varSigma L(2,27)$, and
\item[{\rm (3)}] $LO(G) = 0$ and $a_G = 2$ when $G \cong {\rm Aut}(A_6)$ or
                 $P \varSigma L(2,27)$.
\end{itemize}
\end{thm-A3}

\begin{thm-B1}
Let $G$ be a finite Oliver group of odd order. Then $a_G \geq 2$ and
$$0 \neq LO(G) \subseteq LSm(G) = Sm(G) \subseteq IO(G,G).$$
\end{thm-B1}

\begin{thm-B2}
Let $G$ be a finite Oliver group with a cyclic quotient of order $pq$
for two distinct odd primes $p$ and $q$. Then $a_G \geq 2$ and
$$0 \neq LO(G) \subseteq LSm(G) \subseteq IO(G,G).$$
\end{thm-B2}

\begin{thm-B3}
Let $G$ be a finite nonsolvable gap group not isomorphic to
$P \varSigma L(2,27)$. Then $LO(G) \neq 0$ if and only if $a_G \geq 2$,
$$LO(G) \subseteq LSm(G) \subseteq IO(G,G),$$
and $LSm(G) \neq 0$ if and only if $a_G \geq 2$.
\end{thm-B3}

By \cite[Theorem~A]{Laitinen-Pawalowski}, if $G$ is a finite perfect group,
$LSm(G) \neq 0$ if and only if $a_G \geq 2$. Theorem~B3 extends
this result in two ways. Firstly, it proves the conclusion for a large
class of finite nonsolvable groups $G$, including all finite perfect groups.
Secondly, if $G$ is perfect, it shows that $LSm(G) = IO(G,G)$
(cf.~the Realization Corollary).

If $G$ is as in Theorems~B1 or B2, the Laitinen Conjecture holds by the theorems.
By Theorem~B3, the Laitinen Conjecture holds for any finite nonsolvable gap group
$G$ with $a_G \geq 2$, except when $G \cong P \varSigma L(2,27)$. In
the exceptional case, $LO(G) = 0$ and $a_G = 2$ by Theorem~A3, and thus
${\rm rk}\,IO(G,G) = 1$ by the First Rank Lemma, so that $IO(G,G) \neq 0$.
However, we do not know whether $IO(G,G) \subseteq LSm(G)$, and we are not able
to confirm that $LSm(G) \neq 0$. The same is true when $G \cong {\rm Aut}(A_6)$.
Recall that $P \varSigma L(2,27)$ is a gap group while ${\rm Aut}(A_6)$ is not
a gap group (see \cite[Proposition~4.1]{Morimoto-Sumi-Yanagihara}).

\begin{thm-C1}
Let $G$ be a finite nonabelian simple group.
\begin{itemize}
\item[{\rm (1)}] If $a_G \leq 1$, then $Sm(G) = 0$ and $G$ is isomorphic
                 to one of the groups:

$a_G = 0:$ \ $PSL(2,q)$ for $q = 5, 7, 8, 9, 17$,
	     $PSL(3,4)$, $Sz(8)$, $Sz(32)$, or

$a_G = 1:$ \ $PSL(2,11)$, $PSL(2,13)$, $PSL(3,3)$, $A_7$,
	     $M_{11}$, $M_{22}$.

\item[{\rm (2)}] If $a_G \geq 2$, then $LSm(G) = IO(G,G) \neq 0$, and thus $Sm(G) \neq 0$.
\end{itemize}
\end{thm-C1}

\begin{thm-C2}
Let $G = SL(n,q)$ or $Sp(n,q)$ for $n \geq 2$ where $n$ is even in
the latter case and $q$ is any prime power in both cases.
\begin{itemize}
\item[{\rm (1)}] If $a_G \leq 1$, then $Sm(G) = 0$ and $G$ is isomorphic
                 to one of the groups:

$a_G = 0:$ \ $SL(2,2)$, $SL(2,4)$, $SL(2,8)$, $SL(3,2)$, or

$a_G = 1:$ \ $SL(2,3)$, $SL(3,3)$.

\item[{\rm (2)}] If $a_G \geq 2$, then except for $G = Sp(4,2)$,
                 $LSm(G) = IO(G,G) \neq 0$, and thus $Sm(G) \neq 0$.
                 Moreover, $Sm(G) \neq 0$ for $G = Sp(4,2)$.
\end{itemize}
\end{thm-C2}

\begin{thm-C3}
Let $G = A_n$ or $S_n$ for $n \geq 2$.
\begin{itemize}
\item[{\rm (1)}] If $a_G \leq 1$, then $Sm(G) = 0$ and $G$ is one of the groups:

$a_G = 0:$ \ $A_2$, $A_3$, $A_4$, $A_5$, $A_6$, $S_2$, $S_3$, $S_4$, or

$a_G = 1:$ \ $A_7$, $S_5$.

\item[{\rm (2)}] If $a_G \geq 2$, then $LSm(G) \supseteq LO(G) \neq 0$, and thus
                 $Sm(G) \neq 0$. Moreover, $LSm(G) = LO(G)$ for $G = A_n$.
\end{itemize}
\end{thm-C3}

We recall that $A_n$ is a simple group if and only if $n \geq 5$.
So, except for $A_2$, $A_3$ and $A_4$, every $A_n$ occurs in Theorem~C1.
Moreover, except for $PSL(2,2)$ and $PSL(2,3)$, every $PSL(n,q)$ is a simple group,
and the following holds: $A_5 \cong PSL(2,4) \cong PSL(2,5)$, $A_6 \cong PSL(2,9)$,
and $PSL(2,7) \cong PSL(3,2)$.

The symplectic group $Sp(n,q)$ and the projective symplectic group
$PSp(n,q)$ are defined for any even integer $n \geq 2$ and any prime
power $q$. Except for $PSp(2,2)$, $PSp(2,3)$, and $PSp(4,2)$,
every $PSp(n,q)$ is a nonabelian simple group, and thus occurs
in Theorem~C1. Moreover, in the exceptional cases, the following holds:
$PSp(2,2) \cong PSL(2,2) \cong S_3$, $PSp(2,3) \cong PSL(2,3) \cong A_4$,
and $PSp(4,2) \cong Sp(4,2) \cong S_6$. So, the cases are covered
by Theorem~C3.

\begin{com-D1}{\rm
The conjecture posed in \cite[p.~44]{Dovermann-Suh} asserts that if
$Sm(G) = 0$ for a finite group $G$, then $Sm(H) = 0$ for any subgroup
$H$ of $G$. We are able to give counterexamples to this conjecture. In fact,
according to Theorem~C1 and Example~E1 below, there exist (precisely four)
finite simple groups $G$ with an element of order $8$, such that $Sm(G) = 0$.
But $G$ has a subgroup $H \cong \mathbb{Z}_8$, and we know that $Sm(H) \neq 0$
by \cite{Cappell-Shaneson:1}--\cite{Cappell-Shaneson:3}.}
\end{com-D1}

\begin{com-D2}{\rm
Contrary to the speculation in \cite[Comment~(2), p.~547]{Schultz} that
$Sm(G) \neq 0$ for any finite Oliver group $G$, Theorem~C1 shows that
there exist (precisely fourteen) finite nonabelian simple groups $G$ such
that $Sm(G) = 0$. We recall that any finite nonabelian simple group $G$
is an Oliver group.}
\end{com-D2}

By using Theorems~B1--B3, we can answer the Smith Isomorphism Question
as follows: $Sm(G) \neq 0$ in either of the following cases.
\begin{itemize}
\item[{\rm (1)}] $G$ is a finite Oliver group of odd order (and thus $a_G \geq 2$).
\item[{\rm (2)}] $G$ is a finite Oliver group with a cyclic quotient of order $pq$
                 for two distinct odd primes $p$ and $q$ (and thus $a_G \geq 2$).
\item[{\rm (3)}] $G$ is a finite nonsolvable gap group with $a_G \geq 2$, and
                 $G \not\cong P \varSigma L(2,27)$.
\end{itemize}

In turn, Theorems~C1--C3 allow us to answer the Smith Isomorphism Question
as follows: $Sm(G) = 0$ if and only if $a_G \leq 1$, in either of
the following cases.

\begin{itemize}
\item[{\rm (1)}] $G$ is a finite nonabelian simple group.
\item[{\rm (2)}] $G = PSL(n,q)$ or $SL(n,q)$ for any $n \geq 2$
                 and any prime power $q$.
\item[{\rm (3)}] $G = PSp(n,q)$ or $Sp(n,q)$ for any even $n \geq 2$
                 and any prime power $q$.
\item[{\rm (4)}] $G = A_n$ or $S_n$ for any $n \geq 2$.
\end{itemize}

It follows from \cite[Theorem~B]{Laitinen-Pawalowski} that for $G = A_n$,
$PSL(2,p)$ or $SL(2,p)$ for any prime $p$, $Sm(G) = 0$ if and only if $a_G \leq 1$.
However, while \cite{Laitinen-Pawalowski} considers the realifications of complex
$G$-modules, we deal with real $G$-modules when proving that $Sm(G) \neq 0$ for
$a_G \geq 2$ (cf.~\cite[Corollary~1.8]{Laitinen-Pawalowski}).

By using the Realization Theorem, the Basic Lemma, the First Rank Lemma, and
Theorems~A1--A3, we are able to prove Theorems~B1--B3.

\begin{proof}[Proofs of Theorems~B1--B3]
Let $G$ be as in Theorems~B1--B3. Then, by the Realization Theorem and the Basic Lemma,
$$LO(G) \subseteq LSm(G) \subseteq IO(G,G).$$
If $G$ is as in Theorem~B1 (resp., B2), $a_G \geq 2$ and $LO(G) \neq 0$ by Theorem~A1
(resp., A2). Suppose that $G$ is as in Theorem~B3. According to our assumption,
$G \not\cong P \varSigma L(2,27)$ and $G \not\cong {\rm Aut}(A_6)$ as $G$ is
a gap group while ${\rm Aut}(A_6)$ is not (cf.~\cite[Proposition~4.1]{Morimoto-Sumi-Yanagihara}).
If $a_G \leq 1$, $IO(G,G) = 0$ by the First Rank Lemma, and thus $LO(G) = LSm(G) = 0$.
If $a_G \geq 2$, $LO(G) \neq 0$ by Theorem~A3, and thus $LSm(G) \neq 0$.
\end{proof}

Now, we adopt the following definition for any finite group $G$.
We say that $G$ satisfies the {\it $8$-condition}
if for every cyclic $2$-subgroup $P$ of $G$ with $|P| \geq 8$,
$\dim V^P > 0$ for any irreducible $G$-module $V$. In particular,
if $G$ is without elements of order $8$, $G$ satisfies the $8$-condition.
Recall that in \cite{Laitinen-Pawalowski}, $G$ satisfying the $8$-condition
is called $2$-proper (cf.~\cite[Example~2.5]{Laitinen-Pawalowski}).

If a finite group $G$ satisfies the $8$-condition and $G$ acts smoothly on a homotopy
sphere $\varSigma$ with $\varSigma^G \neq \varnothing$, then the action of $G$ on
$\varSigma$ satisfies the $8$-condition (cf.~Section~0.1), and thus
the following lemma holds (cf.~\cite[Lemma~2.6]{Laitinen-Pawalowski}).

\begin{8-lemma}
For each finite group $G$ satisfying the $8$-condition,
any two Smith equivalent real $G$-modules are also Laitinen--Smith
equivalent; i.e., $Sm(G) \subseteq LSm(G)$, and thus $Sm(G) = LSm(G)$.
\end{8-lemma}

\begin{exa-E1}{\rm
In the following list~(C1), each group $G$ satisfies the $8$-condition
and $a_G = 0$ or $1$, where $G$ is one of the groups:

$a_G = 0:$ \ $PSL(2,q)$ for $q = 2, 3, 5, 7, 8, 9, 17$,
	     $PSL(3,4)$, $Sz(8)$ or $Sz(32)$,

$a_G = 1:$ \ $PSL(2,11)$, $PSL(2,13)$, $PSL(3,3)$, $A_7$,
	     $M_{11}$ or $M_{22}$.

\noindent
If $G = PSL(2,2) \cong S_3$ or $G = PSL(2,3) \cong A_4$, then $a_G = 0$ and
$G$ has no element of order $8$ (cf.~\cite[Proposition~2.4]{Laitinen-Pawalowski}).
In list (C1), except for $PSL(2,2)$ and $PSL(2,3)$, every $G$ is
a nonabelian simple group, and some inspection in \cite{Conway-Wilson}
or \cite{Gorenstein-Lyons-Solomon:3} confirms that $a_G = 0$ for $G = PSL(2,q)$
with $q = 5, 7, 8, 9, 17$, and $a_G = 0$ for $G = PSL(3,4)$, $Sz(8)$ or $Sz(32)$.
Also, $a_G = 1$ corresponding to an element of order $6$ when
$G = PSL(2,11)$, $PSL(2,13)$, $PSL(3,3)$, $A_7$, $M_{11}$ or $M_{22}$.
Further inspection in \cite{Conway-Wilson} or \cite{Gorenstein-Lyons-Solomon:3}
shows that in list (C1), $G$ has an element of order $8$ if and only if
$G = PSL(2,17)$, $PSL(3,3)$, $M_{11}$ or $M_{22}$, and the groups all satisfy
the $8$-condition. All finite groups $G$ without elements of order $8$ also satisfy
the $8$-condition. Therefore, each group $G$ in list (C1) satisfies the $8$-condition.}
\end{exa-E1}

\begin{exa-E2}{\rm
In the following list~(C2), each group $G$ satisfies the $8$-condition
and $a_G = 0$ or $1$, where $G$ is one of the groups:

$a_G = 0:$ \ $SL(2,2)$, $SL(2,4)$, $SL(2,8)$, $SL(3,2)$,
		 $Sp(2,2)$, $Sp(2,4)$ or $Sp(2,8)$,

$a_G = 1:$ \ $SL(2,3)$, $SL(3,3)$ or $Sp(2,3)$.

\noindent
As $Sp(2,q) \cong SL(2,q)$ for any prime power $q$, it sufficies to check
the result for the special linear groups. First, recall that $SL(2,q) \cong PSL(2,q)$
when $q$ is a power of $2$. Clearly, $a_G = 0$ when $G = SL(2,2) \cong PSL(2,2) \cong S_3$,
and by Example~E1, $a_G = 0$ when $G = SL(2,4) \cong PSL(2,4) \cong PSL(2,5) \cong A_5$, or
$G = SL(2,8) \cong PSL(2,8)$, or $G = SL(3,2) \cong PSL(3,2) \cong PSL(2,7)$.
Moreover, for $G = SL(3,3) \cong PSL(3,3)$, $a_G = 1$ corresponding to
an element of order $6$. The same holds for $G = SL(2,3)$ because $G$ has
elements of orders $1$, $2$, $3$, $4$, and $6$, and the elements of order $6$
are all real conjugate in $G$ (cf.~\cite[Proposition~2.3]{Laitinen-Pawalowski}).
By the discussion above and Example~E1, we see that in list (C2),
$G$ has an element of order $8$ if and only if $G = SL(3,3)$, and
$SL(3,3) \cong PSL(3,3)$ satisfies the $8$-condition. So, each group $G$
in list (C2) satisfies the $8$-condition.}
\end{exa-E2}

\begin{exa-E3}{\rm
In the following list (C3), each group $G$ is without elements of order $8$
and $a_G = 0$ or $1$, or $a_G \geq 2$, where $G$ is one of the groups:

$a_G = 0:$ \ $A_2$, $A_3$, $A_4$, $A_5$, $A_6$, $S_2$, $S_3$ or $S_4$,

$a_G = 1:$ \ $A_7$ or $S_5$,

$a_G \geq 2:$ \ $A_8$, $A_9$, $S_6$ or $S_7$.

\noindent
First, we consider the case $G = A_n$. For $n \leq 6$, $a_G = 0$
because each element of $G$ has prime power order. For $n = 7$,
$a_G = 1$ corresponding to the element $(12)(34)(567)$ of order $6$.
For $n \geq 8$, $a_G \geq 2$ because the elements $(12)(34)(567)$
and $(123456)(78)$ have order $6$ and are not real conjugate in $G$.

Now, we consider the case $G = S_n$. For $n \leq 4$, $a_G = 0$
because each element of $G$ has prime power order. For $n = 5$,
$a_G = 1$ corresponding to the element $(12)(345)$ of order $6$.
For $n \geq 6$, $a_G \ge 2$ because the elements $(12)(345)$
and $(123456)$ have order $6$ and are not real conjugate in $G$.

As a result, if $G = A_n$ (resp., $S_n$), $a_G \leq 1$ if and only if
$n \leq 7$ (resp., $n \leq 5$). Moreover, if $G = A_n$ or $S_n$ for
$n \leq 7$, $G$ has no element of order $8$ because any permutation
of order $8$ must involve an $8$-cycle in its cycle decomposition.
Also, if $G = A_8$ or $A_9$, $G$ has no element of order $8$
because an $8$-cycle is not an even permutation. Therefore,
each group $G$ in list (C3) is without elements of order $8$,
and thus $G$ satisfies the $8$-condition.}
\end{exa-E3}

By using the Classification Theorem, Examples~E1--E3, the 8-condition Lemma,
the Basic Lemma, the First Rank Lemma, and Theorems~B1--B3, we are able
to prove Theorems~C1--C3.

\begin{proof}[Proofs of Theorems~C1--C3]
Let $G$ be as in Theorems~C1--C3. Then, by the Classification Theorem
and Examples~E1--E3, $a_G = 0$ or $1$ if and only if $G$ is as in
claims (1) of Theorems~C1--C3.

If $a_G \leq 1$, $G$ satisfies the $8$-condition by Examples~E1--E3, and thus
$$Sm(G) = LSm(G) = IO(G,G) = 0$$
by the 8-condition Lemma, the Basic Lemma and the First Rank Lemma.

If $a_G \geq 2$, $G$ is as in Theorems~B1--B3, and therefore
$$0 \neq LO(G) \subseteq LSm(G) \subseteq IO(G,G).$$ Moreover, except
for $G = S_n$ or $Sp(4,2) \cong S_6$, $G$ is a perfect group, and thus
$LO(G) = LSm(G) = IO(G,G)$ (cf.~the Realization Corollary obtained
from the Realization Theorem and the Basic Lemma).\end{proof}

\section*{0.4 \ \ Second Rank Lemma}

Let $G$ be a finite group. In Sections 0.1 and 0.2, we defined
the following series of free abelian subgroups of $RO(G)$:
$LO(G) \subseteq IO(G,G) \subseteq IO(G)$. Recall that $IO(G)$ consists of
the differences $U - V$ of real $G$-modules $U$ and $V$ which are isomorphic
when restricted to any $P \in \mathcal{P}(G)$, $IO(G,G)$ is obtained from $IO(G)$
by imposing the additional condition that $\dim U^G = \dim V^G$, and $LO(G)$
consists of the differences $U - V \in IO(G)$ such that $U$ and $V$ are both
$\mathcal{L}$-free. Now, for any normal subgroup $H$ of $G$, we put
$$IO(G,H) = IO(G) \cap \, {\rm Ker} \,(RO(G)
	    \overset{{\rm Fix}^H}{\longrightarrow} RO(G/H))$$
where ${\rm Fix}^H (U - V) = U^H - V^H$ and the $H$-fixed point sets $U^H$
and $V^H$ are considered as the canonical $G/H$-modules.
As $RO(G/G) \cong \mathbb{Z}$ and
$${\rm Ker} \,(RO(G) \overset{{\rm Fix}^G}{\longrightarrow} RO(G/G)) =
  {\rm Ker} \,(RO(G) \overset{{\rm Dim}^G}{\longrightarrow} \mathbb{Z}),$$
the two definitions of $IO(G,G)$ coincide. In general, $IO(G,H) \subseteq IO(G,G)$.
In fact, if $U - V \in IO(G,H)$, then $U - V \in IO(G)$ and in addition
$U^H \cong V^H$ as $G/H$-modules, so that
$\dim U^G = \dim (U^H)^{G/H} = \dim (V^H)^{G/H} = \dim V^G$,
proving that $U - V \in IO(G,G)$. Therefore $IO(G,H) \subseteq IO(G,G)$.

Henceforth, we denote by $b_{G/H}$ the number of real conjugacy
classes $(gH)^{\pm 1}$ in $G/H$ of cosets $gH$ containing elements
of $G$ not of prime power order.

In general, $a_G \geq b_{G/H} \geq {a}_{G/H}$.
Clearly, $a_G = b_{G/G} = 0$ when each element of $G$ has
prime power order, and $a_G = b_{G/G} = 1$ when $G$ has
elements not of prime power order and any two such elements are real
conjugate in $G$. Otherwise, $a_G > b_{G/G} = 1$. Therefore,
$a_G = b_{G/G}$ if and only if $a_G = 0$ or $1$.

We compute the rank ${\rm rk} \,IO(G,H)$. For $H = G$,
the computation goes back to \cite[Lemma~2.1]{Laitinen-Pawalowski}
(cf.~the First Rank Lemma in Section~0.1 of this paper).

\begin{second}
Let $G$ be a finite group and let $H \trianglelefteq G$. Then
$${\rm rk} \,IO(G,H) = a_G - b_{G/H} \ \ \hbox{and thus} \ \
  {\rm rk} \,IO(G,G) = a_G - b_{G/G}.$$
In particular, $IO(G,H) = 0$ if $a_G \leq 1$, and $IO(G,G) = 0$
if and only if $a_G \leq 1$.
\end{second}

\begin{proof}
In \cite[Lemma~2.1]{Laitinen-Pawalowski}, the rank of $IO(G)$ is computed
as follows. The rank of the free abelian group $IO(G)$ is equal to
the dimension of the real vector space $\mathbb{R} \otimes_{\mathbb{Z}} IO(G)$
which consists of the real valued functions on $G$ that are constant on
the real conjugacy classes $(g)^{\pm 1}$ and that vanish when $g$ is of
prime power order. Therefore ${\rm rk} IO(G) = a_G$.

Now, for a normal subgroup $H$ of $G$, we compute the rank of the kernel
$$IO(G,H) =  {\rm Ker} \,(IO(G) \overset{{\rm Fix}^H}{\longrightarrow} RO(G/H)).$$
First, for any representation $\rho: G \to {\rm GL}(V)$, consider the representation
${\rm Fix}^H \rho : G/H \to {\rm GL}(V^H)$ given by $({\rm Fix}^H \rho)(gH) = \rho(g)|_{V^H}$
for each $g \in G$. Let $\pi: V \to V$ be the projection of $V$ onto $V^H$, that is,
$$\pi = \frac{1}{| H |}\sum_{h \in H} \rho(h) : \ V \to V.$$
Then the trace of $({\rm Fix}^H \rho) (gH): V^H \to V^H$
is the same as the trace of the endomorphism
$$\rho (g) \circ \pi = \frac{1}{| H |}\sum_{h \in H} \rho(gh): \
  V \to V.$$
So, if $\chi$ is the character of $\rho$, then the character
${\rm Fix}^H \chi$ of ${\rm Fix}^H \rho$ is given by
$$({\rm Fix}^H \chi) (gH) = \frac{1}{| H |}\sum_{h \in H} \chi (gh).$$
This formula extends (by linearity) to $\mathbb{R} \otimes_{\mathbb{Z}}RO(G)$.
Now, consider the basis of $\mathbb{R} \otimes_{\mathbb{Z}}IO(G)$ consisting
of the functions $f_{(g)^{\pm 1}}$ which have the value $1$ on $(g)^{\pm 1}$
and $0$ otherwise, defined for all classes $(g)^{\pm 1}$ represented
by elements $g \in G$ not of prime power order. Then, by the formula
above applied to
$\chi = f_{(g)^{\pm 1}}$,
$$({\rm Fix}^H f_{(g)^{\pm 1}}) (gH) = \frac{| (g)^{\pm 1} \cap gH |}
							 {| H |}$$
and ${\rm Fix}^H f_{(g)^{\pm 1}}$ vanishes outside of $(gH)^{\pm 1}$.
Therefore, the map
$${\rm Fix}^H: IO(G) \to RO(G/H)$$
has image of rank $b_{G/H}$, and its kernel $IO(G,H)$ is of rank
$a_G - b_{G/H}$.
\end{proof}

We wish to note that if $G$ is a finite group and $H \triangleleft G$
(i.e. $H \trianglelefteq G$ and $H \neq G$), then one of the following
conclusions holds:

\begin{itemize}
\item[{\rm (1)}] $a_G = b_{G/H} = 0$ if each $g \in G$ has prime power order, and otherwise
\item[{\rm (2)}] $a_G = b_{G/H} = 1$ (holds, e.g., for $G = S_5$ and $H = G^{{\rm sol}} = A_5$), or
\item[{\rm (3)}] $a_G = b_{G/H} > 1$ (holds, e.g., for $G = {\rm Aut}(A_6)$ and $H = G^{{\rm sol}}$), or
\item[{\rm (4)}] $a_G > b_{G/H} = 1$ (holds, e.g., for $G = S_6$ and $H = G^{{\rm sol}} = A_6$), or
\item[{\rm (5)}] $a_G > b_{G/H} > 1$ (holds, e.g., for $G = A_5 \times \mathbb{Z}_3$ and
                                     $H = G^{{\rm sol}} = A_5$).
\end{itemize}

Let $G$ be a finite group with two subgroups $H \trianglelefteq G$
and $K \trianglelefteq G$. We claim that if $H$ is a subgroup of $K$, $H \leq K$,
then $IO(G,H)$ is a subgroup of $IO(G,K)$. In fact, take an element
$$U - V \in IO(G,H) = {\rm Ker} \,(IO(G) \overset{{\rm Fix}^H}{\longrightarrow} RO(G/H))$$
and consider the $G$-orthogonal complements $U - U^H$ and $V - V^H$ of
the real $G$-modules $U$ and $V$. Then $U - V = (U - U^H) - (V - V^H)$
because $U^H \cong V^H$ as $G/H$-modules, and
$(U - U^H)^K = (V - V^H)^K = \{0\}$ because $H \leq K$. Therefore, it follows that
$$U - V = (U - U^H) - (V - V^H) \in IO(G,K) =
{\rm Ker} \,(IO(G) \overset{{\rm Fix}^K}{\longrightarrow} RO(G/K)),$$
proving the claim that $IO(G,H) \subseteq IO(G,K)$.

For any finite group $G$, we consider the group $IO(G,H)$,
where $H$ is:
$$
\aligned
G^{{\rm sol}}&: \ \hbox{the smallest normal subgroup of $G$
                          such that $G/H$ is solvable,}\\
G^{{\rm nil}}&: \ \hbox{the smallest normal subgroup of $G$
                          such that $G/H$ is nilpotent,}\\
O^p(G)         &: \ \hbox{the smallest normal subgroup of $G$
                          such that $G/H$ is a $p$-group.}
\endaligned
$$

Clearly, $G$ is perfect if and only if $G^{{\rm sol}} = G$, and
$G$ is solvable if and only if $G^{{\rm sol}}$ is trivial. And similarly,
$G$ is nilpotent if and only if $G^{{\rm nil}}$ is trivial. Moreover,
$G^{{\rm sol}} \subseteq G^{{\rm nil}} = \bigcap_p O^p(G)$ taken for
all primes $p \, \big| |G|$.

\begin{subgroup}
Let $G$ be a finite group and let $p$ be a prime.
Then
$$IO(G,G^{{\rm sol}}) \subseteq IO(G,G^{{\rm nil}}) \subseteq LO(G)
		 \subseteq IO(G,O^p(G)) \subseteq IO(G,G).$$
\end{subgroup}

\begin{proof}
By the claim above, $IO(G,G^{{\rm sol}}) \subseteq IO(G,G^{{\rm nil}})$
because $G^{{\rm sol}} \subseteq G^{{\rm nil}}$.
Now, set $H = G^{{\rm nil}}$. For a real $G$-module $V$, consider $V^H$
as a real $G$-module with the canonical action of $G$. Then the $G$-orthogonal
complement $V - V^H$ of $V^H$ in $V$ is $\mathcal{L}$-free because
$H \subseteq O^p(G)$ for each prime $p$. Take an element $U - V \in IO(G,H)$.
Then $U^H \cong V^H$ as $G$-modules, so that
$$U - V = (U - U^H) - (V - V^H) \in LO(G),$$
proving that $IO(G,G^{{\rm nil}}) \subseteq LO(G)$. Any element of $LO(G)$ is
the difference of two real $\mathcal{L}$-free $G$-modules $U$ and $V$ such that
$U - V \in IO(G)$. As $U$ and $V$ are $\mathcal{L}$-free, $\dim U^{O^p(G)} = \dim V^{O^p(G)} = 0$,
and thus $U - V \in IO(G,O^p(G))$, proving that $LO(G) \subseteq IO(G,O^p(G))$.
Clearly, $IO(G,O^p(G)) \subseteq IO(G,G)$ by the claim above.
\end{proof}

By the Subgroup Lemma and the Second Rank Lemma,
$$a_G - b_{G/G^{{\rm nil}}} \leq {\rm rk}\,LO(G) \leq
  \min \{a_G - b_{G/O^p(G)} : \, p \, \big| |G|\}$$
for any finite group $G$. In particular,
$a_G - b_{G/G^{{\rm sol}}} \leq {\rm rk}\,LO(G) \leq a_G - b_{G/G}$.

\begin{exa-E4}{\rm
Let $G = A_n$ for $n \geq 2$. By the First Rank Lemma, we know that
${\rm rk} \,IO(G,G) = 0$ when $a_G \leq 1$, and ${\rm rk} \,IO(G,G) = a_G - 1$
when $a_G \geq 1$. Moreover, by Theorem~C3,
$$Sm(G) \supseteq LSm(G) = LO(G) = IO(G,G).$$
Now, assume that $G = A_8$ or $A_9$. Then $G$ has no element of order $8$,
and thus $Sm(G) = LSm(G)$. By straightforward computation, we check that
$a_G = 3$ (resp., $6$) for $G = A_8$ (resp., $A_9$). As a result, we obtain that
\begin{itemize}
\item[{\rm (1)}] $Sm(A_8) = LSm(A_8) = LO(A_8) = IO(A_8,A_8) \cong \mathbb{Z}^2$ and
\item[{\rm (2)}] $Sm(A_9) = LSm(A_9) = LO(A_9) = IO(A_9,A_9) \cong \mathbb{Z}^5$.
\end{itemize}}
\end{exa-E4}

Generalizing the case where $G = A_8$ or $A_9$, note that the 8-condition Lemma
and the Realization Corollary yield the following corollary.

\begin{8-corollary}
Let $G$ be a finite group satisfying the $8$-condition. If $G$ is perfect,
then $Sm(G) = LSm(G) = LO(G) = IO(G,G)$.
\end{8-corollary}

\begin{exa-E5}{\rm
Let $G = S_n$ and $H = A_n$ for $n \geq 2$. Then
$G^{{\rm sol}} = H = O^2(G)$ and $O^p(G) = G$ for each odd prime $p$.
Therefore, $LO(G) = IO(G,H)$ by the Subgroup Lemma, and
${\rm rk} \,LO(G) = a_G - b_{G/H}$ by the Second Rank Lemma.
It follows from Example~E3 that $b_{G/H} = 0$ for $n = 2$, $3$ or $4$,
$b_{G/H} = 1$ for $n = 5$ or $6$, and $b_{G/H} = 2$ for $n \geq 7$.
Also, $a_G = 0$ for $n = 2$, $3$ or $4$, and $a_G = 1$ for $n = 5$.
Thus, ${\rm rk} \,LO(G) = a_G - b_{G/H} = 0$ for $n = 2$, $3$, $4$ or $5$.
For $n \geq 6$, $a_G \geq 2$ and by Theorem~C3 and the Basic Lemma,
we see that $0 \neq LO(G) \subseteq LSm(G) \subseteq IO(G,G)$.

Now, let $G = S_6$ (resp., $S_7$) and let $H \triangleleft G$ be as above.
By straightforward computation, we check that $a_G = 2$ (resp., $5$).
As we noted above, $b_{G/H} = 1$ (resp., $2$), and thus
${\rm rk} \,LO(G) = a_G - b_{G/H} = 1$ (resp., $3$).
Moreover, by the First Rank Lemma, ${\rm rk} \,IO(G,G) = a_G - 1 = 1$ (resp., $4$).
As $G$ has no element of order $8$, $Sm(G) = LSm(G)$. As a result, we obtain that
\begin{itemize}
\item[{\rm (1)}] $Sm(S_6) = LSm(S_6) = LO(S_6) = IO(S_6,S_6) \cong \mathbb{Z}$ and
\item[{\rm (2)}] $Sm(S_7) = LSm(S_7) \supseteq LO(S_7) \cong \mathbb{Z}^3$ and
                 $IO(S_7,S_7) \cong \mathbb{Z}^4$.
\end{itemize}}
\end{exa-E5}

\section*{0 \ \ Outline of material}

Let $G$ be a finite group. For the convenience of the reader,
we give a glossary of subsets and subgroups (defined above) of
the real representation ring $RO(G)$. First, recall that the following
two subsets of $RO(G)$ consist of the differences $U - V$ of real $G$-modules
$U$ and $V$ such that:
$$
\aligned
Sm(G)  &: \  \hbox{$U$ and $V$ are Smith equivalent}, \\
LSm(G) &: \  \hbox{$U$ and $V$ are Laitinen--Smith equivalent}.
\endaligned
$$

The following four subgroups of $RO(G)$ consist of the differences $U - V$
of real $G$-modules $U$ and $V$ such that $U \cong V$ as $P$-modules for each
$P \in \mathcal{P}(G)$, and:
$$
\aligned
IO(G)   &: \ \hbox{there is no additional restriction on $U$ and $V$}, \\
LO(G)   &: \ \hbox{the $G$-modules $U$ and $V$ are both $\mathcal{L}$-free}, \\
IO(G,G) &: \ \hbox{$\dim U^G = \dim V^G$}, \\
IO(G,H) &: \ \hbox{$U^H \cong V^H$ as $G/H$-modules},
\endaligned
$$
where $IO(G,H)$ is defined for any normal subgroup $H$ of $G$.

In Section~0.1, for a finite group $G$, we recalled the question of
Paul~A.~Smith about the tangent $G$-modules for smooth actions of $G$
on spheres with exactly two fixed points. Then we stated
the Basic Lemma and the First Rank Lemma. Moreover, we restated
the Laitinen Conjecture from \cite{Laitinen-Pawalowski}.

In Section~0.2, we stated the Classification and Realization Theorems
(our main algebraic and topological theorems) and by using the theorems,
we obtained the Classification and Realization Corollaries.

In Section~0.3, we stated Theorems~A1--A3, B1--B3, and C1--C3. We answered
the Smith Isomorphism Question and confirmed that the Laitinen Conjecture
holds for many groups $G$. Then we have proved that Theorems~B1--B3
follow from the Realization Theorem, the Basic Lemma, the First Rank Lemma,
and Theorems~A1--A3. Moreover, we stated the $8$-condition Lemma and we gave
Examples~E1--E3. Finally, we have proved that Theorems~C1--C3 follow from
the Classification Theorem, Examples~E1--E3, the 8-condition Lemma,
the Basic Lemma, the First Rank Lemma, and Theorems~B1--B3.

In Section~0.4, we stated and proved the Second Rank Lemma and the Subgroup Lemma.
We also gave Examples~E4 and E5 with $G = A_n$ and $S_n$, respectively. Moreover,
we obtained the $8$-condition Corollary for any finite perfect group $G$ satisfying
the $8$-condition.

As we pointed out above, the Basic Lemma, the First Rand Lemma, and
the $8$-condition Lemma all three go back to \cite{Laitinen-Pawalowski}.
Therefore, it remains to prove the Classification Theorem, the Realization Theorem,
and Theorems~A1--A3.

In Section~1, we prove Theorems~A1 and A2. To prove Theorem~A1, we obtain
our first major result about the Laitinen number $a_G$. The result asserts
that if $G$ is a finite Oliver group of odd order and without cyclic
quotient of order $pq$ for two distinct odd primes $p$ and $q$, then
$a_G > b_{G/{G^{{\rm nil}}}}$ (Proposition~1.6), and thus $LO(G) \neq 0$
by the Second Rank Lemma and the Subgroup Lemma. If $G$ is a finite group
with a cyclic quotient of order $pq$ for two distinct odd primes
$p$ and $q$, then $a_G \geq 4$ and $LO(G) \neq 0$ by an explicit
construction of two real $\mathcal{L}$-free $G$-modules $U$ and $V$,
which we give at the end of Section~1. This completes the proof
of Theorem~A1, and proves Theorem~A2.

In Section~2, we prove the Classification Theorem by using the fundamental results
of \cite{Gorenstein-Lyons-Solomon:1}--\cite{Gorenstein-Lyons-Solomon:3}, including
those restated in Theorems~2.2--2.4 of this paper, as well as by using
Burnside's $p^aq^b$ Theorem, the Feit--Thompson Theorem, the Brauer--Suzuki Theorem,
and the Classification of the finite simple groups.

In Section~3, we prove Theorem~A3. To present the proof, we analyze first
the cases where $a_G = b_{G/G^{{\rm sol}}}$. As a result, we obtain
our next major result about the Laitinen number $a_G$. The result asserts
that if $G$ is a finite nonsolvable group with $a_G = b_{G/G^{{\rm sol}}}$,
then either $a_G \leq 1$ or $a_G = 2$ and $G \cong {\rm Aut}(A_6)$ or
$P \varSigma L(2,27)$ (Proposition~3.1). By using the Second Rank Lemma,
this allows us to find the cases where $IO(G,G^{{\rm sol}}) \neq 0$
(Corollary~3.13), and then by using the Subgroup Lemma, we are able
to complete the proof of Theorem~A3.

In Section~4, we prove the Realization Theorem. To present the proof,
we recall first in Theorems~4.1 and 4.2 some equivariant thickening
and surgery results which follow from \cite{Morimoto-Pawalowski:2} and
\cite{Morimoto-Pawalowski:3}, respectively. Then, in Theorems~4.3 and 4.4,
we construct smooth actions of $G$ on spheres with prescribed real $G$-modules
at the fixed points. The required proof follows easily from Theorem~4.4.

We use information from \cite{Bredon}, \cite{Dieck}, \cite{Kawakubo}
on transformation group theory and from
\cite{Curtis-Reiner:1}, \cite{Curtis-Reiner:2},
\cite{Gorenstein-Lyons-Solomon:1}--\cite{Gorenstein-Lyons-Solomon:5},
\cite{Huppert-Blackburn}, \cite{James-Liebeck} on group theory and
representation theory.

\begin{ack}{\rm
The first author was supported in part by KBN Grant No.~2\,P03A\,031\,15,
and the second author was supported in part by NSF Grant No.~0070801.
The authors would like to thank Yu-Fen Wu for calling the results
of Bannuscher and Tiedt \cite{Bannuscher-Tiedt} to their attention,
and the first author would like to thank the Department of Mathematics
at The Ohio State University for its hospitality and support during
his visit to the department. Also, both authors would like to express
their thanks to the referee for his critical comments which allowed
the authors to improve the presentation of the material.}
\end{ack}

\section{Proofs of Theorems A1 and A2}

Let $G$ be a finite group. We denote by ${\rm NPP}(G)$ the set of
elements $g$ of $G$ which are not of prime power order, and we refer
to the elements of ${\rm NPP}(G)$ as NPP elements of $G$. Also, we
denote by $\overline{{\rm NPP}}(G)$ the set of real conjugacy classes
which are subsets of ${\rm NPP}(G)$. Therefore, the Laitinen number
$a_G$ is the number of elements in $\overline{{\rm NPP}}(G)$.

Let $H \trianglelefteq G$. Then, by the Second Rank Lemma,
$IO(G,H) \neq 0$ if and only if $a_G > b_{G/H}$. Clearly,
$a_G > b_{G/H}$ if and only if ${\rm NPP}(G)$ contains
two elements $x$ and $y$ not real conjugate in $G$, but such that
the cosets $xH$ and $yH$ are real conjugate in $G/H$.

\begin{lemma}
Let $H \trianglelefteq G$. Then the following three conclusions hold.
\begin{itemize}
\item[{\rm (1)}] Some coset $gH$ meets two members of $\overline{{\rm NPP}}(G)$
                 if and only if $a_G > b_{G/H}$.
\item[{\rm (2)}] If $H$ contains two distinct members of $\overline{{\rm NPP}}(G)$,
                 then $a_G > b_{G/H}$.
\item[{\rm (3)}] If $a_G = b_{G/H}$, then $a_{G/K} = b_{(G/K)/(H/K)}$ for any
                 $K \trianglelefteq G$ with $K \subseteq H$.
\end{itemize}
\end{lemma}

\begin{proof}
The first conclusion is immediate from the remarks above, while
the second one is a special case of the first. To prove
the third conclusion, suppose that $a_{G/K} > b_{(G/K)/(H/K)}$.
Then some coset $\overline{g} (H/K)$ meets two members of
$\overline{{\rm NPP}}(G/K)$. Assume that $x$ and $y$ are two elements
of $G$ such that $\overline{x}$ and $\overline{y}$ are not of prime power
order in $G/K$ and are not in the same real conjugacy class of $G/K$.
Then neither $x$ nor $y$ is of prime power order and $xH = yH$.
If $zxz^{-1} \in \{y, y^{-1}\}$ for an element $z \in G$, then
$\overline{z} \, \overline{x} \, \overline{z}^{-1} \in
\{\overline{y}, \overline{y}^{-1}\}$ contrary to assumption. Therefore,
$x$ and $y$ are not in the same real conjugacy class of $G$, and thus
$a_G > b_{G/H}$ by the first conclusion, proving the third one.
\end{proof}

\begin{lemma}
Let $G$ be a finite group and assume that
$K \trianglelefteq L \subseteq H \subseteq G$ is a sequence of subgroups
of $G$ such that $L/K$ contains NPP elements of two different orders.
Then $H$ contains NPP elements of two different orders, and
$a_G > b_{G/H} \geq a_{G/H}$ when $H \trianglelefteq G$.
\end{lemma}

\begin{proof}
Suppose $xK$ and $yK$ are NPP elements of $L/K$ of different orders.
If the elements $x$ and $y$ have different orders, we are done. If not,
we may assume that the order of $x$ is larger than the order of $xK$,
in which case the cyclic group generated by $x$ contains two NPP elements
of different orders. So in any case, $H$ contains two NPP elements of
different orders. Moreover, if $H \trianglelefteq G$, then $a_G > b_{G/H}$
by Lemma~1.1. Clearly $b_{G/H} \geq a_{G/H}$.
\end{proof}

\begin{lemma}
Let $G$ be a finite group containing a nonsolvable subgroup $B$ and
a cyclic subgroup $C \neq 1$ such that $BC$ is a subgroup of $G$
isomorphic to $B \times C$. Then $G$ has NPP elements of different
orders, and thus $a_G \geq 2$. Moreover, if $B \subseteq G^{{\rm sol}}$,
then $a_G > b_{G/G^{{\rm sol}}}$.
\end{lemma}

\begin{proof}
For a prime divisor $p$ of the order of $C$, choose an element
$g \in C$ of order $p$. Since $B$ is nonsolvable, it follows from
Burnside's $p^aq^b$ Theorem that the order of $B$ has (at least)
three distinct prime divisors $q$, $r$, and $s$, say with $p \ne r$
and $p \ne s$. Choose two elements $x$ and $y$ in $B$ of orders
$r$ and $s$, respectively.
By the assumption, $BC$ is a subgroup of $G$ (which amounts to
saying that $BC = CB$) isomorphic to $B \times C$. Thus, the
elements $gx$ and $gy$ have orders $pr$ and $ps$, respectively,
proving that $a_G \geq 2$. If $B \subseteq G^{{\rm sol}}$, then
the coset $g G^{{\rm sol}}$ contains the elements $gx$ and $gy$ which
are not real conjugate in $G$, and thus $a_G > b_{G/G^{{\rm sol}}}$
by Lemma~1.1.
\end{proof}

\begin{lemma}
Let $G$ be a finite group of odd order and let $H \trianglelefteq G$.
Suppose that $p$ is a prime and $P$ is an abelian $p$-subgroup of $H$
with $P \trianglelefteq G$. Suppose also that $q$ is a prime, $q \neq p$,
and $x \in H$ of order $q$ with $V = C_P(x) \neq 1$. Then $a_G > b_{G/H}$.
\end{lemma}

\begin{proof}
Suppose $a_G = b_{G/H}$. By Lemma~1.1, $H$ contains at most one member
of $\overline{{\rm NPP}}(G)$. On the other hand, every element of
$Vx \smallsetminus \{ x \}$ has order $pq$ and so all
of these elements lie in the same member of $\overline{{\rm NPP}}(G)$.
Take $y, y' \in V \smallsetminus \{1\}$ and set $h = yx$. Now, take
$g \in G$ with $h^g = ghg^{-1} \in \{y'x, (y'x)^{-1}\}$. Then
$y^gx^g \in \{y'x, (y'x)^{-1}x^{-1}\}$ and so $x^g \in \{x, x^{-1}\}$.
As $|G|$ is odd, $x^{-1} \notin x^G$. Hence $x^g = x$, and thus
$g \in C_G(x)$. As $C_G(x)$ normalizes $V = C_P(x)$, $C_G(x)$
transitively permutes the set $Vx \smallsetminus \{x\}$.
But $|Vx \smallsetminus \{x\}| = |V| - 1$ is even, whereas
$|C_G(x)|$ is odd, contradicting Lagrange's Theorem. Thus
$a_G > b_{G/H}$.
\end{proof}

\begin{lemma}
Let $G$ be a finite group of odd order, and let $H \trianglelefteq G$.
Suppose that $a_G = b_{G/H}$. Then $F(H)$ is a $p$-group for some prime $p$,
and the Sylow $q$-subgroups of $H$ are cyclic for all primes $q \neq p$.
\end{lemma}

\begin{proof}
The result is trivial if $|H| = 1$. Otherwise let $p$ be a prime divisor
of $|F(H)|$ and let $P$ be a nontrivial abelian normal subgroup of $H$.
By Lemma~1.4, $C_P(x) = 1$ for all elements $x$ of $H$ of prime order $q$
with $q \neq p$. Therefore $F(H)$ is a $p$-group. Moreover, if $H$
contains a noncyclic abelian $q$-subgroup $A$ for some prime $q \neq p$,
then it follows from Theorem~2.3 below that $C_P(x) \neq 1$ for some
element $x \in A$ of order $q$, a contradiction. Hence, as $|G|$ is odd,
all Sylow $q$-subgroups of $H$ are cyclic for $q \neq p$, as claimed.
\end{proof}

Now, we obtain our first major result about the Laitinen number $a_G$.
First, we wish to recall that a finite group $G$ is an Oliver group
if and only if $G$ does not have subgroups $P \trianglelefteq H \trianglelefteq G$
such that $H/P$ is cyclic, $P$ is a $p$-group and $G/H$ is a $q$-group for some
primes $p$ and $q$, not necessarily distinct.

\begin{proposition}
Let $G$ be a finite Oliver group of odd order. Suppose that each cyclic
quotient of $G$ has prime power order. Then $a_G > b_{G/G^{{\rm nil}}} \geq 1$.
\end{proposition}

\begin{proof}
Set $H = G^{{\rm nil}}$. As $a_G = 0$ if and only if $b_{G/H} = 0$,
we are done once we prove that $a_G > b_{G/H}$. So, assume on the contrary
that $a_G = b_{G/H}$. Then Lemma~1.5 asserts that $F(H) = P$ is a $p$-group
for some odd prime $p$. By assumption, for any prime $q \neq p$, $G$ has
no cyclic quotient of order $pq$. Hence $G/H$ is an $r$-group for some
prime $r$, and thus $H \neq F(H)$ because $G$ is an Oliver group.
By Lemma~1.1\,(3) applied for $K = P$,
$a_{G/P} = b_{(G/P)/(H/P)}$. Thus, by Lemma~1.5 applied to $G/P$,
$F(H/P)$ is a $q$-group for some prime $q$. As $P = F(H)$, we have that
$q \neq p$. Let $F_2$ be the pre-image in $H$ of $F(H/P)$ and write
$F_2 = PQ$ with $Q$ a $q$-group. Again, by Lemma~1.5, $Q$ is cyclic.
Moreover, by Lemma~1.4, $C_P(Q) = 1$ and so $N = N_G(Q)$ is a complement
to $P$ in $G$ by the Frattini argument. As $Q$ is cyclic, ${\rm Aut}(Q)$ is
abelian and so $N/C_G(Q)$ is abelian. Hence $PC_G(Q)$ is a normal subgroup
of $G$ with abelian quotient, whence $H \leq PC_G(Q)$, since $H$ is the
smallest normal subgroup of $G$ with nilpotent quotient. But $QP/P$ is
the Fitting subgroup of $H/P$, whence $C_H(Q) \leq C_H(QP/P) = QP$.
Thus $H = QP$. But then $G$ is not an Oliver group, contrary to assumption.
\end{proof}

\begin{proof}[Proofs of Theorems~A1 and A2]
For $G$ as in Theorems~A1 and A2, we shall prove that $a_G \geq 2$ and 
$LO(G) \neq 0$.

First, assume that $G$ is a finite Oliver group of odd order.
If each cyclic quotient of $G$ has prime power order, then
$a_G > b_{G/G^{{\rm nil}}} \geq 1$ by Proposition~1.6, and thus
$IO(G,G^{{\rm nil}}) \neq 0$ by the Second Rank Lemma in Section~0.4.
In particular, $a_G \geq 2$ and $LO(G) \neq 0$ by the Subgroup Lemma
in Section~0.4.

Now, assume that $G$ is a finite (not necessarily Oliver) group with
a cyclic quotient of order $pq$ for two distinct odd primes $p$ and $q$.
We will prove that $a_G \geq 4$ and $LO(G) \neq 0$. As a result, we will
complete the proof of Theorem~A1 and show that Theorem~A2 holds.

Take $H \trianglelefteq G$ with $G/H \cong \mathbb{Z}_{pq}$ and note that
$\mathbb{Z}_{pq}$ contains $\phi(pq) = (p-1)(q-1)$ elements of order~$pq$,
and hence $\frac{1}{2}(p-1)(q-1)$ real conjugacy classes of elements of
order~$pq$. We may assume that $p \geq 3$ and $q \geq 5$, and as
$a_G \geq b_{G/H} \geq a_{G/H}$, we see that $a_G \geq 4$. We will
prove that $LO(G) \neq 0$ by constructing a nonzero element of $LO(G)$.
Set $n = pq$. Let $\zeta_n$ be the primitive $n$-th root of unity.
Assume that $H = 1$ so that
$G = \mathbb{Z}_n = \langle\, g \, \vert \, g^n = 1 \,\rangle$.
Take $U = U_1 \oplus U_2$ and $V = V_1 \oplus V_2$, where
$U_i$ and $V_i$ ($i=1,2$) are the irreducible $1$-dimensional
complex $G$-modules with characters
$$
\aligned
\chi_U(g)  &= \chi_{U_1}(g) + \chi_{U_2}(g) =
	      \zeta_n + \zeta_n^2 \\
\chi_V(g)  &= \chi_{V_1}(g) + \chi_{V_2}(g) =
	      \zeta_n^a +  \zeta_n^b,
\endaligned
$$
and the integers $a$ and $b$ are chosen so that the following holds:
$$
\aligned
a \equiv 1\pmod p, &\ \ \ a \equiv 2\pmod q \\
b \equiv 2\pmod p, &\ \ \ b \equiv 1\pmod q
\endaligned
$$
(for example, if $p =3$ and $q = 5$, $a = 7$ and $b = 11$).
Then $U$ and $V$ are complex $\mathcal{L}$-free $G$-modules isomorphic
when restricted to $P \cong \mathbb{Z}_p$ or $\mathbb{Z}_q$.
The realifications $r(U)$ and $r(V)$ are not isomorphic as real
$G$-modules (remember $p$ and $q$ are odd) but $r(U)$ and $r(V)$
are isomorphic when restricted to $P \cong \mathbb{Z}_p$ or $\mathbb{Z}_q$.
So, as a result, we obtain that $0 \neq r(U) - r(V) \in LO(G)$.
If $H \neq 1$, the epimorphism $G \to G/H \cong \mathbb{Z}_n$
(mapping large subgroups of $G$ onto large subgroups of $G/H$) allows us
to consider the complex $\mathbb{Z}_n$-modules $U$ and $V$ constructed above
as complex $\mathcal{L}$-free $G$-modules. As before, we obtain that
$0 \neq r(U) - r(V) \in LO(G)$, completing the proofs.
\end{proof}

\section{Proof of the Classification Theorem}

In this section, we wish to classify finite Oliver groups $G$ with Laitinen number
$a_G = 0$ or $1$, in such a way that we obatin a proof of the Classification Theorem
stated in Section~0.2.

\begin{theorem}
Let $G$ be a finite Oliver group. Then $a_G = 0$ or $1$ if and only if one of
the conclusions (1)--(13) in the Classification Theorem holds.
\end{theorem}

In the proof, our analysis will make repeated use of a few basic concepts
and theorems. Recall that for a finite group $H$, the Fitting subgroup $F(H)$
of $H$ is the largest normal nilpotent subgroup of $H$, $E(H)$ denotes the largest
normal semisimple subgroup of $H$, and $F^*(H) = E(H)F(H)$ is the generalized
Fitting subgroup of $H$, as defined by Helmut Bender. In the proof of Theorem~2.1,
we shall use the fundamental results of \cite[Theorems~3.5, 3.6]{Gorenstein-Lyons-Solomon:2}
describing the structure and embedding of $F^*(H)$ for a finite group $H$.

\begin{theorem}
{\rm (Fitting--Bender Theorem)}
For a finite group $H$, the following holds:
$[E(H),F(H)] = 1$ and $C_H(F^*(H)) = Z(F(H))$.
If $H$ is solvable, then $F^*(H) = F(H)$.
\end{theorem}

\begin{theorem}
If $E \cong C_p \times C_p$ acts on an abelian $q$-group $V$,
where $p$ and $q$ are distinct primes, then
$V = \langle C_V(e) : e \in E^\# \rangle$.
\end{theorem}

\begin{theorem}
For two finite groups $H$ and $K$, let $F = K \rtimes H$ be
a Frobenius group with kernel $K$ and complement $H$. If $F$
acts faithfully on a vector space $V$ over a field of characteristic $p$,
where $(p, |K|) = 1$, then $C_V(H) \neq 0$.
\end{theorem}

We shall also use Burnside's $p^aq^b$ Theorem which asserts that
a finite nonsolvable group has order which is always divisible by at
least three distinct primes, the Feit--Thompson Theorem which asserts
that finite groups of odd order are solvable; the Brauer--Suzuki Theorem
which asserts that if $G$ is a finite group with no nontrivial normal
subgroup of odd order and with a Sylow $2$-subgroup of $2$-rank $1$,
then $G$ has a unique involution $z$ lying in $Z(G)$. Finally,
we shall use the Classification of the finite simple groups
(cf.~\cite{Gorenstein-Lyons-Solomon:1}--\cite{Gorenstein-Lyons-Solomon:5}).

The proof of Theorem 2.1 will be accomplished in a sequence of lemmas,
the first two of which will address the following general situation:
finite groups $H$ without NPP elements; that is, each element of $H$
is of prime power order. Such finite groups $H$ are called {\it CP groups},
and CP groups have been studied by several authors including Higman
\cite{Higman}, Suzuki \cite{Suzuki}, Bannuscher--Tiedt \cite{Bannuscher-Tiedt},
and Delgado--Wu \cite{Delgado-Wu}.

We remark that finite simple CP groups were first classified by Michio Suzuki
in a deep paper \cite{Suzuki}, whose main theorem is one of the fundamental
results in the proof of the classification of finite simple groups.

The following lemma goes back to Higman \cite{Higman}.

\begin{lemma}
Let $H$ be a finite solvable CP group. Then one of the following conclusions holds:
\begin{itemize}
\item[{\rm (1)}] $H$ is a $p$-group for some prime $p$; or
\item[{\rm (2)}] $H = K \rtimes C$ is a Frobenius group with kernel $K$ and complement
                $C$, where $K$ is a $p$-group and $C$ is a $q$-group of $q$-rank $1$
                for two distinct primes $p$ and $q$; or
\item[{\rm (3)}] $H = K \rtimes C \rtimes A$ is a $3$-step group, in the sense that
                $K \rtimes C$ is a Frobenius group as in the conclusion~(2)
                with $C$ cyclic, and $C \rtimes A$ is a Frobenius group with
                kernel $C$ and complement $A$, a cyclic $p$-group.
\end{itemize}
\end{lemma}

\begin{corollary}
If $G$ is a finite Oliver CP group, then $F^*(G)$ is nonsolvable.
\end{corollary}

\begin{proof}
As none of the groups in the conclusions of Lemma~2.5 is an Oliver group,
the result follows by Lemma~2.5.
\end{proof}

Now, we analyze the situation where a finite nonabelian simple group $L$
is without NPP elements or all NPP elements of $L$ have the same order.

\begin{lemma}
Let $L$ be a finite nonabelian simple group. Assume that $L$ is without
NPP elements or all NPP elements of $L$ have the same order. Then $L$ is
isomorphic to one of the following groups:
\begin{itemize}
\item[{\rm (1)}] $PSL(2,q)$ with $q \equiv \pm 3 \pmod 8$; or
\item[{\rm (2)}] $PSL(2,q)$ with $q = 9$ or $q$ a Fermat or Mersenne prime; or
\item[{\rm (3)}] $PSL(2,2^n)$ or $Sz(2^n)$, $n \geq 3$; or
\item[{\rm (4)}] $PSL(3,3)$, $PSL(3,4)$, $A_7$, $M_{11}$ or $M_{22}$.
\end{itemize}
\end{lemma}

\begin{proof}
We survey the finite simple groups freely making use of
the information in \cite{Gorenstein-Lyons-Solomon:3} and \cite{Conway-Wilson}.
If $L \cong A_n$ for some $n \geq 8$,
then $L$ contains elements of orders $6$ and $15$, contrary to assumption.
By inspection of \cite[Tables~5.3]{Gorenstein-Lyons-Solomon:3}, we see that
if $L$ is a sporadic simple group, then $L \cong M_{11}$ or $M_{22}$.

Hence, we may assume that $L$ is a finite simple group of Lie type
defined over a field of characteristic $p$. Assume that $L$ is not
isomorphic to $PSL(2,q)$ or $Sz(2^n)$. Then by consideration of
subsystem subgroups (\cite[Section~2.6]{Gorenstein-Lyons-Solomon:3}),
we see that one of the following statements is true about $L$:

\begin{itemize}
\item[{\rm (1)}] $L$ has a subgroup $K$ with $K/Z(K) \cong PSL(4,p)$,
                $PSU(4,p)$, $PSp(6,p)$, $G_2(p)$ or $^2F_4(2)'$; or
\item[{\rm (2)}] $L \cong PSL(3,q)$, $PSU(3,q)$ or $PSp(4,q)$.
\end{itemize}

Suppose first that $p$ is odd. Then $PSp(4,p)$, $PSL(4,p)$, $PSU(4,p)$
and $G_2(p)$ all contain subgroups isomorphic to a commuting product of
$SL(2,p)$ and a cyclic group of order $4$ (see \cite[Table~4.5.1]{Gorenstein-Lyons-Solomon:3}).
Hence, each contains elements of order $6$ and $12$. Thus, we are reduced
to the cases $L \cong PSL(3,q)$ or $L \cong PSU(3,q)$. In both cases,
$L$ contains a subgroup isomorphic to $SL(2,q)$. If $q > 3$, then
$SL(2,q)$ contains an element of odd prime order $r > 3$ and hence
elements of orders $6$ and $2r$, contrary to assumption. Thus, we may
assume that $L \cong PSL(3,3)$ or $PSU(3,3)$. We readily check that
$PSU(3,3)$ contains elements of order $12$ (cf.~\cite{Conway-Wilson}),
completing the case when $p$ is odd.

Now suppose that $p = 2$. Now $SL(3,2^n)$ contains a subgroup isomorphic
to $GL(2,2^n)$, whence $PSL(3,2^n)$ contains $H = J \times C$,
where $J \cong SL(2,2^n)$ and $C$ is cyclic of order $2^n - 1$ or
$\frac {2^n - 1}{3}$. If $n > 2$, this contradicts Lemma~1.3.
Similarly, $SU(3,2^n)$ contains a subgroup isomorphic to $GU(2,2^n)$,
whence $PSU(3,2^n)$ contains $H_1 = J_1 \times C_1$ with
$J_1 \cong SL(2,2^n)$ and $C_1$ cyclic of order $2^n + 1$ or
$\frac{2^n + 1}{3}$. If $n > 1$, this again contradicts Lemma~1.3.
Finally, $PSp(4,2^n) = Sp(4,2^n)$ contains a subgroup isomorphic
to $GL(2,2^n)$, again giving a contradiction when $n > 1$.

We know that $PSL(4,2) \cong A_8$, $PSU(4,2) \cong PSp(4,3)$ and
$G_2(2)' \cong U_3(3)$. By inspection in \cite{Conway-Wilson}, $Sp(6,2)$
and $^2F_4(2)'$ have elements of orders $6$ and $10$. We conclude that
the only examples with $p = 2$ are $PSL(3,2) \cong PSL(2,7)$,
$PSL(3,4)$ and $Sp(4,2)' \cong A_6$.

Finally, suppose that $L \cong PSL(2,q)$ and $q \equiv \varepsilon \pmod 8$,
$\varepsilon = \pm 1$.  Then $L$ has a cyclic subgroup of order
$\frac {q - \varepsilon}{2}$. If $r$ is an odd prime divisor of
$q - \varepsilon$, then $L$ has elements of order $2r$ and $4r$,
contrary to assumption. Hence, $q$ is a Fermat or Mersenne prime,
or $q = 9$, completing the proof.
\end{proof}

\begin{lemma}
Suppose that $F^*(G) = L$ is a finite nonabelian simple group
and $a_G = b_{G/L}$.
Then $G$ is isomorphic to one of the following groups:
\begin{itemize}
\item[{\rm (1)}] $PSL(2,q)$, $q \in \{ 5,7,8,9,11,13,17 \}$; or
\item[{\rm (2)}] $Sz(8)$, $Sz(32)$, $A_7$, $PSL(3,3)$, $PSL(3,4)$, $M_{11}$
                 or $M_{22}$; or
\item[{\rm (3)}] $PGL(2,5)$, $PGL(2,7)$, $P \varSigma L(2,8)$, $M_{10}$,
                 ${\rm Aut}(A_6)$, $P \varSigma L(2,27)$ or the extension
                 $PSL(3,4)^* = PSL(3,4) \rtimes \langle u \rangle$ of $PSL(3,4)$
                 by an involutory graph-field automorphism $u$ of order $2$.
\end{itemize}
If $G$ is a CP group, then $G$ is isomorphic
to one of the following groups: $PSL(2,q)$, $q \in \{ 5,7,8,9,17 \}$;
or $Sz(8)$, $Sz(32)$, $PSL(3,4)$ or $M_{10}$. Moreover, if $G = {\rm Aut}(A_6)$
or $P \varSigma L(2,27)$, then $a_G = 2$. In all other cases, $a_G \leq 1$.
\end{lemma}

\begin{proof}
Certainly $L$ is one of the groups listed in Lemma~2.7. First,
suppose that $G \neq L$. Note that the hypotheses imply that for any
$x \in G$, all NPP elements of the coset $Lx$ have the same order.
By easy inspection (or making use of \cite{Conway-Wilson}), we see that if
$$L \in \{ PSL(3,3), PSL(3,4), A_7, M_{11}, M_{22} \},$$
then $L \cong PSL(3,4)$ and $G$ is as described.

Suppose next that $L \cong Sz(2^n)$ and let $x \in G \smallsetminus L$
be of prime order $p$. Then $x$ induces a field automorphism on $L$ and
$p$ divides $n$, whence $p$ is odd. But $C_L(x)$ has a subgroup
$H \cong Sz(2)$, which has a cyclic subgroup of order $4$.
Hence $G$ has elements of orders $2p$ and $4p$, contrary to assumption.

Suppose that $L \cong PSL(2,p^n)$. If $x \in G \smallsetminus L$
has prime order $r$ and induces a field automorphism on $L$, then $C_L(x)$
contains a subgroup $H \cong PSL(2,p)$. If $p > 3$, then $Lx$ contains
elements of orders $2r$, $3r$ and $pr$, at least two of which are not
prime powers, a contradiction. Hence $p \in \{2,3 \}$ and $Lx$ contains
elements of orders $2r$ and $3r$, whence $r \in \{ 2,3 \}$ and $C_L(x)$
is a $\{2,3 \}$-group. Hence $C_L(x) \cong PSL(2,2)$, $PSL(2,3)$ or
$PGL(2,3)$. Thus $p^n \in \{ 4,8,9,27 \}$. Now by inspection, we get the
cases listed in Lemma~2.8.

Suppose now that $L \cong PSL(2,q)$ for $q > 9$, and $G$ has no element
inducing a non-trivial field automorphism on $L$. As $L \neq G$, it
follows that $q$ is odd. Then by Lemma~2.7, $q$ is an odd power of
a prime and so $G \cong PGL(2,q)$. Let $q \equiv \varepsilon \pmod 4$,
$\varepsilon = \pm 1$. Then $G \smallsetminus L$ has an element $x$
of order $q + \varepsilon$, and two elements $y$, $y'$ in
$\langle x \rangle$ are $G$-conjugate if and only if $y' = y^{-1}$.
However $Lx$ contains $\varphi(q + \varepsilon)$ elements of order
$q + \varepsilon$, whence $\varphi(q + \varepsilon) = 2$ and
$q + \varepsilon = 6$, contrary to the assumption that $q > 9$.

As $PSL(2,4) \cong PSL(2,5)$, we conclude the following: if $L \cong PSL(2,q)$,
then $q \in \{ 5,7,8,9,27 \}$, as claimed. The precise possibilities
for $G$ as stated in the proposition may then be inferred easily
from \cite{Conway-Wilson}.

Now suppose that $G = L \cong PSL(2,q)$. First we make a numerical remark.
Suppose $2^n + 1 = 3^m$ for some natural numbers $m$ and $n$. If $m$ is odd,
then $3^m - 1 \equiv 2 \pmod 8$ and so $m = 1$. If $m = 2r$, then
$2^n = (3^r - 1)(3^r + 1)$ and so $m = 2$.

Now suppose that $G \cong PSL(2,q)$ with $q = 2^n > 8$. Then $G$ has
cyclic subgroups $D_1$ and $D_2$ of orders $2^n - 1$ and $2^n + 1$
respectively. If $n$ is odd, then $3$ divides $2^n + 1$, but $2^n + 1$
is not a power of $3$ by the first paragraph. Hence
${\rm NPP}(G) \cap D_2$ contains $\varphi(2^n + 1)$ elements of
order $2^n + 1$, lying in $\frac {\varphi(2^n + 1)}{2}$ real classes.
Thus, $\varphi(2^n + 1) \leq 2$ because $a_G \leq 1$, whence $2^n + 1 = 3$,
a contradiction. Thus $n = 2s$ is even and $3$ divides
$2^n - 1 = (2^r - 1)(2^r + 1)$. As $n > 2$, $D_1$ is not a $3$-group
and, as above, $\frac {\varphi(2^n + 1)}{2} \leq 1$. Then $2^n - 1 = 3$,
again a contradiction.

Finally suppose that $G \cong PSL(2,q)$ with $q$ odd and $q > 17$.
Then $G$ has cyclic subgroups $T_1$ and $T_{-1}$ of orders
$\frac {q-1}{2}$ and $\frac {q+1}{2}$ respectively.
If $q \equiv \pm 3 \pmod 8$ and $q \equiv \varepsilon \pmod 4$,
then $T_{\varepsilon} \cap {\rm NPP}(G)$ has
$\varphi(\frac{q - \varepsilon}{2})$ elements of order
$\frac {q - \varepsilon}{2}$ lying in
$\frac {\varphi (\frac {q - \varepsilon}{2})}{2}$ real classes.
Hence $\varphi (\frac {q - \varepsilon}{2}) = 2$, whence
$\frac {q - \varepsilon}{4} = 3$, a contradiction.

Hence by Lemma~2.7, $q$ is a Fermat or Mersenne prime.
Again assume that $q \equiv \varepsilon \pmod 4$. As $q \neq 3$,
$3$ divides $q + \varepsilon$. Suppose that $q + \varepsilon = 2 \cdot 3^m$
for some $m \geq 2$. As $q - \varepsilon = 2^k$, we have
$2^k + 2 \varepsilon = 2 \cdot 3^m$. Hence $2^{k-1} = 3^m - \varepsilon$.
If $\varepsilon = -1$, then $3^m + 1 \equiv 2$ or $4 \pmod 8$, whence
$q \leq 9$, contrary to assumption.  If $\varepsilon = 1$,
then by the first paragraph, $m \leq 2$ and $q \leq 17$,
again a contradiction. It follows that $\frac {q + \varepsilon}{2}$
is not a prime power. But then $T_{- \varepsilon} \cap {\rm NPP}(G)$
has $\varphi (\frac {q + \varepsilon}{2})$ elements of order
$\frac {q + \varepsilon}{2}$ lying in
$\frac {\varphi(\frac {q + \varepsilon}{2})}{2}$ real classes.
As usual this implies that $\varphi(\frac{q + \varepsilon}{2}) = 2$,
again a contradiction.

Finally suppose that $G = L \cong Sz(2^n)$ with $n \geq 7$ and
set $q = 2^n$. Then $G$ has cyclic subgroups $T_{\varepsilon}$
with $|T_{\varepsilon}| = q + \varepsilon \sqrt{2q} + 1$, for
$\varepsilon = \pm 1$. As $q = 2^n$, $n$ odd, we have that $5$ divides
$q^2 + 1 = |T_1||T_{-1}|$. Thus $5$ divides $|T_{\varepsilon}|$ for some
$\varepsilon$. We shall argue that $|T_{\varepsilon}|$ is not a power of
$5$ when $n \geq 7$. For suppose that it is. Let $n = 2m+1$. Then
$$q + \varepsilon \sqrt{2q} + 1 = 2^{2m+1} + \varepsilon 2^{m+1} + 1 = 5^k$$
for some $k \geq 1$. Consideration of the $2$-part of $5^k - 1$
shows that the smallest positive $k$ for which $5^k - 1$ is
divisible by $2^{m+1}$, $m \geq 1$, is $k = 2^{m-1}$.
But $2^{2m+1} + 2^{m+1} + 1 < 2^{2m+2}$, while
$5^{2^{m-1}} > 4^{2^{m-1}} = 2^{2^m}$. Thus for equality to hold,
we must have $2^m < 2m+2$, which holds only for $m \leq 2$, i.e.
for $n \leq 5$.

Thus for $G \cong Sz(2^n)$, $n \geq 7$, the cyclic subgroup
$T_{\varepsilon}$ is generated by elements in ${\rm NPP}(G)$.
Let $h = |T_{\varepsilon}|$. As $N_G(T_{\varepsilon})/T_{\varepsilon}$
has order $4$ and $a_G = 1$, we must have $\varphi(h) = 4$,
whence $h = 5$, a final contradiction.
\end{proof}

Henceforth, we assume that $G$ is a finite Oliver group. Moreover,
we assume that $G$ satisfies the following two properties:
\begin{itemize}
\item[{\rm (1)}] all elements of ${\rm NPP}(G)$ have the same order; and
\item[{\rm (2)}] if $K \trianglelefteq G$ and $K \cap {\rm NPP}(G) \neq \varnothing$,
                 then ${\rm NPP}(G) \subseteq K$.
\end{itemize}
We shall call $G$ an {\it EP group} if the properties (1) and (2) above hold.
Of course, both of these properties hold when $a_G \leq 1$. By Lemma~1.2,
the class of EP groups is closed under taking subgroups and homomorphic images.

\begin{lemma}
Suppose that $G$ is an EP group and $F(G)$ is not a $p$-group.
Then $G$ is solvable and the following conclusions hold:
\begin{itemize}
\item[{\rm (1)}] $F(G) = P \times Q$ with $P$ an elementary abelian $p$-group
                 of order $p^a$ and $Q$ an elementary abelian $q$-group of order
                 $q^b$; and either
\item[{\rm (2)}] $G/F(G)$ is an $r$-group of $r$-rank $1$, $r$ prime; with
                 $r \notin \{ p,q \}$; or
\item[{\rm (3)}] $G/F(G)$ is a nonabelian metacyclic Frobenius group of order
                 $p^c q^d$.
\end{itemize}
\end{lemma}

\begin{proof}
Clearly ${\rm NPP}(G) \subseteq Z(F(G))$, whence $F(G) = P \times Q$
with $P$ and $Q$ elementary abelian, as in (1). If $G = F(G)$,
then $G$ is not an Oliver group, contrary to assumption.
Therefore $G \neq F(G)$. Set $\overline G = G/F(G)$.
If $\overline G$ has $r$-rank greater than $1$ for some prime~$r$, then
from Theorem~2.3 it follows that $G \smallsetminus F(G)$ contains elements
of order $rs$ for $s \in \{p,q \} \smallsetminus \{ r \}$, a contradiction
to ${\rm NPP}(G) \subseteq F(G)$. So $\overline G$ has $r$-rank $1$
for every prime divisor $r$ of $|\overline G|$.

Now, suppose that $F(\overline G) = 1$. Then by the Feit--Thompson Theorem,
$\overline G$ has no nontrivial normal subgroup of odd order. Moreover
$\overline G$ has $2$-rank $1$, whence by the Brauer--Suzuki Theorem,
$1 \neq Z(\overline G) \leq F(\overline G)$, a contradiction. Thus
$F(\overline G) \neq 1$ and, as ${\rm NPP}(\overline G) = \varnothing$,
$F(\overline G) = \overline R$ is an $r$-group of $r$-rank $1$ for some
prime $r$. Moreover, note that $C_{\overline G}(Z(\overline R))$ is
a normal $r$-subgroup of $\overline G$, whence
$C_{\overline G}(Z(\overline R)) = \overline R$. Moreover, note that
$\overline G/\overline R$ is isomorphic to a cyclic $r'$-subgroup
of ${\rm Aut}(Z(\overline R))$. If $\overline G = \overline R$, then
$r \notin \{ p,q \}$ because $G$ is an Oliver group, and (2) holds.
Otherwise $\overline G = \overline R \rtimes \overline C$ with both $\overline R$
and $\overline C$ cyclic. As in Lemma 2.5, $C$ is an $s$-group for some
prime $s \neq r$. Choose $s \in \{p,q \} \smallsetminus \{ r \}$. As
$\overline R \rtimes \overline C$ is a Frobenius group,
$C_P(\overline C) \neq 1$ and so $s = p$. If also $r \neq q$,
then the same argument would yield $s = q$, a contradiction.
Hence $r = q$ and $s = p$ with $p < q$, whence (3) holds.
\end{proof}

\begin{lemma}
Suppose that $G$ is a finite Oliver group with $a_G \leq 1$ and
with $F(G)$ not a $p$-group. Then $F(G) \cong C_2^2 \times C_3$ and
one of the following holds:
\begin{itemize}
\item[{\rm (1)}] $G \cong {\rm Stab}_{A_7}(\{1,2,3 \})$; or
\item[{\rm (2)}] $G \cong C_2^2 \rtimes D_9$.
\end{itemize}
\end{lemma}

\begin{proof}
We continue the notation of Lemma~2.9, and we note that the following holds:
$|{\rm NPP}(G)| = (p^a - 1)(q^b - 1)$. Moreover, ${\rm NPP}(G)$
is a union of one or two $G$-classes of equal cardinality.

If $|\overline G| = r^c$ for some prime $r$, and $R \in {\rm Syl}_r(G)$,
then $C_R(P) = 1 = C_R(Q)$, whence
$r^c \leq \min \{ p^a - 1,q^b - 1 \}$.
On the other hand, $(p^a - 1)(q^b - 1) \leq 2r^c$, which is a contradiction.

Hence $|\overline G| = p^cq^d$ with $p < q$.
Let $F(\overline G) = \overline R$ with $|\overline R| = q^d$, and let
$R$ be the full preimage of $\overline R$ in $G$. Then, as $F(G)$ is
abelian, $1 \neq C_Q(R) \triangleleft G$. As all elements of $F(G)$
of order $pq$ lie in the same real $G$-conjugacy class, this forces
$C_Q(R) = Q$. Moreover $\overline R$ acts semi-regularly on $P^\#$,
whence $q^d$ divides $p^a - 1$.
Let $\overline S \in {\rm Syl}_p(\overline G)$.
Then $\overline S$ acts semi-regularly on $Q^\#$, whence $p^c$ divides
$q^b - 1$. Finally $(p^a - 1)(q^b - 1) = p^cq^d$ or $2p^cq^d$. If both
$p$ and $q$ are odd, then $4$ divides the left-hand side of the equation
but not the right. Hence $p = 2$. Then both $q^b - 1$ and $q^d + 1$ are
powers of $2$, whence $q = q^d = 3$ and $p^a = 4$. As $\overline G$ acts
faithfully on $P$, it follows that $\overline G \cong S_3$. In particular
$p^c = 2$, whence $q^b = 3$. Thus $F(G) \cong C_2^2 \times C_3$. Moreover,
if $R_0 \in {\rm Syl}_3(G)$, then $G = P \rtimes N_G(R_0)$ and $R_0$ is
inverted by an involution in $N_G(R_0)$. Thus either (1) or (2) holds.
\end{proof}

We keep our assumption that $G$ is a finite Oliver group $G$. Recall that
$G$ is an EP group if all elements of ${\rm NPP}(G)$ have the same order,
and the following holds: if $K \trianglelefteq G$ and
$K \cap {\rm NPP}(G) \neq \varnothing$, then ${\rm NPP}(G) \subseteq K$.

\begin{lemma}
If $G$ is an EP group, then one of the conclusions holds:
\begin{itemize}
\item[{\rm (1)}] $F^*(G)$ is a $p$-group for some prime $p$; or
\item[{\rm (2)}] $F^*(G)$ is one of the nonabelian simple groups listed in Lemma 2.7; or
\item[{\rm (3)}] $G$ is solvable and satisfies the conclusions of Lemma 2.9.
\end{itemize}
Moreover if $a_G \leq 1$, then either $F^*(G)$ is a $p$-group or one of
the conclusions of Lemma 2.8 or 2.10 holds.
\end{lemma}

\begin{proof}
Suppose that $L$ is a normal quasisimple subgroup of $F^*(G)$. Then,
by Burnside's $p^aq^b$ Theorem, there exist distinct primes $p$, $q$
and $r$ dividing $|L|$. Thus if $C_{F^*(G)}(L) \neq 1$, then
${\rm NPP}(G)$ contains elements of two distinct orders,
a contradiction. Hence $C_{F^*(G)}(L) = 1$, whence $L = F^*(G)$ is
a nonabelian simple group and one of the conclusions of Lemma 2.7 (resp. 2.8)
holds. On the other hand, if $F^*(G) = F(G)$, then either $F^*(G)$ is
a $p$-group or one of the conclusions of Lemma 2.9 (resp. 2.10) holds, as
claimed.
\end{proof}

Henceforth, we shall assume that $F^*(G) = P$ is a $p$-group.
Clearly $G \neq P$ and we set
$\overline G = G/P$. Also we let $L$ be the full pre-image in $G$
of $F^*(\overline G)$.

\begin{lemma}
Suppose that $G$ is an EP group. Then either $\overline L$ is a $q$-group
for some prime $q \neq p$ or $\overline L$ is a nonabelian simple group.
\end{lemma}

\begin{proof}
Suppose that the conclusion of Lemma 2.12 does not hold. Then, by arguing
as in Lemmas 1.1 and 1.2, we see that $\overline G$ is an EP group,
provided $\overline G$ is an Oliver group. If $L$ is nonsolvable,
then $\overline G$ is an Oliver group, whence Lemma 2.11 applied to
$\overline G$ yields that $\overline L$ is a nonabelian simple group.

Suppose that $L$ is solvable but $\overline L$ is not a $q$-group for
any prime $q$. Since all elements of ${\rm NPP}(G)$ have the same order,
$\overline L = \overline{Q} \times \overline{R}$ where $Q$ is an
elementary abelian $q$-group and $R$ is an elementary abelian $r$-group
for some primes $q$, with $r$ different from $p$. Moreover, $G$ contains no
elements of order $pq$ or $pr$, whence $|Q| = q$ and $|R| = r$.
Since $\overline L = F^*(\overline G)$ and
$\overline L$ is cyclic of order $qr$, we conclude that $G/L$ is abelian.
As ${\rm NPP}(G) \subseteq L$, in fact $G/L$ is abelian of prime power
order. But then as $P$ is a $p$-group and $L/P$ is cyclic, $G$ is not
an Oliver group, contrary to assumption.
\end{proof}

Henceforth, we assume that $G$ is a finite Oliver group with
Laitinen number $a_G \leq 1$. In the next three lemmas,
we treat the case where $\overline L$ is a $q$-group.

\begin{lemma}
If $\overline L$ is a $q$-group of $q$-rank $1$, then $q = 2$ and $G = P \rtimes K$,
where $K \cong SL(2,3)$ or $\hat{S}_4$ and $P$ is an abelian $p$-group of odd order
inverted by the unique involution of $K$.
\end{lemma}

\begin{proof}
Suppose that $\overline L$ is a cyclic $q$-group. As
$\overline L = F^*(\overline G)$, $\overline G/\overline L$
is a cyclic $q'$-group and $\overline G$ is a metacyclic
Frobenius group with kernel $\overline L$.
If $\overline x \in \overline G \smallsetminus \overline L$,
$C_P(\overline x) \neq 1$ whence $G$ has elements of order
$pr$, where $r$ is the order of $x$ and $(r,p) = 1$.
In this case $\overline G/\overline L$ has prime order $r$, and
otherwise $\overline G/\overline L$ is a $p$-group. In either case, as
$\overline L$ is cyclic, $G$ is not an Oliver group, a contradiction.

Hence $\overline L \cong Q_8$ and
$[\overline G,\overline G] \cong SL(2,3)$.
As ${\rm NPP}(\overline G) \neq \varnothing$, $P$ is inverted
by $z$ for any involution $z$ of $L$. As $P = C_G(P)$, it follows
that $G = P \rtimes K$ and $K$ contains a unique involution,
whence the lemma holds.
\end{proof}

\begin{lemma}
Suppose that $\overline L$ is a $q$-group of $q$-rank greater than $1$.
Then $\overline G$ is a solvable group without NPP elements. Moreover,
$P$ is a finite elementary abelian $p$-group and
$H = P \rtimes Q \rtimes C$, where $L = P \rtimes Q$,
$Q \in {\rm Syl}_q(G)$ and $N_G(Q) = Q \rtimes C$ is
a Frobenius group with kernel $Q$ and complement $C$
such that $C$ is a $p$-group.
\end{lemma}

\begin{proof}
Let $Z = \{ z \in Z(P) : z^p = 1 \}$. Then $Z$ is a nontrivial
elementary abelian normal $p$-subgroup of $G$. Let $E$ be an
elementary $q$-subgroup of $L$ of $q$-rank greater than $1$. Then,
according to Theorem~2.3, ${\rm NPP}(L)$ contains an element $x$
of order $pq$ with $x^q \in Z$. Hence ${\rm NPP}(G/Z) = \varnothing$
and so, applying Theorem~2.3 again in $G/Z$, we conclude that $P = Z$.
If $L = G$, then $G$ is not an Oliver group, contrary to assumption.
Thus $L \neq G$.

Let $L = P \rtimes Q$, $Q \in {\rm Syl}_q(L)$. Suppose that
$C_P(Q) = A \neq 1$. Then every element $x \in {\rm NPP}(G)$
satisfies $x^q \in A$. But then $C_P(e) = A$ for all $e \in E^\#$,
whence $P = A$ by Theorem~2.3. But then $Q \leq C_G(P) = P$,
a contradiction. Therefore $C_P(Q) = 1$ and by making use of a Frattini argument,
we see that $N_G(Q)$ is a complement to $P$ in $G$.

Let $N = N_G(Q)$. Then $N$ is without NPP elements. Suppose that $r$
is a prime divisor of $|N|$ with $r \notin \{ p,q \}$ and let $R$ be
a nontrivial $r$-subgroup of $N$. As ${\rm NPP}(N) = \varnothing$,
it follows that $Q \rtimes R$ is a Frobenius group with kernel $Q$
acting faithfully on $P$. Hence $C_P(R) \neq 1$, a contradiction.
Hence $N$ is a $\{p,q \}$-group. In particular, $N$ is solvable by
Burnside's theorem and so Lemma 2.5 applies to $N$, yielding that
either $N = Q \rtimes C$ is a Frobenius group with kernel $Q$ and
complement $C$ a $p$-group, as claimed, or $N = Q \rtimes C \rtimes A$
with $C$ a cyclic $p$-group and $A$ a cyclic $q$-group disjoint from
$Q$ (as $CA = N_N(C)$ is a complement to $Q$ in $N$).
Suppose the latter and let $y \in A$ of order $q$ and
$z \in Q \cap Z(QA)$ of order $q$. Then
$U = \langle y,z \rangle \cong C_q \times C_q$ and so
$P = \langle C_P(u) : u \in U^\# \rangle$ by Theorem~2.3.
However ${\rm NPP}(G) \subseteq L$ and $U \cap L = \langle z \rangle$,
whence $C_P(u) = 1$ for all $u \in U \smallsetminus \langle z \rangle$.
Thus $P = C_P(z)$, a~contradiction.
\end{proof}

Finally, we complete the analysis of the case when $\overline L$
is a $q$-group.

\begin{lemma}
Suppose that $G = P \rtimes Q \rtimes C$ as in Lemma 2.14.
Then one of the following conclusions hold:
\begin{itemize}
\item[{\rm (1)}] $P \cong C_3^3$ and $QC \cong A_4$; or
\item[{\rm (2)}] $P \cong C_2^4$, $PQ \cong A_4 \times A_4$ and $C \cong C_4$; or
\item[{\rm (3)}] $P \cong C_2^8$ and $QC \cong (C_3 \times C_3) \rtimes C_8$; or
\item[{\rm (4)}] $P \cong C_2^8$ and $QC \cong (C_3 \times C_3) \rtimes Q_8$.
\end{itemize}
\end{lemma}

\begin{proof}
Let $x \in Q$ of order $q$ with $C_P(x) = V$ of maximum order.
The elements of $V^\# x$ are in ${\rm NPP}(G)$.
As $(vx)^{-1} \in Vx^{-1}$, either $C_G(x)$ is transitive on $V^\# x$
or $q = 2$ and $C_G(x)$ has two equal-sized orbits on $V^\# x$.
In any case, as $QC$ is a Frobenius group, $C_G(x) = VC_Q(x)$ and
$V\langle x \rangle$ acts trivially on $V^\# x$. Hence
$|V^\#| = |V^\# x| = q^b$ for some $b \geq 0$. Let $|V| = p^a$.
Then either $p = 2$ and $b = 1$ or $q = 2$ and $p^a \in \{ p,9 \}$.
In all cases we set $Z = \{ z \in Z(Q) : z^q = 1 \}$. Thus $Z$ is
a normal elementary abelian $q$-subgroup of $QC$ and $PZ \triangleleft G$.

Suppose first that $q = 2$. As $ZC$ is a Frobenius group, $Z$ contains
a Klein $4$-subgroup $U$ and so ${\rm NPP}(G) \leq PZ \triangleleft G$.
Thus in particular $x \in Z$. Suppose further that $|V| = p$. Then
$|P| \leq p^3$. Moreover, as $P$ is a faithful $QC$-module, it follows that
$|C| \leq dim V \leq 3$, whence $p = 3 = \dim V$ and $P \cong C_3^3$.
Now, $QC$ is isomorphic to a Frobenius $\{2,3 \}$-subgroup of $SL(3,3)$,
and thus $QC \cong A_4$. Therefore (1) holds.

Next suppose that $q = 2$ and $|V| = 9$. Note that $C_Q(v) \leq Z$
for all $v \in P^\#$. In particular $C_Q(V) = A$ is an elementary
abelian $q$-group with $C_P(a) = V$ for all $a \in A^\#$, by maximal
choice of $V$. If $A \neq \langle x \rangle$, then by Theorem~2.3,
$P = V$, a contradiction. Hence $C_Q(V) = \langle x \rangle$ and
$Q/\langle x \rangle$ is isomorphic to a subgroup of $GL(V) \cong GL(2,3)$.
In particular $|Q| \leq 32$. As $C$ is fixed point free on $Q$,
$|Q| = 2^m$, $m$ even. As $|Q/Z| \leq 4$, $[Q,Q]$ is cyclic, whence
$[Q,Q] = 1$ and $Q \cong C_2^2$, $C_2^3$ or $C_4 \times C_4$.
On the other hand, $Q/\langle x \rangle$ is isomorphic to an abelian
subgroup of $GL(2,3)$ of order $8$, whence $Q/\langle x \rangle \cong C_8$,
a~contradiction.

Finally suppose that $p = 2$ and $|V| = 2^a = q+1$. As $QC$ is
a Frobenius group with $C$ a $2$-group, $Q$ is abelian.
Note that $C$ permutes the set
$$\mathcal{Z} = \{ z \in Z^\# : C_P(z) \neq 1 \}$$
in one or two equal-sized orbits. As $z \in \mathcal{Z}$ if and only if
$\langle z \rangle^\# \leq \mathcal{Z}$, it follows that $|\mathcal{Z}| = k(q-1)$,
$k \geq 1$. But also, as $|C|$ is a power of $2$, $|\mathcal{Z}| = 2^c$
for some $c \geq 1$. Hence $q = 2^d + 1 = 2^a - 1$, whence $q = 3$
and $|V| = 4$. Now $P$ is a completely reducible $Z$-module and as
before $C_P(Z) = 1$, whence $P = P_1 \oplus \dots \oplus P_r$ with
$P_i$ an irreducible $Z$-module and $|P_i| = 4$ for all $i$. Then
$C_Z(P_1 \oplus P_2) = 1$, whence $Q$ acts faithfully on
$P_1 \oplus P_2 \cong C_2^4$. Hence $Q = Z \cong C_3 \times C_3$.
As $QC$ is a Frobenius group, either $C \cong Q_8$ or $C$ is cyclic
with $|C| \leq 8$. In any case the involution of $C$ inverts $Q$.
By Theorem~2.3, at least two cyclic subgroups of $Q$ have non-trivial
fixed points on $P$ and as $a_G = 1$, $C$ permutes these subgroups
transitively. Hence $|C| \geq 4$. If $|C| = 4$, then only two cyclic
subgroups of $Q$ have non-trivial fixed points on $P$ and so $|P| = 16$
and  $L \cong A_4 \times A_4$. Therefore (2) holds.

If $|C| = 8$, then $\dim V \geq 8$. On the other hand, $\dim V \leq 8$
as $Q$ has only four cyclic subgroups of order $3$. Therefore, equality
holds and $G$ is as described either in (3) or (4),
completing the proof.
\end{proof}

We have now completed the analysis of the case where $\overline L$
is a $q$-group. Thus, for the remainder of the analysis in the proof
of Theorem~2.1, we may assume that $\overline L$ is a nonabelian simple group.

\begin{lemma}
If the $p$-group $P$ is of odd order, then $P$ is elementary
abelian and one of the following conclusions holds:
\begin{itemize}
\item[{\rm (1)}] $p = 3$ and $\overline G \cong PSL(2,q)$, $q \in \{ 5,7,9,17 \}$, or
                 $\overline G \cong M_{10}$; or
\item[{\rm (2)}] $p = 7$ and $\overline G \cong SL(2,8)$ or $Sz(8)$; or
\item[{\rm (3)}] $p = 31$ and $\overline G \cong Sz(32)$.
\end{itemize}
\end{lemma}

\begin{proof}
Let $Z = \{ z \in Z(P) : z^p = 1 \}$. Then by Theorem~2.3,
$G$ contains an element $x$ of order $2p$ with $x^2 \in Z$ and so
every element of ${\rm NPP}(G)$ has this property. In particular
$P = Z$ is elementary abelian, as usual. Moreover
${\rm NPP}(\overline G) = \varnothing$ and so
$\overline G \cong PSL(2,q)$, $q \in \{ 5,7,8,9,17 \}$,
$M_{10}$, $Sz(8)$, $Sz(32)$ or $PSL(3,4)$. If $\overline G$ contains
a subgroup isomorphic to $C_3 \times C_3$ or to $A_4$, then $G$ contains
elements of order $3p$ and so $p = 3$.

Thus if $p > 3$, then $G \cong SL(2,8)$, $Sz(8)$ or $Sz(32)$. Moreover,
consideration of the Borel subgroups of $\overline G$ in these cases
shows that $p = 7$, $7$ or $31$, respectively, as claimed.

Suppose finally that $p = 3$. If $\overline G \cong SL(2,8)$, $Sz(8)$,
$Sz(32)$ or $PSL(3,4)$, then $\overline L$ contains a Frobenius group
with kernel a $2$-group and complement of order $7$, $7$, $31$ or $5$,
respectively. But then $G$ would contain elements of order $21$, $21$,
$93$ or $15$, respectively, a~contradiction.
\end{proof}

\begin{lemma}
$F^*(G)$ is a $2$-group.
\end{lemma}

\begin{proof}
Suppose first that $p > 3$. Let $U \in {\rm Syl}_2(G)$.
Then by consideration of the Frobenius group
$\overline B = N_{\overline G}(\overline U)$, we have $\dim(P) \geq p$.
Hence if $V = C_P(z)$ for $z \in Z(U)^\#$, then
$\dim(V) \geq \frac {1}{3}p > 2$. On the other hand
$C_G(z)/V\langle z \rangle$ permutes $V^\#$ in at most two equal orbits.
Hence $p^3 - 1 \leq |U|$, which is false in all cases.

Thus we may assume that $p = 3$ and $G \cong PSL(2,5)$, $PSL(2,7)$,
$PSL(2,9)$, $PSL(2,17)$ or $M_{10}$ and with $\dim(P) \geq 4$, $6$,
$4$, $16$, $8$, respectively. Again let $U \in {\rm Syl}_2(G)$,
$z \in Z(U)^\#$ and $V = C_P(z)$. Then $C_G(z) = VU$. Thus, as above,
if $\dim(V) = d$, then $3^d - 1$ is a power of $2$, whence $d \leq 2$
and so $\dim(P) \leq 3d \leq 6$. Hence $\dim(V) = 2$, $\dim(P) \leq 6$
and $|U| \geq 8$. Thus $\overline G \cong PSL(2,7)$ or $PSL(2,9)$.
However, in both cases, $U/\langle z \rangle \cong C_2 \times C_2$,
which cannot act semiregularly on $V^\#$ by Theorem~2.3,
a~contradiction.
\end{proof}

\begin{lemma}
One of the following conclusions holds:
\begin{itemize}
\item[{\rm (1)}] $\overline G \cong PSL(2,q)$, $q \in \{ 5,7,8,9,17 \}$; or
\item[{\rm (2)}] $\overline G \cong Sz(8)$, $Sz(32)$, $PSL(3,4)$, $PGL(2,5)$
                 or $M_{10}$.
\end{itemize}
\end{lemma}

\begin{proof}
If ${\rm NPP}(\overline G) = \varnothing$, $\overline G$ is
listed above. So, assume that ${\rm NPP}(\overline G) \neq \varnothing$.
Then $G$ has no element $x$ of odd order, such that $C_P(x) \neq 1$.
In particular, it follows from Theorem~2.3 that $G$ has a cyclic Sylow
$3$-subgroup, whence $\overline G \cong PGL(2,5)$, $PGL(2,7)$, $PSL(2,11)$
or $PSL(2,13)$. Note that in the last three cases, $\overline G$
contains a Frobenius subgroup of order $21$, $55$ or $39$, respectively,
whence $G$ contains elements $x$ of order $3$, $5$, $3$, respectively,
with $C_P(x) \neq 1$, which is a contradiction. Hence $\overline G \cong PGL(2,5)$.
\end{proof}

\begin{lemma}
If $\overline G$ is one of the groups $PSL(2,7)$, $PSL(2,9)$, $PSL(2,17)$,
$PSL(3,4)$, or $M_{10}$, then one of the following conclusions holds:
\begin{itemize}
\item[{\rm (1)}] $P \cong C_2^3$ and $\overline G \cong GL(3,2)$; or
\item[{\rm (2)}] $P \cong C_2^4$ and $\overline G \cong A_6 \cong Sp(4,2)'$; or
\item[{\rm (3)}] $P \cong C_2^8$ and $\overline G \cong M_{10}$.
\end{itemize}
\end{lemma}

\begin{proof}
Suppose that $\overline G \cong PSL(2,9)$, $M_{10}$ or
$PSL(3,4)$. Let $T \in {\rm Syl}_3(G)$. Then $T \cong C_3 \times C_3$
and $N_{\overline G}(\overline T) = \overline T \rtimes \overline Q$
with $\overline Q \cong C_4$, $Q_8$, $Q_8$, respectively, and with
$\overline T \rtimes \overline Q$ a Frobenius group. Hence
$\dim(P) \geq 4$, $8$, $8$, respectively. On the other hand, if
$x \in T^\#$ and $V = C_P(x)$, then $C_G(x) = VT$ acts transitively
on $V^\#$, whence $\dim(V) \leq 2$. As $\overline Q$ transitively
permutes the set $\mathcal{T}$ of non-identity cyclic subgroups
$\langle y \rangle$ of $T$ with $C_P(y) \neq 1$, we have that
$|\mathcal{T}| = 2$, $4$,$4$, respectively. Hence $\dim(P) \leq 4$, $8$, $8$,
respectively, whence equality holds in all cases. But then if
$\overline G \cong PSL(3,4)$ and $g \in G$ of order $7$, then
$C_G(y) \neq 1$, whence $G$ contains elements of order $6$ and $14$,
a contradiction.

Next suppose that $\overline G \cong PSL(2,7)$. Let $x \in G$ be
an element of order $3$. Then $\langle x \rangle$ is a Frobenius
complement in a subgroup $F$ of order $21$, whence $C_P(x) \neq 1$.
Thus $G$ contains elements of order $6$ and therefore $G$ contains
no elements of order $14$. So $P$ is a sum of faithful $F$-module,
hence a sum of free $\langle x \rangle$-modules. On the other hand,
as $C_G(x) = C_P(x)\langle x \rangle$, we must have $|C_P(x)| = 2$.
Hence $P$ is a single free $\langle x \rangle$-module, i.e. $|P| = 8$.
Finally suppose that $\overline G \cong PSL(2,17)$. By inspection of
the $2$-modular character table for $\overline G$, we see that if
$x \in G$ of order $3$, then $\dim(C_P(x)) \geq 3$. But
$C_G(x) = C_P(x)X$ with $|X| = 9$, whence $C_G(x)$ is not
transitive on $C_P(x)^\#$, a contradiction.
\end{proof}

\begin{completion} {\rm
Now, we complete the proof of Theorem~2.1 as follows.
The possibilities for $\overline G$ listed in Lemma~2.18 and
not discussed in Lemma~2.19 are precisely those groups which
are listed in the final conclusion (13) of the Classification Theorem.
For each of these cases, if $x \in G$ is of odd order, then
$\widetilde C = C_G(x)/\langle C_V(x), x \rangle$
must transitively permute $C_P(x)^\#$. However $|\widetilde C| \leq 2$,
except in the cases when $\overline G = SL(2,8)$ or $Sz(32)$ and
both $x$ and $\widetilde C$ have order $p = 3$ or~$5$, respectively.
Consideration of the $2$-modular representations of these two groups
shows that if $W$ is an nontrivial irreducible $2$-modular representation
of $\overline G$ with $C_W(x) \neq 0$, then $|C_W(x)| > p + 1$, and so
$\widetilde C$ cannot act transitively on $C_P(x)^\#$. Thus in all of
these cases, we must have $C_V(x) = 0$ for all $x \in G$ of odd order.
Let $H = G^2$. Thus $H = G$ in all cases, except when
$\overline G = \varSigma L(2,4) \cong S_5$. Then by the above remarks,
$H$ is a CP group and so by the theorem of Bannuscher--Tiedt
\cite{Bannuscher-Tiedt}, the structure of $H$ is as specified in
the Classification Theorem, completing the proof of Theorem~2.1.\qed}
\end{completion}

\section{Proof of Theorem~A3}

The main goal of this section is to prove the following proposition
which contains our next major result about the Laitinen number $a_G$
(cf.~Proposition~1.6).

\begin{proposition}
Let $G$ be a finite nonsolvable group. If $a_G = b_{G/G^{{\rm sol}}}$,
then either $a_G \leq 1$ or $a_G = 2$ and $G \cong {\rm Aut}(A_6)$ or
$P \varSigma L(2,27)$.
\end{proposition}

By inspection in \cite{Conway-Wilson}, we see that for $G = {\rm Aut}(A_6)$,
$a_G = 2$ corresponding to elements of order $6$ and $10$, and for
$G = P \varSigma L(2,27)$, $a_G = 2$ corresponding to elements of order $6$ and $14$.

Below, we assume that $G$ is a finite nonsolvable group and we set
$H = G^{{\rm sol}}$. As we know, we always have $a_G \geq b_{G/H}$. We shall
analyze the situation where $a_G = b_{G/H}$. Clearly, in this situation each coset
$gH$ meets at most one real conjugacy class $(x)^{\pm 1}$ with $x \in {\rm NPP}(G)$.

The proof of Proposition~3.1 will proceed via a sequence of lemmas.
As the arguments are very similar to those in Section 2, we shall be
a bit sketchy. First, we remark that if $H = G$ (which amounts to saying
that $G$ is perfect), then $b_{G/H} \leq 1$ and there is nothing to prove
because $a_G = 0$ (resp., $1$) if and only if $b_{G/H} = 0$ (resp., $1$).
So, we may assume that $H < G$. Let $S$ denote the solvable radical of
$G$ (i.e., $S$ is the largest normal solvable subgroup of $G$).

\begin{lemma}
$S \leq H$ and $G/S \cong PGL(2,5)$, $PGL(2,7)$, $P \varSigma L(2,8)$,
$M_{10}$, ${\rm Aut}(A_6)$, $P \varSigma L(2,27)$ or $PSL(3,4)^*$.
\end{lemma}

\begin{proof}
Let $S_0$ be the solvable radical of $H$. Set $\overline G = G/S_0$
and note that as $G$ is nonsolvable, $\overline G$ has a subnormal
nonabelian simple subgroup $\overline L \leq \overline H$.
By Lemmas 2.5\,(3) and 2.7, we see that $C_{\overline G}(\overline L) = 1$.
Hence $\overline L = F^*(\overline G) \neq \overline G$. Then the possibilities
for $\overline G$ follow from Lemma~2.8. As $\overline S = 1$, we see that
$S_0 = S$ and the proof is complete.
\end{proof}

\begin{lemma}
Either $S = 1$ or $S$ is a $p$-group for some prime~$p$.
\end{lemma}

\begin{proof}
Suppose that $S \neq 1$. Now, by Lemma 2.9, $F(G)$ is a $p$-group
for some prime $p$. Let $\overline G = G/F(G)$ and
$\overline L = F^*(\overline G)$. Suppose that $\overline L$ is
a nonabelian simple group. As $G$ has only one nonabelian
composition factor by Lemma 3.2, $H \leq L$, where $L$ is the
pre-image of $\overline L$ in $G$. Then $S = F(G)$ and therefore
$S$ is a $p$-group, as claimed.

Now, by Lemma~2.12, we may assume that $\overline L$ is a $q$-group
for some prime $q \neq p$.
As $C_{\overline G}(\overline L) \leq \overline L$, it follows that
${\rm Aut}(\overline L)$ is nonsolvable, whence $\overline L$ has $q$-rank
at least $2$. Thus $L$ contains elements of order $pq$ and so every
NPP element of $H$ lies in $L$. Since either $p$ or $q$ is odd, it
follows that $H/L$ has $2$-rank $1$. But then by the Brauer--Suzuki
Theorem, $H/L$ contains NPP elements, whence so does $H \smallsetminus L$,
a contradiction.
\end{proof}

Now, we wish to prove that $S = 1$. In order to prove it, we assume
the contrary and argue to a contradiction. Henceforth, we set
$\overline{G} = G/S$ and $\overline{H} = H/S$ (remember $S \leq H$ by Lemma~3.2).

\begin{lemma}
$S$ is either a $2$-group or an elementary abelian $p$-group for some odd
prime~$p$, and in the latter case, every NPP element of $H$ has order $2p$.
\end{lemma}

\begin{proof}
As $\overline H$ contains a Klein $4$-group, we may apply the usual argument
to obtain the result.
\end{proof}

\begin{lemma}
$\overline H$ is not isomorphic to $PSL(2,27)$.
\end{lemma}

\begin{proof}
Suppose that the contrary claim holds: $\overline H$ is isomorphic to $PSL(2,27)$.
Then $H$ contains an element $x$ of order $14$ with $x^{14} \notin S$. Hence $H$ has
no NPP element $y$ with $y^r \in S$ for $r \in \{ 2,3 \}$. However, as
$\overline H$ has $2$-rank $2$ and $3$-rank $3$, this contradicts Theorem~2.3
for $r \neq p$.
\end{proof}

\begin{lemma}
$S$ is a $2$-group and $\overline G$ is not isomorphic to $M_{10}$.
\end{lemma}

\begin{proof}
Suppose first that $S$ is a $2$-group and $\overline G \cong M_{10}$.
Then every element of $G \smallsetminus H$ is a $2$-element, whence
$b_{G/H} \leq 1$, contrary to hypothesis.

Thus it remains to prove that $S$ is a $2$-group. Suppose $S$ is not
a $2$-group. Then $S$ is an elementary abelian $p$-group for some odd
prime $p$. Suppose that there is an involution $x \in G \smallsetminus H$.
Then by inspection $x$ centralizes a coset $Hy$ of odd order, whence $Hx$
contains NPP elements outside $Sx$. But then $Sx$ contains no NPP elements,
i.e. $x$ inverts $S$ and $H = [H,x]$ centralizes $S$, a contradiction.
Thus there is no involution in $G \smallsetminus H$, and therefore
we have the following two possibilities:
$\overline G \cong P \varSigma L(2,8)$ or $M_{10}$.

Suppose now that $\overline G \cong P \varSigma L(2,8)$. As $\overline H$
contains a Frobenius subgroup of order $56$, $S$ must be a $7$-group
and $H$ must contain elements of order $14$. Thus a $3$-element of
$H$ acts without fixed points on $S$. But then by Theorem~2.3,
some element $x \in G \smallsetminus H$ of order $3$ must have
fixed points on $S$ and so $Hx$ contains elements of orders $6$ and $21$,
a contradiction.

Finally suppose that $\overline G \cong M_{10}$. Then $\overline H$
contains an $A_4$-subgroup and threfore $S$ is a $3$-group and the
NPP elements of $H$ have order $6$. Let $t$ be an involution of $H$.
As $H$ contains only one real $G$-class of NPP elements, $C_G(t)$
permutes transitively the nonidentity elements of $C_S(t)$.
Now $|C_G(t)/C_S(t)| = 16$ and $C_S(t)$ acts trivially on itself
by conjugation. Hence $|C_S(t)| \leq 9$. Moreover, as $H$ contains
elements of order $6$, $H$ contains no elements of order $15$.
Hence, an element of $H$ of order $5$ acts fixed point freely on $S$,
whence $\dim(S)$ is a multiple of $4$. By Theorem~2.3, on the other
hand, $\dim(S) \leq 3 \dim(C_S(t))$, whence $|C_S(t)| = 9$.
Let $T$ be a Sylow $2$-subgroup of $G$ with $t \in Z(T)$.
As $t$ acts trivially on $C_S(t)$, $T/\langle t \rangle$ must
act regularly on the eight elements of $C_S(t) \smallsetminus \{ 1 \}$.
But $T/\langle t \rangle$ is a dihedral group of order $8$ and
hence has no such regular action, a final contradiction.
\end{proof}

In Lemmas~3.7--3.11 below, $x$ will be an element of $H$ of order $3$
with $U = C_S(x)$ and with $|U| = 2^a > 1$, if possible.
Set $\widetilde C = C_G(x)/ \langle U,x \rangle$.

\begin{lemma}
$U$ is elementary abelian and $\widetilde C$ transitively permutes
the set $Ux \smallsetminus \{ x \}$ of cardinality $2^a - 1$.
Moreover no chief $H$-factor in $S$ is a trivial $\overline H$-module.
\end{lemma}

\begin{proof}
Note that $|H/S|$ is divisible by at least two odd primes and
$H$ has NPP elements of order $2p$ for at most one odd prime $p$.
Therefore, no chief $H$-factor in $S$ is a trivial $\overline H$-module.

If $U = 1$, the lemma holds trivially. Suppose $U \neq 1$. As all
NPP elements of $H$ have order $6$, all elements of $U$ have order
$2$ and so $U$ is elementary abelian. Moreover all elements of
$Ux \smallsetminus \{ x \}$ are $G$-conjugate, hence $C_G(x)$-conjugate
and since $\langle U,x \rangle$ is contained in the kernel of the
conjugation action on $Ux$, the result follows.
\end{proof}

\begin{lemma}
$\overline G$ is not isomorphic to $PGL(2,5)$.
\end{lemma}

\begin{proof}
Suppose $\overline G \cong PGL(2,5)$. As $b_{G/H} = 2$, $H$ must contain
NPP elements. By using Lemma~3.7, we obtain that every chief $H$-factor
of $S$ is isomorphic to a $4$-dimensional irreducible $H/S$-module. Thus
if $y \in H$ of order $5$, we have $C_S(y) = 1$. Hence $H$ must have
elements of order $6$ and so with $x$ and $U$ as above, $U \neq 1$.
Indeed some chief $H$-factor of $S$, say $V$, is a permutation module
for $\overline H$ with $|C_V(x)| = 4$. Thus $|U| \geq 4$. But
$|\widetilde C| = 2$ and so $\widetilde C$ does not act transitively
on $Ux \smallsetminus \{ x \}$, contrary to Lemma~3.7.
\end{proof}

\begin{lemma}
$\overline H$ is not isomorphic to $PSL(2,7)$.
\end{lemma}

\begin{proof}
Suppose that $\overline H \cong PSL(2,7)$.  As $\widetilde C$ is a
$2$-group acting transitively on the involutions of $U$, $|U| \leq
2$. According to \cite[2.8.10]{Gorenstein-Lyons-Solomon:3}, the
nontrivial irreducible $GL(3,2)$-modules are the standard
$3$-dimensional module $V$, its dual $V^*$ and the Steinberg module,
which is the nontrivial constituent of $V \otimes V^*$. As a
$3$-element of $GL(3,2)$ has $1$-dimensional fixed point space on $V$
and $V^*$ and $2$-dimensional fixed point space on the Steinberg
module, it follows that $S$ has a unique irreducible composition
factor and this has dimension $3$, i.e. $S \cong V$ or $V^*$ as
$\overline H$-module. But then, as $C_G(S) = S$ and ${\rm Aut}(S)
\cong H/S$, we obtain that $G = H$, a contradiction.
\end{proof}

\begin{lemma}
$\overline H$ is not isomorphic to $PSL(2,8)$.
\end{lemma}

\begin{proof}
Suppose $\overline H \cong PSL(2,8)$. Let $y \in G \smallsetminus H$ be
an element of order $3$. Then $yS$ lies in a complement of a Frobenius
subgroup of $G/S$ of order $21$. Hence there exists $ty \in Hy$ of order
$6$ with $(ty)^3 \in S$. However there also exists $sy \in Hy$ of order
$6$ with $(sy)^3 \in H \smallsetminus S$, a contradiction.
\end{proof}

\begin{lemma}
$\overline H$ is not isomorphic to $PSL(2,9)$ or $PSL(3,4)$.
\end{lemma}

\begin{proof}
Suppose $\overline H \cong PSL(2,9)$ or $PSL(3,4)$.
Then either $\overline G \cong {\rm Aut}(A_6)$ or
$\overline H \cong PSL(3,4)$ with $|G:H| = 2$. In either case,
$|\widetilde C| = 6$. Again as $\widetilde C$ acts transitively on the
involutions of $U$, we conclude that $|U| \leq 4$. Let $E$ be
a Sylow $3$-subgroup of $H$. Then $E \cong \mathbb{Z}_3 \times \mathbb{Z}_3$
and $N_G(E)$ transitively permutes the elements of $E$ of order $3$.
Hence $|C_T(y)| \leq 4$ for all $y \in E \smallsetminus \{1\}$ and so,
by Theorem~2.3, $|S| \leq 2^8$.

Suppose that $\overline H \cong PSL(3,4)$. Then $G \smallsetminus H$
contains an element $\gamma$ such that $\overline \gamma$ has order $2$
and centralizes $N_{\overline H}(\overline E)$ and
$N_{\overline G}(\overline E) = \overline E \overline Q \times
\langle \overline \gamma \rangle$ for some quaternion group $Q$ of
order $8$ transitively permuting the nonidentity elements of $E$.
Now $\overline E \overline Q$ acts faithfully on $C_S(\gamma)$ by
Thompson's $A \times B$ Lemma (see~\cite[11.7]{Gorenstein-Lyons-Solomon:2}).
However, by Clifford Theory, a faithful $\overline E \overline Q$ module must
have dimension at least $8$. As $\dim S \leq 8$, this would force $C_S(\gamma) = S$,
which is absurd.
Hence $\overline G$ is isomorphic to ${\rm Aut}(A_6)$. Again $N_{\overline G}(\overline E)$
contains a subgroup $\overline E \overline Q$ with $\overline Q \cong Q_8$,
as above. Therefore $\dim(S) = 8$ and $C_S(E) = 1$, and $S$ is a faithful
irreducible $N_G(E)$-module. In particular, $E$ acts nontrivially
on $U = C_T(x) \cong \mathbb{Z}_2 \times \mathbb{Z}_2$ and so
$UE$ is isomorphic to $A_4 \times \mathbb{Z}_3$.
As $\overline G$ contains a subgroup isomorphic to $S_6$, by inspection
we see that $N_G(E)$ contains an involution $t$ centralizing $x$
such that $E \langle t \rangle$ is isomorphic to $S_3 \times \mathbb{Z}_3$.
Then $UE \langle t \rangle \cong S_4 \times \mathbb{Z}_3$ with
$x \in Z(UE \langle t \rangle)$. But then the coset $Ht$ contains
elements of order $6$ and $12$, whence $a_G > b_{G/H}$,
contrary to assumption.
\end{proof}

\begin{completion} {\rm
In the case where $H = G^{{\rm sol}} < G$, we have studied $G/S$, where
$S$ is the solvable radical of $G$. Having exhausted all possible structures
for $G/S$, we conclude that $S = 1$ and Proposition~3.1 may be readily verified.
In fact, as $S = 1$, the possibilities for $G$ are enumerated in Lemma~3.2.
In the case where $H \cong PSL(2,5)$, $PSL(2,7)$, $PSL(2,8)$, $PSL(2,9)$ or $PSL(3,4)$,
every element of $H$ has prime power order. Hence, $b_{G/H}\leq 1$ for every
$G$ in Lemma~3.2, unless $G \cong {\rm Aut}(A_6)$ or $P \varSigma L(2,27)$.
As the two exceptional groups are covered by the comments following the statement
of Proposition~3.1, we have completed the proof of Proposition~3.1.\qed}
\end{completion}

Now, by using the Laitinen number $a_G$, we are able to determine completely
the cases where $IO(G, G^{{\rm sol}}) \neq 0$ for finite nonsolvable groups $G$.

\begin{corollary}
Let $G$ be a finite nonsolvable group. Then
\begin{itemize}
\item[{\rm (1)}] $IO(G,G^{{\rm sol}}) = 0$ for $a_G \leq 1$,
\item[{\rm (2)}] $IO(G,G^{{\rm sol}}) \neq 0$ for $a_G \geq 2$, except when
                 $G \cong {\rm Aut}(A_6)$ or $P \varSigma L(2,27)$,
\item[{\rm (3)}] $IO(G,G^{{\rm sol}}) = 0$ and $a_G = 2$ when
                 $G \cong {\rm Aut}(A_6)$ or $P \varSigma L(2,27)$.
\end{itemize}
\end{corollary}

\begin{proof}
Set $H = G^{{\rm sol}}$. By the Second Rank Lemma in Section~0.4,
we know that ${\rm rk} \,IO(G,H) = a_G - b_{G/H}$.
If $a_G \leq 1$, then $a_G = b_{G/H}$, and thus $IO(G,H) = 0$. In turn,
if $a_G \geq 2$, then except when $G \cong {\rm Aut}(A_6)$ or
$P \varSigma L(2,27)$, $a_G > b_{G/H}$ by Proposition~3.1, and thus
$IO(G,H) \neq 0$. In the exceptional cases, we know that $a_G = b_{G/H} = 2$,
and thus $IO(G,H) = 0$.
\end{proof}

\begin{proof}[Proof of Theorem~A3]
Let $G$ be a finite nonsolvable group. We shall prove that $LO(G) = 0$
for $a_G \leq 1$, and $LO(G) \neq 0$ for $a_G \geq 2$, except when
$G \cong {\rm Aut}(A_6)$ or $P \varSigma L(2,27)$, and in the exceptional
cases, we shall prove that $LO(G) = 0$ (we already know that $a_G = 2$).

By the Subgroup Lemma in Section~0.4, the following holds:
$$IO(G, G^{{\rm sol}}) \subseteq LO(G) \subseteq IO(G,O^p(G)) \subseteq IO(G,G)$$
for any prime $p$. If $a_G \leq 1$, then $IO(G,G) = 0$ by the
First Rank Lemma in Section~0.1, and thus $LO(G) = 0$. If $a_G \geq 2$,
then except when $G \cong {\rm Aut}(A_6)$ or $P \varSigma L(2,27)$,
Corollary~3.13 asserts that $IO(G,G^{{\rm sol}}) \neq 0$, and thus $LO(G) \neq 0$.

For $G = {\rm Aut}(A_6)$, $O^2(G) = A_6 = G^{{\rm sol}}$
(and $O^p(G) = G$ for any prime $p \neq 2$). Hence $IO(G,O^2(G)) = 0$
by Corollary~3.13, and thus $LO(G) = 0$.

For $G = P \varSigma L(2,27)$, $O^3(G) = PSL(2,27) = G^{{\rm sol}}$
(and $O^p(G) = G$ for any prime $p \neq 3$). Hence
$IO(G,O^3(G)) = 0$ by Corollary~3.13, and thus again $LO(G) = 0$,
completing the proof.
\end{proof}

\section{Proof of the Realization Theorem}

In this section, we shall prove the Realization Theorem stated in Section~0.2;
i.e., we shall prove that $LO(G) \subseteq LSm(G)$ for any finite Oliver gap
group $G$. The proof follows from a number of results which we collect below.
The key results are obtained in Theorems~4.3 and 4.6.

Let $G$ be a finite group. Following \cite{Laitinen-Morimoto},
consider the real $G$-module
$$V(G) = (\mathbb{R}[G] - \mathbb{R}) - \bigoplus_{p\,|\,|G|}
	 (\mathbb{R}[G]^{O^p(G)} - \mathbb{R})$$
where $\mathbb{R}[G]$ denotes the real regular $G$-module, $\mathbb{R}[G]^{O^p(G)}$
has the canonical action of $G$, and $G$ acts trivially the subtracted
summands $\mathbb{R}$. The family of the isotropy subgroups in
$V(G) \smallsetminus \{0\}$ consists of subgroups $H$ of $G$ such that
$H$ is not large in $G$; i.e., $H \notin \mathcal{L}(G)$ (cf.~\cite{Laitinen-Morimoto}).
In particular, $V(G)$ is $\mathcal{L}$-free.

By arguing as in \cite[the proof of Theorem~0.3]{Morimoto-Pawalowski:2}
in the case $G$ is an Oliver group, we obtain the following theorem which
allows us to construct Oliver equivalent real $\mathcal{L}$-free $G$-modules
(cf.~\cite[Theorem~0.4]{Oliver:3}).

\begin{theorem}  {\rm (cf.~\cite{Morimoto-Pawalowski:2})}
Let $G$ be a finite Oliver group. Let $V_1,\dots, V_k$ be real
$\mathcal{L}$-free $G$-modules all of dimension $d \geq 0$, such that
$V_i - V_j \in IO(G)$ for all $1 \le i, j \le k$.
Set $n = d + \ell\,\dim V(G)$ for an integer $\ell$.
If $\ell$ is sufficiently large,
there exists a smooth action of $G$ on the $n$-disk $D$ with
$D^G = \{x_1, \dots, x_k\}$ and $T_{x_i}(D) \cong V_i \oplus \ell\,V(G)$
for all $1 \leq i \leq k$.
\end{theorem}

By using equivariant surgery developed in \cite{Bak-Morimoto:1},
\cite{Bak-Morimoto:2}, \cite{Laitinen-Morimoto-Pawalowski},
\cite{Laitinen-Morimoto}, \cite{Morimoto:1}--\cite{Morimoto:3},
so called ``deleting--inserting''
theorems are obtained in \cite[Theorem~2.2]{Laitinen-Morimoto-Pawalowski}
for any finite nonsolvable group $G$, and in \cite[Theorems~0.1 and 4.1]{Morimoto:3}
for any finite Oliver group $G$. Under suitable conditions, these theorems allow us
to modify a given smooth action of $G$ on a sphere $S$ (resp., disk $D$) with
fixed point set $F$, in such a way that the resulting smooth action of $G$ on $S$
(resp., $D$) has a fixed point set obtained from $F$ by deleting or inserting
a number of connected components of $F$. We restate only the ``deleting part''
of \cite[Theorem~0.1]{Morimoto:3} in a modified form presented in
\cite[Theorem~18]{Morimoto-Pawalowski:3}, where the $G$-orientation
condition of \cite{Morimoto:3} is replaced by the weaker
$\mathcal{P}$-orientation condition of \cite{Morimoto-Pawalowski:3}.

Let $G$ be a finite group. Then a real $G$-module $V$ is called
{\it $G$-oriented} if $V^H$ is oriented for each $H \leq G$, and
the transformation $g: V^H \to V^H$ is orientation preserving for each
$g \in N_G(H)$. More generally, a real $G$-module $V$ is called {\it $\mathcal{P}$-oriented}
if $V^P$ is oriented for each  $P \in \mathcal{P}(G)$, and also the transformation
$g: V^P \to V^P$ is orientation preserving for each $g \in N_G(P)$.

For example, the realification $r(U)$ of a complex $G$-module $U$ is $G$-oriented.
If $V$ is a real $G$-module, then the $G$-module $2V = V \oplus V$ is the realification
of the complexification of $V$, and thus $2V$ is $G$-oriented.

For a smooth manifold $F$ with the trivial action of $G$, a real $G$-vector
bundle $\nu$ over $F$ is called {\it $\mathcal{L}$-free} if each fiber of
$\nu$ is $\mathcal{L}$-free (as a real $G$-module).

Let $M$ be a smooth $G$-manifold. We denote by $\mathcal{F}_{{\rm iso}}(G;M)$
the family of the isotropy subgroups $G_x$ of $G$ occurring at points $x \in M$.
For $H \leq G$, the set $M^H$ (resp., $M^{=H}$) consists of points
$x \in M$ with $G_x \geq H$ (resp., $G_x = H$). In general, $M^H$
(resp., $M^{=H}$) may have connected components of different dimensions.
Henceforth, by $\dim M^H$ (resp., $\dim M^{=H}$) we mean the maximum of
the dimensions of the connected components of $M^H$ (resp., $M^{=H}$).

We denote by $\mathcal{PC}(G)$ the family of subgroups $H$ of $G$ such that
$H/P$ is cyclic for some $P \trianglelefteq H$ with $P \in \mathcal{P}(G)$.
Clearly, $\mathcal{P}(G) \subseteq \mathcal{PC}(G)$. Moreover, if $G$ is a finite
Oliver group, then $\mathcal{PC}(G) \cap \mathcal{L}(G) = \varnothing$, and thus
$\mathcal{P}(G) \cap \mathcal{L}(G) = \varnothing$ (cf.~\cite{Laitinen-Morimoto}).
The family $\mathcal{PC}(G)$ was considered for the first time by Oliver \cite{Oliver:1},
and it was denoted by $\mathcal{G}^1(G)$.

Now, we state an equivariant surgery result which allows us to construct
smooth actions of $G$ on spheres with prescribed fixed point sets.
The result is a special case of \cite[Theorem~18]{Morimoto-Pawalowski:3}
(cf.~\cite[Theorem~36]{Morimoto-Pawalowski:3}).

\begin{theorem} {\rm (cf.~\cite[Theorem~18]{Morimoto-Pawalowski:3})}
Let $G$ be a finite Oliver group acting smoothly on a homotopy sphere
$\varSigma$. Let $F$ be a union of connected components of the fixed
point set $\varSigma^G$. Suppose that the following five conditions hold.
\begin{itemize}
\item[{\rm (1)}] $\dim \varSigma^P > 2 \dim \varSigma^H$ for all subgroups
                 $P < H \leq G$ with $P \in \mathcal{P}(G)$.
\item[{\rm (2)}] $\dim \varSigma^P \ge 5$ and $\dim \varSigma^{=H} \geq 2$
                 for any $P \in \mathcal{P}(G)$ and $H \in \mathcal{PC}(G)$.
\item[{\rm (3)}] $\varSigma^P$ is simply connected for any $P \in \mathcal{P}(G)$.
\item[{\rm (4)}] The tangent $G$-module $T_x(\varSigma)$ is $\mathcal{P}$-oriented
                 for some $x \in F$.
\item[{\rm (5)}] The equivariant normal bundle $\nu_{F \subset \varSigma}$
                 is $\mathcal{L}$-free.
\end{itemize}
Then there exists a smooth action of $G$ on the sphere $S$ of the same dimension
as $\varSigma$, and such that $S^G = F$ and $\nu_{F \subset S} \cong \nu_{F \subset \varSigma}$.
Moreover,  $\dim S^P = \dim \varSigma^P$ for each $P \in \mathcal{P}(G)$.
\end{theorem}

Let $G$ be a finite group. Then a pair $(P,H)$ of subgroups $P$ and $H$
of $G$ is called {\it proper} if $P \in \mathcal{P}(G)$ and $P < H \leq G$.
Following \cite{Morimoto-Sumi-Yanagihara}, for a real $G$-module $V$ and
a proper pair $(P,H)$ of subgroups of~$G$, we set
$$d_V(P,H) = \dim V^P - 2 \dim V^H.$$

A real $G$-module $V$ is called a {\it gap} $G$-module if $d_V(P,H) > 0$
for each proper pair $(P,H)$ of subgroups of $G$. Therefore, by
the definition of gap group recalled in Section~0.2, a finite group $G$
is a gap group if and only if $\mathcal{P}(G) \cap \mathcal{L}(G) = \varnothing$
and $G$ has a real $\mathcal{L}$-free gap $G$-module.

Now, by using Theorems~4.1 and 4.2, we obtain a result for actions on
spheres similar to that one obtained in Theorem~4.1 for actions on disks.

\begin{theorem}
Let $G$ be a finite Oliver gap group. Let $V$ be a real $\mathcal{P}$-oriented
$\mathcal{L}$-free gap $G$-module containing $V(G)$ as a direct summand.
Let $V_1, \dots, V_k$ be real $\mathcal{P}$-oriented $\mathcal{L}$-free
$G$-modules all of dimension $d \geq 0$, and such that $V_i - V_j \in IO(G)$
for all $1 \leq i, j \leq k$. Set $n = d + \ell\,\dim V$ for some integer $\ell$.
If $\ell$ is sufficiently large, there exists a smooth action of $G$
on the $n$-sphere $S$ with $S^G = \{x_1, \dots, x_k\}$ and
$T_{x_i}(S) \cong V_i \oplus \ell\,V$ for all $1 \leq i \leq k$.
\end{theorem}

\begin{proof}
Let $\mathcal{S}(G)$ be the family of all subgroups of $G$.
By \cite{Laitinen-Morimoto}, we know that
$\mathcal{F}_{{\rm iso}}(G; V(G) \smallsetminus \{0\})
				  = \mathcal{S}(G) \smallsetminus \mathcal{L}(G)$
and $\mathcal{PC}(G) \cap \mathcal{L}(G) = \varnothing$. Therefore
$$\mathcal{PC}(G) \subseteq \mathcal{F}_{{\rm iso}}(G; V(G) \smallsetminus \{0\}).$$
As $V$ contains $V(G)$ as a direct summand, $\dim V^{=H} \geq \dim V(G)^{=H} \geq 1$
for each $H \in \mathcal{PC}(G)$. Now, for each $i = 1, \dots, k$, consider the invariant
unit sphere
$$\varSigma_i = S(V_i \oplus \ell\,V \oplus \mathbb{R}),$$
where $G$ acts trivially on $\mathbb{R}$. The fixed point set $\varSigma_i^G$
consists of exactly two points, say $a_i$ and $b_i$, at which
$T_{a_i}(\varSigma_i) \cong T_{b_i}(\varSigma_i) \cong V_i \oplus \ell\,V$.
Set $F_i = \{b_i\}$. We note that
$n = d + \ell\,\dim V = \dim V_i + \ell\,\dim V = \dim \varSigma_i$.

We claim that the conditions (1)--(5) in Theorem~4.2 all hold for
the sphere $\varSigma_i$, provided $\ell$ is sufficiently large.
As $d_V(P,H) > 0$, we can choose $\ell$ so that
$$\ell\,d_V(P,H) > - d_{V_i}(P,H)$$
for each proper pair $(P,H)$ of subgroups of $G$. Then
$$d_{V_i \oplus \ell\,V}(P,H) = d_{V_i}(P,H) + \ell\,d_V(P,H) >
				    d_{V_i}(P,H) - d_{V_i}(P,H) = 0,$$
and thus $\dim \varSigma_i^P > 2 \dim \varSigma_i^H$, proving that
the condition (1) holds.

As $\dim V^{=H} \geq 1$ for each $H \in \mathcal{PC}(G)$, we see that
the following holds:
$$\dim \varSigma_i^H \geq \dim (\ell\,V)^H = \ell\,\dim V^H
  \geq \ell\,\dim V^{=H} \geq \ell$$
and similarly $\dim \varSigma_i^{=H} \geq \ell\,\dim V^{=H} \geq \ell$.
Hence, if $\ell \geq 5$, the condition (2) holds and the sphere $\varSigma_i^P$
is simply connected for each $P \in \mathcal{P}(G)$, proving that the
condition (3) also holds. As $V_i$ and $V$ are $\mathcal{P}$-oriented and
$T_{b_i}(\varSigma_i) \cong V_i \oplus \ell\,V$, the condition (4) holds.
Similarly, as $V_i$ and $V$ are $\mathcal{L}$-free and
$\nu_{F_i \subset \varSigma_i}$ has just one fiber $V_i \oplus \ell\,V$,
the condition (5) holds. As a result, the conditions (1)--(5) in Theorem~4.2
all hold, proving the claim.

Thus, we may apply Theorem~4.2 to obtain a smooth action of $G$ on a copy $S_i$
of the $n$-sphere such that
$S_i^G = F_i = \{b_i\}$ and $T_{b_i}(S_i) \cong V_i \oplus \ell\,V$,
provided $\ell$ is sufficiently large.

As $V_i - V_j \in IO(G)$ for all $1 \leq i,j \leq k$, Theorem~4.1 asserts
that there exists a smooth action of $G$ on the $n$-disk $D_0$ such that
$D_0^G = \{x_1, \dots, x_k\}$ and $T_{x_i}(D_0) \cong V_i \oplus \ell\,V$
for all $1 \leq i \leq k$, provided $\ell$ is sufficiently large.

The equivariant double $S_0 = \partial (D_0 \times [0,1])$ of $D_0$ is
a copy of the $n$-sphere equipped with a smooth action of $G$ such that
$S_0^G = \{x_1, y_1, \dots, x_k, y_k\}$ and
$$T_{x_i}(S_0) \cong T_{y_i}(S_0) \cong V_i \oplus \ell\,V \cong
  T_{b_i}(S_i)$$
for all $1 \leq i \leq k$. Now, consider the equivariant connected sum
$$S = S_0 \# S_1 \# \dots \# S_k$$
of the $n$-spheres $S_0, S_1, \dots, S_k$ formed by connecting sufficiently
small invariant disk neighborhoods of the points $y_i \in S_0$ and $b_i \in S_i$
for all $1 \leq i \leq k$. Then $S$ is the $n$-sphere with a smooth action of
$G$ such that $S^G = \{x_1, \dots, x_k\}$ and
$T_{x_i}(S) \cong V_i \oplus \ell\,V$ for all $1 \leq i \leq k$.
\end{proof}

We wish to remark that by using the methods of \cite{Morimoto:3},
\cite{Morimoto-Pawalowski:1}, \cite{Morimoto-Pawalowski:2},
\cite{Oliver:3}, and \cite{Pawalowski:2}, we can prove more general results
than that presented in Theorem~4.3. In fact, the results of
\cite[Theorems~27 and 28]{Morimoto-Pawalowski:3} show that each
isolated fixed point in Theorem~4.3 can be replaced by a smooth manifold
which is simply connected or stably parallelizable. However, instead of
using \cite{Morimoto-Pawalowski:3}, we decided to give an independent
proof of Theorem~4.3 due to simplifications which occur in the case
where the fixed point set is a discrete space.

Let $G$ be a finite group. A proper pair $(P,H)$ of subgroups of $G$ is called
{\it odd} if $|H:P| = |H\,O^2(G) : P\,O^2(G)| = 2$ and $P\,O^p(G) = G$ for all
odd primes $p$. Moreover, $(P,H)$ is called {\it even} if $(P,H)$ is not odd.

It follows from \cite[Theorem~2.3]{Laitinen-Morimoto} that for a proper pair
$(P,H)$ of subgroups of a finite group $G$, the following holds:
\begin{itemize}
\item[{\rm (1)}] $d_{V(G)}(P,H) = 0$ when $(P,H)$ is odd, and
\item[{\rm (2)}] $d_{V(G)}(P,H) > 0$ when $(P,H)$ is even.
\end{itemize}

Recall that by definition a real $G$-module $V$ is a gap $G$-module
if $d_V(P,H) > 0$ for each proper pair $(P,H)$ of subgroups of $G$.
If $O^p(G) \neq G$ and $O^q(G) \neq G$ for two distinct odd primes
$p$ and $q$, or $O^2(G) = G$, then any proper pair $(P,H)$ of subgroups
of $G$ is even by \cite{Morimoto-Sumi-Yanagihara}, and thus $V(G)$ is
a gap $G$-module.

In order to ensure (stably) the $\mathcal{P}$-orientability of any real
$G$-modules $V_1, \dots, V_k$ satisfying the condition that $V_i - V_j \in IO(G)$,
we use the following lemma whose proof is given at the end of this section
(cf.~\cite[Lemma~15]{Morimoto-Pawalowski:3}).

\begin{key}
Let $G$ be a finite group. Let $U$ and $V$ be two real $G$-modules
such that $U - V \in IO(G)$. Then the real $G$-module $U \oplus V$ is
$\mathcal{P}$-oriented.
\end{key}

The Key Lemma allows us to obtain the following modification of Theorem~4.3,
which we will use to prove the Realization Theorem stated in Section~0.2.

\begin{theorem}
Let $G$ be a finite Oliver gap group and let $V_1, \dots, V_k$ be real
$\mathcal{L}$-free $G$-modules with differences $V_i - V_j \in IO(G)$
for all $1 \leq i, j \leq k$. Then there exists a smooth action of $G$
on a sphere $S$ such that $S^G = \{x_1, \dots, x_k\}$ and
$T_{x_i}(S) \cong V_i \oplus W$ for all $1 \leq i \leq k$ and some
real $\mathcal{L}$-free $G$-module $W$. Moreover, $S^P$ is connected
for each $P \in \mathcal{P}(G)$.
\end{theorem}

\begin{proof}
As $G$ is a gap group, there exists a real $\mathcal{L}$-free gap
$G$-module~$U$ and so, in particular, $d_U(P,H) > 0$ for each proper
pair $(P,H)$ of subgroups of $G$. Set $V = 2U \oplus 2V(G)$.
As $d_{V(G)}(P,H) \geq 0$ by \cite[Theorem~2.3]{Laitinen-Morimoto}
$$d_V(P,H) = d_{2U \oplus 2V(G)}(P,H)
			= 2 d_U(P,H) + 2 d_{V(G)}(P,H) > 0,$$
proving that $V$ is a gap $G$-module. Clearly, $V$ is $\mathcal{P}$-oriented,
and $V$ is $\mathcal{L}$-free as so are $U$ and $V(G)$.
Moreover, $V$ contains $V(G)$ as a direct summand.

Let $V_0$ be one of the $G$-modules $V_1, \dots, V_k$. So, by assumption,
the difference $V_i - V_0$ is in $IO(G)$ for each $1 \leq i \leq k$, and thus
$V_i \oplus V_0$ is $\mathcal{P}$-oriented by the Key Lemma. Clearly, each
$G$-module $V_i \oplus V_0$ is $\mathcal{L}$-free. Again by assumption,
$(V_i \oplus V_0) - (V_j \oplus V_0) \in IO(G)$ for all $1 \leq i, j \leq k$.

Now, we may apply Theorem~4.3 to conclude that there exists a smooth action
of $G$ on a sphere $S$ such that $S^G = \{x_1, \dots, x_k\}$ and
$T_{x_i}(S) \cong V_i \oplus V_0 \oplus \ell\,V$ for all $1 \leq i \leq k$,
where $\ell$ is some sufficiently large integer. Set $W = V_0 \oplus \ell\,V$.
Then $W$ is $\mathcal{L}$-free. Moreover, $\dim W^P > 0$ for each
$P \in \mathcal{P}(G)$, as $W$ contains $V(G)$ as a direct summand.
By Smith theory, $S^P$ has $\mathbb{Z}_p$-homology of a sphere for any
$p$-subgroup $P$ of $G$. By the Slice Theorem, $\dim S^P \geq \dim W^P > 0$
and thus $S^P$ is connected for each $P \in \mathcal{P}(G)$.
\end{proof}

\begin{proof}[Proof of the Realization Theorem]
Let $G$ be a finite Oliver gap group. We shall prove that
$LO(G) \subseteq LSm(G)$. So, take an element $U - V \in LO(G)$,
the difference of two real $\mathcal{L}$-free $G$-modules $U$ and $V$
with $U - V \in IO(G)$. Then Theorem~4.4 asserts that there exists a smooth
action of $G$ on a sphere $S$ such that $S^G = \{x,y\}$ for two points
$x$ and $y$ at which $T_x(S) \cong U \oplus W$ and $T_y(S) \cong V \oplus W$
for some real $\mathcal{L}$-free $G$-module $W$, and $S^P$ is connected for
each $P \in \mathcal{P}(G)$. In particular, the action of $G$ on $S$
satisfies the $8$-condition. Consequently, the $G$-modules $U \oplus W$
and $V \oplus W$ are Laitinen--Smith equivalent, and thus
$$U - V = (U \oplus W) - (V \oplus W) \in LSm(G),$$
completing the proof.\end{proof}

In order to obtain Theorem~4.4 from Theorem~4.3 we have used the Key Lemma
asserting that given two real $G$-modules $U$ and $V$ such that $U - V \in IO(G)$,
the $G$-module $U \oplus V$ is $\mathcal{P}$-oriented, where $G$ is an arbitrary
finite group. By using some deep topological results about the existence of
specific group actions, a proof of the assertion is presented in
\cite[Lemma~15]{Morimoto-Pawalowski:3}. In the remaining part of this section,
we prove the Key Lemma using only algebraic arguments.

\begin{lemma}
Let $G$ be a finite group and let $T = \langle t \rangle$ be the cyclic
subgroup of $G$ generated by an element $t \in G$ of $2$-power order.
Let $U$ and $V$ be two real $G$-modules of the same dimension.
If $\dim U^T \equiv \dim V^T \pmod 2$, then the determinants
of the transformations $t: U \to U$ and $t: V \to V$ agree,
$\det (t|_U) = \det (t|_V)$.
\end{lemma}

\begin{proof}
If $W$ is a $2$-dimensional irreducible real $T$-module, then
the eigenvalues for $t$ on $W$ form a complex conjugate pair
and so $\det (t|_W) = 1$.

Let $m_U$ and $m_V$ be the dimensions of the $(-1)$-eigenspace
for $t$ on $U$ and $V$, respectively. Clearly, the hypothesis that
$\dim U^T \equiv \dim V^T \pmod 2$ implies that $m_U \equiv m_V \pmod 2$.
Therefore
$$\det (t|_U) = (-1)^{m_U} = (-1)^{m_V} = \det (t|_V),$$
as claimed.
\end{proof}

The next lemma is used in an inductive step of the proof of the Key Lemma.

\begin{lemma}
Let $G$ be a finite group such that $G = PT$ for some normal $p$-subgroup $P$
($p$ odd) and some cyclic $2$-subgroup $T = \langle t \rangle$. Let $U$ and $V$
be two non-zero real $G$-modules with $U^G = V^G = \{0\}$. If $U \cong V$ as
$P$-modules, then the determinants of the transformations $t: U \to U$ and
$t: V \to V$ agree, $\det (t|_U) = \det (t|_V)$.
\end{lemma}

\begin{proof}
We proceed by induction on $|P| + \dim U$. By assumption, $U \cong V$
as $P$-modules, and thus $\dim U = \dim V$. Therefore, by Lemma~4.5,
it will suffice to prove that the congruence
$$\dim U^T \equiv \dim V^T \pmod 2$$
occurring in Lemma~4.5 holds. Clearly, if $P = 1$, then $\dim U^T = \dim V^T = 0$
by hypothesis, and we are done.

Suppose now that $P \neq 1$. Let $K$ be the kernel of the $P$-action on $U$ (and $V$).
If $K \neq 1$, we are done by induction in $G/K$. Therefore we may assume that $K = 1$.
Let $E$ be a minimal normal subgroup of $G$ with $E \leq P$. Suppose that
$\dim U^E > 0$. Then $U^E \cong V^E$ and $U - U^E \cong V - V^E$ as $P$-modules
and all four of these are $G$-modules. Hence induction yields that
$\det(t|_{U^E}) = \det(t|_{V^E})$ and
$\det(t|_{U - U^E}) = \det(t|_{V - V^E})$, and we are done.

Therefore we may assume that $\dim U^E = \dim V^E = 0$. Now, if $E \neq P$,
we are done by induction in the group $ET$. As a result, we may assume that $P$
is an elementary abelian $p$-group and that $P$ is a minimal normal subgroup
of $G$. Also, $\dim U^P = \dim V^P = 0$. If $t' \in T$, then the centralizer
$C_P(t')$ is normal in $G$, hence is $1$ or $P$. As a result, either $G = P \times T$
is cyclic or the center $Z = Z(G)$ is a proper subgroup of $T$ and the quotient
$G/Z$ is a Frobenius group with kernel $PZ/Z$ and complement $T/Z$.

By \cite[Chapter VII, Theorem~1.18]{Huppert-Blackburn}, if $W$ is an irreducible
$\mathbb{R}[G]$-module, then there are the following two possibilities
for the $\mathbb{C}[G]$-module $W \otimes_{\mathbb{R}} \mathbb{C}$:
\begin{itemize}
\item[{\rm (1)}] $W \otimes_{\mathbb{R}} \mathbb{C}$ is irreducible, and we say that $W$
                 is {\it absolutely irreducible}, or
\item[{\rm (2)}] $W \otimes_{\mathbb{R}} \mathbb{C} = W_1 \oplus W_2$, where $W_1$ and $W_2$
                 are irreducible $\mathbb{C}[G]$-modules which are complex conjugate
                 (i.e., Galois conjugate).
\end{itemize}

If $G$ is cyclic, the condition on $U$ and $V$ that $U^P = V^P = \{0\}$ ensures that
$$U \otimes_{\mathbb{R}} \mathbb{C} = U_1 \oplus U_2 \ \ \hbox{and} \ \
  V \otimes_{\mathbb{R}} \mathbb{C} = V_1 \oplus V_2$$
where $U_2$ (resp., $V_2$) is the complex conjugate module of $U_1$
(resp., $V_1$).

As $\dim (U_1)^T = \dim (U_2)^T$ and
   $\dim (V_1)^T = \dim (V_2)^T$, it follows that
$$\dim U^T \equiv 0 \equiv \dim V^T \pmod 2,$$
completing the case where $G$ is cyclic. Therefore, we may assume that $G/Z$
is a Frobenius group with $|T/Z| = 2^d$ for some integer $d \geq 1$.

By Clifford theory, we know that if $W$ is an irreducible
$\mathbb{C}[G]$-module whose kernel does not contain $P$, then $\dim W$
is divisible by $2^d$. Thus in fact if $W$ is any $\mathbb{C}[G]$-module
with $W^P = \{0\}$, then $\dim W$ is divisible by $2^d$.

Suppose that $M$ is an absolutely irreducible $\mathbb{R}[G]$-module.
Then the group $Z$ maps into the group $\{I, -I\}$ of the real scalar
transformations $I$ and $-I$ of $M$. In fact, $Z$ maps into the
multiplicative group of the ring ${\rm End}_{\mathbb{R}[G]} (M) \cong \mathbb{R}$
of the endomorphisms of $M$, regarded as the ring of scalar linear
transformations acting on $M$. Since $Z$ is a $2$-group, $Z$ maps into
the group of real $2^m$th roots of $1$, which is just $\{1,-1\}$.
So we may assume that $|Z| \leq 2$.
If $|Z| = 1$, then we can replace $G$ with a larger group,
so that in fact we may assume without loss that $|Z| = 2$.
We shall argue that $Z$ acts trivially on $M$ by computing
the Frobenius--Schur indicator $\nu(\chi)$ of the character $\chi$
afforded by the absolutely irreducible $\mathbb{R}[G]$-module $M$.
By definition,
$$\nu(\chi) = \frac{1}{|G|} \sum_{g \in G} \chi(g^2).$$
Note that $\chi = {\rm Ind}_{PZ}^{\ G} (\lambda)$ for some irreducible character
$\lambda$ of $PZ$ such that ${\rm Res}^{PZ}_{\,P}(\lambda) \neq 1_P$. Since $PZ$
is a normal subgroup of $G$, thus $\chi(g) = 0$ for all
$g \in G \smallsetminus PZ$. Hence, in the displayed sum, all the terms
are $0$ except when $g^2 \in PZ$. Let $v \in T$ with $v^2 = z$. Then
$g^2 \in PZ$ if and only if $g \in P \langle v \rangle$, which is
a union of two cosets of $PZ$. Consider the squaring map on $PZ$.
This is a two-to-one map of $PZ$ onto $P$ (if $x \in P$, then
$x^2 = (xz)^2 \in P$). Since $P \langle v \rangle = PZ \cup PZ v$,
we have
$$\nu(\chi) = \frac{1}{|G|} \left(2 \sum_{g \in P} \chi(g) +
	      \sum_{g \in PZ v} \chi(g^2)\right).$$
Now
$\frac{1}{|P|} \sum_{g \in P} \chi(g) = \langle {\rm Res}^G_P(\chi), 1_P \rangle$,
the inner product of ${\rm Res}^G_P(\chi)$ and $1_P$; i.e., it is
the multiplicity of $1_P$ as a constituent of ${\rm Res}^G_P(\chi)$,
which is exactly the dimension of $M^P$, which is $0$ by assumption.
So $2 \sum_{g \in P} \chi(g) = 0$ and
$$\nu(\chi) = \frac{1}{|G|} \sum_{g \in PZ v} \chi(g^2).$$
Let $x \in P$. Then $vxv^{-1} = x^{-1}$. As $v^2 = z$,
$vxvx = vxv^{-1}v^2x = x^{-1}v^2 x = z$.
Also $vxzvxz = vxvx = z$. So $g^2 = z$ for all $g \in PZ v$. Thus
$$\nu(\chi) = \frac {|PZ | \chi(z)}{|G|} = \frac {\chi(z)}{\chi(1)}.$$
As $\chi$ is afforded by the absolutely irreducible $\mathbb{R}[G]$-module $M$,
$\nu(\chi) = 1$ and so $\chi(z) = \chi(1)$, as claimed.

Suppose now that $M$ is a sum of absolutely irreducible
$\mathbb{R}[G]$-modules such that $M^P = \{0\}$. Then $M$ may be regarded as
a faithful $\mathbb{R}[G/Z]$-module, and thus $M$ is a free $\mathbb{R}[T/Z]$-module
by the representation theory of Frobenius groups.

Now consider the decomposition $U = M_U \oplus N_U$, where $M_U$ is
the sum of all the absolutely irreducible $\mathbb{R}[G]$-summands of $U$.
Then, as $\mathbb{C}[G]$-modules,
$$N_U \otimes_{\mathbb{R}} \mathbb{C} = X_U \oplus Y_U$$
where $Y_U$ is the complex conjugate module of $X_U$, so that in particular,
we have $\dim (X_U)^T = \dim (Y_U)^T$. By the previous paragraph,
$M_U$ may be regarded as the sum of $m_U$ free $\mathbb{R}[T/Z]$-modules for
$m_U = \dim (M_U)^T$. As we know that $\dim (N_U)^T = 2 \dim (X_U)^T$,
it follows that
$$\dim U^T = \dim (M_U)^T + 2 \dim (X_U)^T \equiv m_U \pmod 2.$$
Now we may do a similar analysis for $V = M_V \oplus N_V$ and
$N_V \otimes_{\mathbb{R}} \mathbb{C} = X_V \oplus Y_V$ with obvious notations.
Therefore, it suffices to show that $m_U \equiv m_V \pmod 2$ for
$m_V = \dim (M_V)^T$. Note that
$$\dim U = 2^d \,m_U + 2 \dim X_U = 2^d \,m_V + 2 \dim X_V = \dim V.$$
By an earlier remark, both $\dim X_U$ and $\dim X_V$ are divisible by $2^d$.
So, dividing by $2^d$, we see that $m_U \equiv m_V \pmod 2$, completing the proof.
\end{proof}

\begin{proof}[Proof of the Key Lemma]
Let $G$ be a finite group. Let $U$ and $V$ be two real $G$-modules
such that $U - V \in IO(G)$. We shall prove that the $G$-module $U \oplus V$
is $\mathcal{P}$-oriented. It suffices to show that for each $P \in \mathcal{P}(G)$
and each $g \in N_G(P)$, the determinants of the transformations $g: U^P \to U^P$
and $g: V^P \to V^P$ agree,
$$\det (g|_{U^P}) = \det (g|_{V^P}),$$
because then $\det (g|_{(U \oplus V)^P}) = 1$, as required.

Let $t \in G$ be an element of $2$-power order. If $g = tx = xt$ for an
element $x \in G$ of odd order, then $\det(x) = 1$, and therefore $\det(g) = \det(t)$.
Thus it suffices to prove the claim for $g = t$. By induction on the order
of $G$, we may assume that $G = PT$ for some normal $p$-subgroup $P$ of $G$
and some cyclic $2$-subgroup $T$ of $G$. Let $t$ be a generator of $T$.

If $p = 2$, $G$ is a $2$-group and then by using the hypothesis that
$U - V \in IO(G)$, we see that $U \cong V$ as $G$-modules. Therefore,
the result is clear for $p = 2$.

Assume that $p$ is odd. As $U - V \in IO(G)$, $U \cong V$ both as
$P$-modules and $T$-modules. Write $U = U^P \oplus (U - U^P)$ and
$V = V^P \oplus (V - V^P)$, and note that $\det(t|_{U^P}) = \det(t|_{V^P})$
if and only if $\det(t|_{U - U^P}) = \det(t|_{V - V^P})$.
Since $U - U^P \cong V - V^P$ as $P$-modules, we may apply Lemma~4.6
to the $G$-modules $U - U^P$ and $V - V^P$ to conclude that
$\det(t|_{U - U^P}) = \det(t|_{V - V^P})$, and thus
$\det(t|_{U^P}) = \det(t|_{V^P})$, completing the proof.  \end{proof}

\Addresses\recd

\end{document}